\documentclass[11pt,centertags,oneside]{amsart}
\usepackage[margin=1.2 in]{geometry}
\usepackage{amsmath,amssymb,amsthm,mathtools,mathrsfs,mathdots}
\usepackage{enumitem}
\usepackage{microtype}
\usepackage{xcolor}
\usepackage{kotex}
\usepackage[
colorlinks=true,
citecolor=blue,
linkcolor=blue,
urlcolor=blue
]{hyperref}
\usepackage{graphicx}
\usepackage{float}
\usepackage{multirow, makecell}

\newcommand{\R}{\mathbb{R}}
\newcommand{\N}{\mathbb{N}}
\newcommand{\Z}{\mathbb{Z}}
\newcommand{\C}{\mathbb{C}}
\newcommand{\Q}{\mathbb{Q}}

\newcommand{\Zprim}{\mathbb{Z}_\textup{prim}}
\newcommand{\SL}{\operatorname{SL}}
\newcommand{\GL}{\mathrm{GL}}
\newcommand{\PGL}{\mathrm{PGL}}

\newcommand{\PSp}{\operatorname{PSp}}

\newcommand{\SO}{\mathrm{SO}}

\newcommand{\Sp}{\mathrm{Sp}}
\newcommand{\g}{\mathfrak g}
\newcommand{\f}{\mathfrak f}

\newcommand{\Ad}{\operatorname{Ad}}

\newcommand{\Sym}{\operatorname{Sym}}

\newcommand{\Gr}{\operatorname{Gr}}
\newcommand{\diag}{\operatorname{diag}}

\newcommand{\Rep}{\operatorname{R}^\infty}
\newcommand{\hatRep}{\widehat{\operatorname{R}}^\infty}
\newcommand{\Lie}{\operatorname{Lie}}
\newcommand{\tr}{\operatorname{Trace}}
\newcommand{\ad}{\operatorname{ad}}
\newcommand{\stab}{\operatorname{Stab}}
\newcommand{\vol}{\operatorname{Vol}}
\newcommand{\GG}{\operatorname{\mathbf{G}}}
\newcommand{\FF}{\operatorname{\mathbf{F}}}
\newcommand{\TT}{\operatorname{\mathbf{T}}}

\newcommand{\PP}{\operatorname{\mathbf{P}}}

\newcommand{\Tr}{\operatorname{Tr}}
\newcommand{\Hom}{\operatorname{Hom}}
\newcommand{\cov}{\operatorname{Cov}}

\newcommand{\V}{\boldsymbol{V}}

\newcommand{\sgn}{\operatorname{sgn}}
\newcommand{\upE}{\operatorname{E}}
\newcommand{\upH}{\operatorname{H}}
\newcommand{\upS}{\operatorname{S}}
\newcommand{\red}{\operatorname{red}}

\newcommand{\Id}{\mathrm{id}}
\newcommand{\dd}{\mathrm{d}}

\newcommand{\Uk}{\V_k^-}
\newcommand{\Vk}{\V_k^+}
\newcommand{\PiK}{\Pi_k^+}

\newcommand{\boldP}{\boldsymbol{P}}
\newcommand{\Wp}{W_{\boldP}}

\newcommand{\piL}{\pi_{\textup{left}}}
\newcommand{\piR}{\pi_{\textup{right}}}

\newcommand{\Weyl}{\mathcal W}
\newcommand{\F}{\mathcal F}
\newcommand{\CC}{\mathcal C}
\newcommand{\VV}{\mathcal V}

\newcommand{\calI}{\mathcal I}
\newcommand{\calS}{\mathcal S}
\newcommand{\calP}{\mathcal P}
\newcommand{\calF}{\mathcal F}
\newcommand{\calM}{\mathcal M}
\newcommand{\calO}{\mathcal O}
\newcommand{\calH}{\mathcal H}
\newcommand{\calD}{\mathcal D}
\newcommand{\calE}{\mathcal E}
\newcommand{\calB}{\mathcal B}

\newcommand{\bGamma}{\boldsymbol{\Gamma}}
\newcommand{\bG}{\boldsymbol{G}}
\newcommand{\bg}{\boldsymbol{g}}
\newcommand{\bw}{\boldsymbol{w}}

\newcommand{\bb}{\boldsymbol{b}}
\newcommand{\ba}{\boldsymbol{a}}
\newcommand{\bx}{\boldsymbol{x}}
\newcommand{\by}{\boldsymbol{y}}
\newcommand{\bv}{\boldsymbol{v}}
\newcommand{\be}{\boldsymbol{e}}
\newcommand{\boldf}{\boldsymbol{f}}
\newcommand{\bu}{\boldsymbol{u}}
\newcommand{\bB}{\boldsymbol{B}}
\newcommand{\bT}{\boldsymbol{T}}
\newcommand{\bU}{\boldsymbol{U}}
\newcommand{\bF}{\boldsymbol{F}}
\newcommand{\bcalO}{\boldsymbol{\calO}}
\newcommand{\bdelta}{\boldsymbol{\delta}}
\newcommand{\balpha}{\boldsymbol{\alpha}}

\newcommand{\gen}[1]{\left\langle #1 \right\rangle}
\newcommand{\set}[1]{\left\{ #1 \right\}}
\newcommand{\paren}[1]{\left( #1 \right)}
\newcommand{\lie}[1]{\mathfrak{#1}}
\newcommand{\floor}[1]{\left\lfloor{#1}\right\rfloor}
\newcommand{\bkt}[1]{[{#1}]}

\newcommand{\stdlist}{{\hypersetup{linkcolor= black}
		\hyperlink{lab:StandardRealisation}{Standard Admissible}}}

\newcommand{\de}{\delta}

\newcommand{\norm}[1]{\left\lVert #1 \right\rVert}
\numberwithin{equation}{section}
\numberwithin{figure}{section}
\numberwithin{table}{section}

\theoremstyle{plain}
\newtheorem{theorem}{Theorem}[section]
\newtheorem{lemma}[theorem]{Lemma}

\newtheorem{proposition}[theorem]{Proposition}
\newtheorem{corollary}[theorem]{Corollary}

\theoremstyle{definition}
\newtheorem{definition}[theorem]{Definition}

\theoremstyle{remark}
\newtheorem{remark}[theorem]{Remark}

\theoremstyle{remark}
\newtheorem*{remark*}{Remark}

\usepackage{fancyhdr}

\title{Irreducible $\SL_2(\R)$-Action on 
projective $k$-Lattices and Frobenius Numbers}
\author{Michael Bersudsky$^*$}
\thanks{$^*$Max Planck Institute for Mathematics in the Sciences
Geometry, Groups, and Dynamics,
Inselstrasse 22,
04103 Leipzig,
Germany. E-mail: Michael.Bersudsky@mis.mpg.de}

\author{Seongmin Kim$^{\dagger}$}
\thanks{$^{\dagger}$Department of Mathematical Sciences, Seoul National University, Seoul, 08826, Republic of Korea. E-mail: seongmin.kim@snu.ac.kr}

\author{Seonhee Lim$^{\ddagger}$}
\thanks{$^{\ddagger}$Department of Mathematical Sciences and Research Institute of Mathematics, Seoul National University, Seoul, 08826, Republic of Korea. E-mail: slim@snu.ac.kr}

\date{}

\begin{document}
	
	\begin{abstract}
The irreducible representation of \(\SL_2(\R)\) on \(\R^n\) induces an action on the space \(X_k(\R^n)\) of projective classes of  rank \(k\)-discrete subgroups in \(\R^n\). For every lattice subgroup of \(\SL_2(\R)\) and every starting point in \(X_k(\R^n)\) and \(2\leq k<n\), we prove convergence of the orbit averages over expanding norm balls and give a complete, explicit description of their limits.

In an overview, every limit is supported on the subbundle over the unique compact \(\SL_2(\R)\)-orbit in the Grassmannian. The conditional measures on the fibers are either uniform on finite sets or convex combinations of at most two homogeneous measures associated with split simple subgroups of \(\PGL_k(\R)\) containing the principal \(\SL_2(\R)\).

Another objective of our work is to relate weak-\(*\) limits of probability measures on \(X_{n-1}(\R^n)\)  supported on the projective classes of lattices orthogonal to integer
vectors with limiting distributions of the corresponding normalized
Frobenius numbers. This correspondence gives limit laws for vectors
sampled from certain proper subvarieties, extending J.~Marklof's work
on sampling from full-dimensional domains.
Specifically, together with our results for \(\SL_2(\Z)\) orbits in \(X_{n-1}(\R^n)\), this
correspondence gives limit laws for primitive integer vectors sampled
from \(\SL_2(\R)\)-orbits. Using an earlier result of the first author
and Shapira, we obtain analogous laws for integer vectors sampled on
level sets of certain quadratic forms, as the level tends to infinity.
\end{abstract}

\maketitle
    \tableofcontents
	\section{Introduction}

The moduli space \(X_n\cong \PGL_n(\R)/\PGL_n(\Z)\) of unimodular
\emph{full-rank} lattices in \(\R^n\) is an important tool in number
theory. Reformulating arithmetic problems as problems on this space makes it
possible to apply analytic, dynamical, and geometric methods. A landmark
example is Margulis's resolution of the Oppenheim conjecture
\cite{Margulis1987}.
The natural non-full-rank counterparts of \(X_n\) are the spaces
\(X_k(\R^n)\) of projective classes of rank-\(k\) discrete subgroups of
\(\R^n\), where \(1\leq k<n\); see
Section~\ref{subsec:orbit-limits} for the definition.

In this article, we demonstrate that \(X_{n-1}(\R^n)\) provides a natural
tool to study the distribution of Frobenius numbers. We recall that the Frobenius
number \(F(a_1,\dots,a_n)\), for coprime natural numbers
\(a_1,\dots,a_n\geq 2\),  is the largest number that cannot be
obtained as the sum of a collection of coins of denominations \(a_1,\dots,a_n\). For
\(k<n-1\), the diagonal Frobenius problem \cite{AlievHenk2010} suggests
a similar role for \(X_k(\R^n)\). We restrict our Frobenius applications
in this paper to the classical problem, which uses lattices of
codimension one.

Kannan \cite{Kannan1992} gave an exact formula for the Frobenius number
in terms of the covering radius of a simplex with respect to an
associated lattice, see \eqref{eq:FrobNumerFormulaCoveringRad}. This
geometric connection has been used to obtain upper and lower bounds
\cite{FukshanskyRobins2007,AlievGruber2007}. For \(n\geq 3\), Marklof
\cite{Marklof2010} proved that suitably normalized Frobenius numbers
have a limiting distribution when primitive positive integer vectors
are sampled from expanding dilates of a bounded set, in the nonnegative
orthant, of \emph{positive Lebesgue measure} whose boundary has measure
zero. He identified the limit with the distribution of the covering
radius of the standard simplex with respect to a random unimodular
lattice.

We show that Kannan's correspondence translates naturally
to the space \(X_{n-1}(\R^n)\). Proposition~\ref{prop:PushforwardTheta}
gives a general criterion for passing from limiting distributions of
the projective classes of primitive rank-\((n-1)\) sublattices of
\(\Z^n\) to limiting distributions of the corresponding normalized
Frobenius numbers. This connects tools from homogeneous dynamics
to study the distribution of these lattice classes and thereby obtain
limiting laws for Frobenius numbers of integral vectors in ``thin" sets, see Theorem~\ref{thm:frob_on_quad_var}, which uses \cite{BersudskyShapira2023}, and see Theorem~\ref{thm:1.2}, which applies the results of this work.

To our knowledge, previous work on distributions of Frobenius numbers for sampling from lower-dimensional subsets has treated only the sampling on the boundary of an expanding box. In this case, Ustinov \cite{Ustinov2010} proved a limit law with rate of convergence in this setting for three variables, and in higher dimensions, a limiting distribution follows by combining the equidistribution theorem for fixed-denominator families of multidimensional Farey points \cite[Theorem~1.3]{EinsiedlerMozesShahShapira2016} with  \cite{Marklof2010}.

The main work in this paper concerns the action of \(\SL_2(\R)\) on
\(X_k(\R^{n+1})\), for \(2\leq k\leq n\), induced by its irreducible
representation on \(\R^{n+1}\). Given a lattice
\(\Gamma\subset\SL_2(\R)\) and a starting point \(x\in X_k(\R^{n+1})\),
we determine the limiting distribution of \(\gamma x\) when
\(\gamma\) is sampled uniformly among the elements of \(\Gamma\)
in a growing norm ball in \(\SL_2(\R)\). 

This problem belongs to the general study of the distribution of
lattice orbits for continuous Lie group actions on manifolds, see e.g.
\cite{GorodnikWeiss2007}. We note that the spaces \(X_k(\R^d)\) belong to a
broader class of homogeneous bundles over flag varieties, with
fibers given by quotients of a Lie group by a discrete subgroup, and the study of the discrete subgroup 
dynamics on such spaces has received attention only recently.

The work of Sargent and Shapira \cite{SargentShapira2019} on random
walks on \(X_2(\R^3)\) motivated the study of dynamics on these
bundles. More recently, the authors of \cite{GLS24} have studied
random walks of Zariski-dense discrete subgroups of \(\SL_2(\R)\)
acting through representations of \(\SL_2(\R)\) on general homogeneous
bundles over flag varieties. In the irreducible case, their equidistribution result concerns random walks
starting in the subbundle over the unique compact orbit in the base
\cite[Proposition~4.1]{GLS24}. They explicitly leave aside the
equidistribution question for walks starting
outside this subbundle. In this work, we study norm-ball averages
for lattice actions on \(X_k(\R^{n+1})\), allowing an arbitrary
starting point. 

The problem of norm-ball averages for lattice subgroup actions on such spaces was recently treated in
\cite{BersudskyXing2024,BersudskyShah}. The work \cite{BersudskyXing2024}
studies the action of lattices in \(\SL_n(\R)\) on \(X_{n-1}(\R^n)\).
In this setting, every lattice orbit is dense, and the same analysis
applies to every starting point. The more recent work
\cite{BersudskyShah} treats the more refined setting of lattices in
\(\mathrm{SO}(n,1)^\circ\) acting on pairs of perpendicular projective
lattices in \(X_k(\R^{n+1})\times X_{n+1-k}(\R^{n+1})\). In this
setting, orbits accumulate on a subbundle over subspaces tangent to
the light cone. In dimension three, a commensurability condition
singles out exceptional starting lattices whose limiting measures
are supported on finite sets in the fibers.

The actions studied in this article overlap with those in \cite{BersudskyShah}
only in the case of \(X_2(\R^3)\), where the irreducible
\(\SL_2(\R)\)-action agrees with the \(\SO(2,1)^\circ\)-action via
the double cover.

In higher dimensions, each \(\SL_2(\R)\)-orbit in the Grassmannian
has dimension at most three, so there are uncountably many such
orbits whenever \(k(n+1-k)>3\). Moreover, there are more possibilities
for the limiting measures. Our main result shows that, for any
base point \([\Lambda]\in X_k(\R^{n+1})\), the norm-ball averages
converge to a measure concentrated on the subbundle over the unique
compact \(\SL_2(\R)\)-orbit in the Grassmannian, which is a circle.
The fiber measures either have finite support, as described in
Theorem~\ref{thm:MainTheoremGammaSpecial}, or are convex combinations
of homogeneous measures for a simple split Lie group containing
the principal \(\SL_2(\R)\) in \(\GL_k(\R)\).

Similar to \cite{BersudskyShah}, we
reduce lattice sums to integrals over norm balls. The resulting
problem is to analyze the measures on
\[
(\SL_2(\R)/\Gamma)\times(\PGL_k(\R)/\PGL_k(\Z))
\]
defined in \eqref{eq:DefJointMeasure} and supported on Wronskian
curves in the product, expanded by a diagonal one-parameter subgroup. We will use \cite{Yang2020} to analyze these weak-\(*\) limits. 

We proceed to present our results.

\subsection{Distribution of Frobenius numbers in lower dimensional varieties}

For a primitive vector \(\ba=(a_1,\dots,a_n)\in\N^n\), every sufficiently
large integer can be written as \(\bx\cdot\ba\) for some
\(\bx\in\Z_{\geq 0}^n\). We define the Frobenius number
\begin{equation}\label{eq:DefFrobNum}
    F(a_1,\dots,a_n)
    :=\max\bigl(\Z\smallsetminus
    \set{\bx\cdot\ba : \bx\in\Z_{\geq 0}^n}\bigr).
\end{equation}
As \(F(\ba)=-1\) if some \(a_i=1\), we restrict ourselves to the nontrivial case \(a_i\geq 2\) for all \(1\leq i\leq n\). 

Denote the standard simplex in \(\R^{n-1}\) by
\[
    \calS:=\set{\bx\in\R_{\geq 0}^{n-1}:
    \sum_{i=1}^{n-1}x_i\leq 1}.
\]
For a rank-$(n-1)$ lattice \(\Lambda\subset\R^{n-1}\), the covering
radius of \(\calS\) with respect to \(\Lambda\) is
\[
    \chi(\Lambda):=\inf\set{\alpha>0:
    \Lambda+\alpha\calS=\R^{n-1}}.
\]

Let \(X_{n-1}\cong\PGL_{n-1}(\R)/\PGL_{n-1}(\Z)\) be the space of homothety classes of full-rank
lattices in \(\R^{n-1}\), and let \(\mu_{X_{n-1}}\) be the \(\PGL_{n-1}(\R)\)-invariant probability
measure on \(X_{n-1}\). The normalized covering
radius is the function
\begin{equation}\label{eq:DefNormalizedCoveringRad}
    \hat\chi:X_{n-1}\to \R_{\geq 0},\quad
    \hat\chi([\Lambda])
    =\chi\bigl(\cov(\Lambda)^{-1/(n-1)}\Lambda\bigr)
    =\frac{\chi(\Lambda)}{\cov(\Lambda)^{1/(n-1)}}.
\end{equation}
 For \(\alpha\in\R\), set
\begin{equation}\label{eq:DefComulativeFrob}
    \Psi_n(\alpha):=
    \mu_{X_{n-1}}\bigl(\set{[\Lambda]\in X_{n-1}:
    \hat\chi([\Lambda])\geq\alpha}\bigr).
\end{equation}

\begin{theorem}\label{thm:frob_on_quad_var}
   Choose \(d\geq 4\) and let \(Q\) denote one of the following quadratic forms 
		\[Q_1(\boldsymbol{x}):=x_1^2+\dots+x_d^2,\;\;Q_2(\boldsymbol{x}):=-(x_1^2+\dots+x_{d-1}^2)+x_d^2.\]
		Fix a compact subset
			\(\calD\subset \set{\boldsymbol{x}\in\R^d_{>0}:Q(\boldsymbol{x})=1}\)
			with \(\mu_\calH(\calD)>0\) and \(\mu_\calH(\partial \calD)=0\), where \(\mu_\calH\) is the $\SO_Q(\R)$-invariant measure on \(\calH:=\set{\boldsymbol{x}\in\R^d:Q(\boldsymbol{x})=1}\).
			
			Let \(\set{N_i}_{i=1}^\infty\subseteq\N\) be a sequence such that \(N_i \to \infty\), and assume in addition that \(N_i \notin8\N\) for all \(i\) if \(d=4\) and \(Q=Q_1\). Then, for every \(\alpha\in \R\), 
			\[\lim_{i\to\infty}\frac{1}{\#(\sqrt{N_i}\calD\cap\Zprim^d)} \#\set{\boldsymbol{a}\in \sqrt{N_i}\calD\cap\Zprim^d:\,\frac{F(\boldsymbol{a})}{(a_1\dots a_d)^{1/(d-1)}}\geq \alpha}=\Psi_{d}(\alpha).\]
\end{theorem}
The proof of the above result,  see Section~\ref{sec:proof-of-shap-bers-applic},
allows us to derive from \cite[Theorem~4.8]{BersudskyShapira2023} more general results for rational definite quadratic
forms and for rational indefinite quadratic forms of signature
\((1,d-1)\), for vectors on double-sheeted hyperboloids. 

In Theorem~\ref{thm:frob_on_quad_var}, the level set varies with the
sampling parameter. In the next application, we sample from growing
regions of a fixed \(\SL_2(\R)\)-orbit in \(\R^{n+1}\), using elements
of \(\SL_2(\Z)\) in growing norm balls.

To motivate the result we present below, we note the following formula for Frobenius
numbers of geometric sequences obtained in
\cite{OngPonomarenko2008Frobenius}. For coprime positive integers
\(a,b\) with \(a\ne b\) and \(n\geq 1\), we have
\[
F(a^n,a^{n-1}b,a^{n-2}b^2,\dots,b^n)
=b^{n-1}(ab-a-b)
+\frac{(b-1)a^2(a^{n-1}-b^{n-1})}{a-b}.
\]
Interestingly, the vectors
\((a^n,a^{n-1}b,a^{n-2}b^2,\dots,b^n)\in\Zprim^{n+1}\) are obtained
from an \(\SL_2(\Z)\)-orbit in an irreducible representation of
\(\SL_2(\R)\) on \(\R^{n+1}\), see Remark~\ref{rem:GeomSeqAndSpecialPoints}.
Using our results, it is possible to derive the distribution of such Frobenius numbers, for sampling primitive vectors \((a,b)\in \Zprim^2\cap \R^2_{>0}\) from growing regions. However, the distribution in this case can also be derived by more elementary considerations using the above formula.

We shall apply our results to show that the limiting distributions
of normalized Frobenius numbers of sums of geometric sequences
follow the same law as in \cite{Marklof2010}. We restrict only to the even-\(n\) result for sums, since the analysis required contains the main idea and is simpler in this case.
\begin{theorem}\label{thm:1.2}
   Let \(n\) be an even natural number,  and for \(a,b,c,d\in \N\), we let
		\[\bv(a,b,c,d):=(a^n+b^n,a^{n-1}c+b^{n-1}d,\dots,c^n+d^n),\]
		and let \(R(a,b,c,d)\) be the product of the entries of \(\bv(a,b,c,d)\),
		\[R(a,b,c,d):=\prod_{k=0}^{n}
		\left(a^{n-k}c^k + b^{n-k}d^k\right).\]
		Then, for \(p \geq 1\) and \(\alpha\in\R\),
		\[\lim_{T\to\infty}\frac{\#\set{\begin{smallmatrix}
		    (a,b,c,d)\in \N^4,\;ad-bc=1,\\
            (|a|^p+|b|^p+|c|^p+|d|^p)^{1/p}\leq T
       \end{smallmatrix}:\,\frac{F(\bv(a,b,c,d))}{R(a,b,c,d)^{1/n}}\geq \alpha}
		}{c_p T^2}=\Psi_{n+1}(\alpha)\]
		with the constant $c_p=\frac{3\Gamma(1/p)^4}{\pi^2p^2\Gamma(2/p)^2},$ where \(
\Gamma(x) = \int_{0}^{\infty} t^{x-1} e^{-t}\,dt,
\,x > 0.
\)
\end{theorem}

We stress that the above sampled vectors lie
in a lower-dimensional subset of \(\R^{n+1}\) of codimension bounded below by \(n+1-3\).

\subsection{Limiting distributions of lattice orbits}
\label{subsec:orbit-limits}
For a subset \(\calS\) of a real vector space \(\V\), and  a subring \(R\subset \R\), we denote
	\[\gen{\calS}_R:=\set{\sum_{i=1}^\ell a_i v_i:\,\ell\in \N,\;a_i\in R,\,v_i\in\calS}.\]
  
	\begin{definition}\label{sec:SpaceOFkLattices}
	Let \(\VV\) be an $m$-dimensional real vector space, and let \(1\leq k\leq m\).   A \emph{$k$-lattice} \(\Lambda\) in \(\VV\) is a $\Z$-span of \(k\) $\R$-linearly independent vectors.
    
	For a \(k\)-lattice \(\Lambda\subset \VV\), we define its projective class by \([\Lambda]:=\R^\times \Lambda\)  and we define the \emph{space of projective classes of \(k\)-lattices in \(\VV\)} by \(X_k(\VV)\), that is
	\[X_k(\VV)=\set{[\Lambda]:\,\Lambda=\gen{v_1,\dots,v_k}_\Z,\;\set{v_1,\dots,v_k}\subset\VV\text{ are $\R$-linearly independent}}.\]
	\end{definition}
	Let \(\Gr_k(\VV)\) be the space of \(k\)-dimensional subspaces of \(\VV\). There is a natural map
	\[\pi:X_k(\VV)\to \Gr_k(\VV),\;\;[\Lambda]\mapsto \gen{\Lambda}_\R,\]
	which sends a class of \(k\)-lattices to the subspace containing the $k$-lattices.
	The map \(\pi\) endows \(X_k(\VV)\) with a fiber-bundle structure over \(\Gr_k(\VV)\), with fibers isomorphic to the moduli space \(X_k \simeq \PGL_k(\R)/\PGL_k(\Z)\) of classes of full-rank lattices in \(\R^k\). We will denote the fiber above \(\V\in \Gr_k(\VV)\) by
	\[X(\V):=\pi^{-1}(\V)=\set{[\Lambda]:\,\Lambda\subset \V}.\]
    
We proceed to define the irreducible \(\SL_2(\R)\)-action on such spaces. 

For \(n\in \N\), we denote the vector space of binary homogeneous polynomials of degree \(n\) by $$\VV_n := \Sym^n(\R^2) = \gen{x^n,\, x^{n-1}y,\dots
		,\,y^n}_{\R}\cong \R^{n+1},$$ 
	and we consider the irreducible \(\GL_2(\R)\)-action on \(\VV_n\), defined as 
	\begin{equation}\label{eq:DefOfIrredAction}
		g\cdot p(x,y) :=p\bigl((x,y)g\bigr),\quad g\in\GL_2(\R), \,p\in \VV_n. 
	\end{equation}
	Restricting the above representation to \(\SL_2(\R)\) gives the $(n+1)$-dimensional irreducible representation of $\SL_2(\R)$, unique up to isomorphism.
	
	   Via the above linear action, \(\SL_2(\R)\) acts on \(X_k(\VV_n)\) by
	\[g\cdot[\Lambda]:=[g\cdot \Lambda],\;g\in\SL_2(\R),\;[\Lambda]\in X_k(\VV_n).\]

For $1\le k\le n$, set
\begin{equation}\label{eq:DefDiagHighestWeights}
		\Vk:=\gen{x^n,\,x^{n-1}y,\dots,\,x^{n-(k-1)}y^{k-1} }_{\R},\quad
\CC_k:=\SL_2(\R)\cdot\Vk\subset\Gr_k(\VV_n).
	\end{equation}

Note that \(\CC_k\cong \SL_2(\R)/\begin{psmallmatrix}
    *&*\\0&*
\end{psmallmatrix}\cong \R/\pi\Z \) is the \(\SL_2(\R)\)-equivariant embedding of the projective line \(\Gr_1(\VV_1)\cong \PP^1(\R^2)\), described naturally by the map \[[ax+by]\mapsto [(ax+by)^{n-(k-1)}\Sym^{k-1}(\R^2)].\]

We parametrize \(\CC_k\) by
	\begin{equation}\label{eq:CircleIdentification}
		\sigma :\R/\pi\Z\to \CC_k,\;\;\sigma(\theta):=\kappa(\theta)\cdot \Vk, 
	\end{equation}
    where \(\kappa(\theta):=
\begin{psmallmatrix}
\cos\theta&-\sin\theta\\
\sin\theta&\cos\theta
\end{psmallmatrix}.\)

Let $\|\cdot\|$ be a norm on $M_2(\R)$, let
$\Gamma\subset\SL_2(\R)$ be a lattice, and put
\[
\Gamma_T:=\set{\gamma\in\Gamma:\|\gamma\|\le T}.
\]
Given \(\bkt\Lambda\in X_k(\VV_n)\), sufficiently large $T$ such that \(\Gamma_T\neq \emptyset\), we define the probability measures
\begin{equation}\label{eq:DefNormBallAverages}
		\mu_{T,[\Lambda]}(f):=\frac{1}{\#\Gamma_T} \sum_{\gamma\in \Gamma_T} f(\gamma\cdot[\Lambda]),\quad f\in C_c(X_k(\VV_n)),\;[\Lambda]\in X_k(\VV_n).
	\end{equation}

\begin{theorem}\label{thm:GeneralForm}
For every $2\le k\le n$ and $[\Lambda]\in X_k(\VV_n)$,
there exists a family of probability measures
$\{\mu_{\theta}^{[\Lambda]}\}_{\theta\in\R/\pi\Z}$,
with $\mu_{\theta}^{[\Lambda]}$ supported on
$X(\kappa(\theta)\cdot\Vk)$ for \(\theta\in \R/\pi\Z\), such that for every
$f\in C_c(X_k(\VV_n))$,
\[
\lim_{T\to\infty}
\mu_{T,[\Lambda]}(f)
=
\int_0^\pi
\mu_\theta^{[\Lambda]}(f)\,
w(\theta)\,\dd\theta.
\]
Here
\begin{equation}\label{eq:def-wtheta}
w(\theta_0):=
\frac{
\displaystyle\int_{-\infty}^{\infty}
\|\kappa(\theta_0)\delta(\infty)\nu(s)\|^{-2}
\,\dd s
}{
\displaystyle\int_0^\pi
\int_{-\infty}^{\infty}
\|\kappa(\theta)\delta(\infty)\nu(s)\|^{-2}
\,\dd s\,\dd\theta
},\;\;\theta_0\in \R/\pi\Z,
\end{equation}
where \(\nu(s)=
\begin{psmallmatrix}
1&s\\
0&1
\end{psmallmatrix},
\quad
\delta(\infty)=
\begin{psmallmatrix}
1&0\\
0&0
\end{psmallmatrix}.\)
\end{theorem}
\begin{remark}
    In this paper we restrict to the range \(2\leq k \leq n\) due to our interest. 

    The analogous result for \(k\in \set{1,n+1}\) follows directly from this work.
\end{remark}
Observe that the density \(w\) simplifies to the constant density \(w\equiv \frac1\pi\) whenever the norm is left \(\SO_2(\R)\)-invariant.

We shall describe explicitly the fiber-measures \(\{\mu_{\theta}^{[\Lambda]}\}_{\theta\in\R/\pi\Z}\) for any given \(\bkt\Lambda\in X_k(\VV_n)\). Since further definitions are required, we defer the discussion to Section~\ref{sec:2}.

At this point, let us further expand on the statement of Theorem~\ref{thm:GeneralForm}  and  overview the description of the fiber measures. 

First, notice that the weak-* limits of \(\mu_{T,\bkt\Lambda},\;T\to\infty\)   are supported on the subbundle over \(\CC_k\), i.e.,  
	\[ X(\CC_k):=\pi^{-1}(\CC_k)=\set{[\Lambda]\in X_k(\VV_n):\,\gen{\Lambda}_\R\in\CC_k},\]
    regardless of the starting point \(\bkt\Lambda\). The subbundle \(X(\CC_k)\subset X_k(\VV_n)\) is in general a submanifold of high codimension. 

    The essential reason that \(X(\CC_k)\) plays the role of an attractor is due to the following. Denote
    \begin{equation}\label{eq:DefIwasawaCoordinates}
		\kappa(\theta):=\begin{pmatrix}
			\cos\theta&-\sin\theta\\\sin\theta&\cos\theta
		\end{pmatrix},\quad   \delta(t):=\begin{pmatrix}
			e^{t/2}&0\\0&e^{-t/2}
		\end{pmatrix}, 
		\quad
		\nu(s):=\begin{pmatrix}1&s\\0&1\end{pmatrix}.
	\end{equation}
    
    Observe that $\Vk\in \Gr_k(\VV_n)$, see \eqref{eq:DefDiagHighestWeights}, is generated by weight vectors of the largest first \(k\)-weights of the diagonal group \(\left\{\delta(t)\right\}_{t\in \R}\). Then, given \(\V\in \Gr_k(\VV_n)\), in view of the Iwasawa decomposition \(g=\kappa(\theta)\delta(t)\nu(s)\) for elements \(g\in \SL_2(\R)\), for generic \(s\) and large \(t\), we have that \[g\cdot \V \approx \kappa(\theta)\cdot \V^+_k\in \CC_k,\]
    see e.g. Proposition~\ref{prop:AsympAction}.

    Now let us overview the description  \(\{\mu_{\theta}^{[\Lambda]}\}_{\theta\in\R/\pi\Z}\), and discuss their dependence on \(\bkt\Lambda\in X_k(\VV_n)\) explicitly. 

    It is instructive to first consider the starting points \(\bkt\Lambda\in X_k(\VV_n)\) with \(\Lambda\subset \Vk\). Observe that
	\[\V_k^+= x^{n-(k-1)}\gen{x^{k-1},x^{k-2}y,\dots,y^{k-1}}_{\R}= x^{n-(k-1)}\Sym^{k-1}(\R^2).\]
	Thus, for \([\Lambda]\in X(\V_k^+)\),  we can write  \(\Lambda= x^{n-(k-1)} \Lambda'\), where \(\Lambda'\) is a full rank lattice in \(\VV_{k-1}=\Sym^{k-1}(\R^2)\). Moreover, the action factors as 
	\[g\cdot [x^{n-(k-1)}\Lambda']=[(g\cdot x^{n-(k-1)}) (g\cdot \Lambda')],\,g\in \SL_2(\R),\]
	where the actions on \(x^{n-(k-1)}\) and \(\Lambda'\) are the irreducible \(\SL_2(\R)\)-action on \(\VV_{n-(k-1)}\) and \(\VV_{k-1}\), respectively. Writing \(g=\kappa(\theta)\delta(t)\nu(s)\), the term \(\R^\times( g\cdot x^{n-(k-1)})=\R^\times(\kappa(\theta)\cdot x^{n-(k-1)})\), indicates the corresponding angle of the \(k\)-dimensional subspace containing \(g\cdot [x^{n-(k-1)}\Lambda']\), while  \([g\cdot \Lambda']\) indicates how the \(k\)-lattice $\Lambda'$ transforms in the subspace.
	
	Let \(\Gamma \subset \SL_2(\R)\) be a lattice.  As \(X(\V_k^+)\cong \PGL_k(\R)/\PGL_k(\Z)\), it follows from \cite{Shah1998} (building on \cite{Rat91}) that the closure of \([\Gamma \cdot \Lambda']\) is either finite, or is a finite volume orbit of a closed reductive subgroup \(F\subseteq\PGL_k(\R)\), which contains the image \(\varrho(\SL_2(\R))\) of the irreducible \(\SL_2(\R)\)-representation.
    
	The classification of the intermediate reductive subgroups \(\varrho(\SL_2(\R))\subsetneq F\subsetneq\PGL_k(\R)^\circ\) is well known. The Lie algebra  \(\Lie(F)\)  of such a subgroup is a split simple Lie algebra corresponding to one of (i) the symplectic algebra \(\lie{sp}_{2m}(\R)\)  when \(k=2m\), (ii) the orthogonal algebra \(\lie{so}_{m,m+1}(\R)\) when \(k=2m+1\), (iii) the split exceptional algebra \(\g_{2(2)}\subseteq \lie{so}_{3,4}\) when \(k=7\) (see Appendix Section~\ref{sec:DynkingAppend}).

    We would like to stress that while the orbits for starting points in \(X(\CC_k)\) are described in a simple manner, and while \(X(\CC_k)\) is an attractor in a sense, the analysis of the asymptotic behavior of trajectories for starting points \([\Lambda]\in X_k(\VV_n)\smallsetminus X(\CC_k)\) is more challenging, and does not follow in this work by an approximation arguments using   trajectories in \(X(\CC_k)\). Nonetheless, the resulting fiber-measures
\(\{\mu_{\theta}^{[\Lambda]}\}_{\theta\in\R/\pi\Z}\)
in Theorem~\ref{thm:GeneralForm} are always either uniform on a finite
set or convex combinations of at most two homogeneous measures
associated with connected split simple subgroups of \(\PGL_k(\R)\)
containing the principal image of \(\SL_2(\R)\).

    The dependence of \(\{\mu_{\theta}^{[\Lambda]}\}_{\theta\in\R/\pi\Z}\) on the acting lattice \(\Gamma\)   occurs only for \(\Lambda\) lying in subspaces generated by homogeneous polynomials having a maximal degree common factor, which is described by a certain commensurability condition with respect to \(\Gamma\), see Section~\ref{sec:multisec}. In such cases, the fiber measures are uniform on finitely many points.  
    
    As in \cite[Definition~1.5]{SargentShapira2019}, we call a probability
measure on \(X_k(\R^n)\) an \(m\)-extension if its conditional measure
on almost every fiber is uniform on exactly \(m\) points. In their study of random walks on \(X_2(\R^3)\), Sargent and
Shapira \cite[Section~1.3]{SargentShapira2019} discovered a
\(1\)-extension through a computer experiment. They further realized
it explicitly by a continuous equivariant section over the circle
of isotropic planes. Understanding this example was one of the
motivations for \cite{GLS24},
who gave a more general explanation for irreducible
\(\SL_2(\R)\)-actions on fiber bundles over flag varieties, where the fibers are lattice quotients. Using the description in terms of the action on homogeneous polynomials, the criterion for the existence of an equivariant section and its description become especially transparent, see Proposition~\ref{prop:EquivariantSection}.

In other situations, the fiber-measures \(\{\mu_{\theta}^{[\Lambda]}\}_{\theta\in\R/\pi\Z}\) have no sensitivity to \(\Gamma\), and are determined by rationality conditions on a \(\Z\)-basis generating \(\Lambda\) with respect to a distinguished basis of \(\gen{\Lambda}_{\R}\), see Definition~\ref{def:StandardBasis}.

For example, for subspaces \(\V\in \Gr_k(\VV_n)\) whose generating polynomials don't satisfy certain isotropy conditions with respect to bilinear forms, see Remark~\ref{rem:GenericSubspaces}, it holds that for every \(\Lambda\subset \V\), the fiber-measures \(\mu^{[\Lambda]}_{\theta}\) are given by the unique \(\PGL_k(\R)\)-invariant probability measures on \(X(\kappa(\theta)\cdot \Vk)\cong \PGL_k(\R)/\PGL_k(\Z)\), for every \(\theta\in \R/\pi\Z\).

\medskip
\noindent\textbf{Organization of the paper.}
In Section~\ref{sec:2}, we give the explicit classification of the limiting
fiber measures in Theorem~\ref{thm:GeneralForm}.
In Section~\ref{sec:3}, we reduce norm-ball averages of lattice
orbits to the equidistribution of expanding Wronskian
curves. Section~\ref{sec:4} classifies the lower-dimensional
obstructions to equidistribution and establishes
criteria for standard bases.
In Section~\ref{sec:5}, we prove the linear stability properties
of Wronskian curves needed to apply Yang's
equidistribution theorem.
In Section~\ref{sec:6}, we establish the joint equidistribution
of expanding Wronskian curves and prove the
orbit-limit theorems stated in Section~\ref{sec:2}.
Finally, in Section~\ref{sec:7}, we relate projective
$(n-1)$-lattices to Frobenius numbers through
Kannan's correspondence and prove
Theorems~\ref{thm:frob_on_quad_var}
and~\ref{thm:1.2}.
The appendix provides the classification of reductive
Lie subalgebras used in Section~\ref{sec:4}.

    \medskip
\noindent\textbf{Declaration of generative AI use.} 
The authors used OpenAI’s ChatGPT (GPT-5.6 Sol and GPT-6.1 Sol) to search the literature, revise the manuscript, and assist in proving technical lemmas. The authors take full responsibility for the entire manuscript, including the accuracy of its references and the correctness of all mathematical results and proofs.

\medskip
\noindent\textbf{Acknowledgments.} During the collaboration, all the authors were members of the DASOM (Dynamics and Arithmetic Structures on Manifolds) Center, a global basic research lab supported by the NRF of Korea under Project no. RS-2025-02293115.
The second author is partially supported by BK plus 21 SNU Mathematical Sciences Division.
The second and third authors are supported by the NRF of Korea under Project no. RS-2025-00515082. The third author is an associate member of KIAS (Korea Institute of Advanced Study).

	\section{The classification of the limiting  fiber measures}\label{sec:2}

To describe our results, we need to identify the fiber \(X(\Vk)\) with
\(X_k:=\PGL_k(\R)/\PGL_k(\Z)\) in  a concrete manner described below. In \eqref{eq:DefOfWronskian}, we use
this identification to associate to each \(\bkt{\Lambda}\in X_k(\VV_n)\)
a Wronskian curve \(W_{\bkt{\Lambda}}\) taking values in \(X_k\). We use
this curve to describe the fiber-measures \(\mu^{\bkt\Lambda}_{\theta}\) of Theorem~\ref{thm:GeneralForm}.

\subsection{Linear coordinates}

Recall that \(\Vk=x^{n-(k-1)}\VV_{k-1}\). We equip \(\VV_{k-1}\)
with the ordered basis
\begin{equation}\label{eq:DefOfOrderedBasisForRinfty}
    \hat e_i(x,y):=\frac{x^{k-i}y^{i-1}}{(i-1)!},
    \quad i=1,\ldots,k.
\end{equation}
Writing \(\be_1,\dots,\be_k\) for the standard basis of \(\R^k\),
we define
\[
    \hatRep:\R^k\to \VV_{k-1},
    \quad \hatRep(\be_i):=\hat e_i(x,y).
\]
We use the same notation for the induced isomorphism
\(\hatRep:\GL_k(\R)\to\GL(\VV_{k-1})\). Namely, for \(g\in\GL_k(\R)\),
the matrix of \(\hatRep(g)\) in the basis
\eqref{eq:DefOfOrderedBasisForRinfty} is \(g\), so
\begin{equation}\label{eq:DefRep}
    \hatRep(g)(\hat e_j(x,y))=\sum_{i=1}^k g_{ij}\hat e_i(x,y),
    \quad j=1,\ldots,k.
\end{equation}

Multiplication by \(x^{n-(k-1)}\) then gives the identification
\begin{equation}\label{eq:DefRepForVector}
    \Rep:\R^k\to \Vk, \;\; \Rep(\bv):=x^{n-(k-1)}\hatRep(\bv),\;\bv\in\R^k,
\end{equation}
and we get the induced isomorphism \(\Rep:\GL_k(\R)\to\GL(\Vk)\)
\begin{equation}\label{eq:RepOrderedBasisIrredAction}
    \Rep(g)\bigl(x^{n-(k-1)}\hat P(x,y)\bigr)
    :=x^{n-(k-1)}\hatRep(g)(\hat P(x,y)),
    \;\; g\in\GL_k(\R),\; \hat P(x,y)\in\VV_{k-1}.
\end{equation}

Passing to projective classes of full lattices gives the identification
\begin{equation}\label{eq:DefRepOrderedBasis}
    \Rep:X_k\to X(\Vk),
    \quad \Rep(gy_0):=\bkt{\Rep(g\Z^k)},
\end{equation}
where \(y_0\) is the identity coset in \(X_k\) and \(g\in\GL_k(\R)\)
acts on \(X_k\) through its projective class. 

Finally, let \(\varrho_{k-1}:\GL_2(\R)\to\GL_k(\R)\) be the matrix representation
of the action \eqref{eq:DefOfIrredAction} in the basis
\eqref{eq:DefOfOrderedBasisForRinfty}. For \(h\in\GL_2(\R)\), we have
\begin{equation}\label{eq:DefOfIrredRep}
    h\cdot\hat e_j
    =\sum_{i=1}^k\bigl(\varrho_{k-1}(h)\bigr)_{ij}\hat e_i(x,y),
    \;\; j=1,\dots,k.
\end{equation}
Equivalently, \(\hatRep(\varrho_{k-1}(h))(\hat P)=h\cdot\hat P\) for
every \(\hat P\in\VV_{k-1}\).

Set \(\delta_{k-1}(t):=\varrho_{k-1}(\delta(t))\) and
\(\nu_{k-1}(s):=\varrho_{k-1}(\nu(s))\). Applying these actions to the
basis \eqref{eq:DefOfOrderedBasisForRinfty} shows
\begin{equation}\label{eq:DefOfrhonum}
		\delta_{k-1}(t)=\begin{pmatrix}
			e^{\frac{(k-1)t}{2}}&0&0&\cdots&0\\
			0&e^{\frac{(k-3)t}{2}}&0&\cdots&0\\
			0&0&e^{\frac{(k-5)t}{2}}&\ddots&\vdots\\
			\vdots&\vdots&\vdots&\ddots&0\\\
			0&0&0&\cdots&e^{\frac{-(k-1)t}{2}}
		\end{pmatrix},\;\;
        	\nu_{k-1}(s)=\begin{pmatrix}
			1&s&\frac{s^2}{2!}&\cdots&\frac{s^{k-1}}{(k-1)!}\\
			0&1&s&\cdots&\frac{s^{k-2}}{(k-2)!}\\
			0&0&1&\ddots&\vdots\\
			\vdots&\vdots&\vdots&\ddots&s\\\
			0&0&0&\cdots&1
		\end{pmatrix}.
	\end{equation}

    \subsection{Admissible subgroups}
	Identify \(\lie{pgl}_k(\R)\cong\lie{sl}_k(\R)\) and let \(\be_{i,j}\) denote the \(k\times k\) matrix whose \((i,j)\)- entry equals \(1\) and zero elsewhere.  
	Define the following elements of \(\lie{pgl}_{k}(\R)\) 
	\begin{equation}\label{eq:DefPrincSL2Trip}
		\upE_{k-1}:=\sum_{i=1}^{k-1}\be_{i,i+1},\;\upH_{k-1}=\diag(k-1,k-3,\dots,-(k-1)),\;\upS_{k-1}=\sum_{i=1}^{k-1}i(k-i)\be_{i+1,i},
	\end{equation}
	Let \(\dd\varrho_{k-1}:\lie{sl}_2(\R)\to\lie{pgl}_k(\R)\) be the differential of \eqref{eq:DefOfIrredRep}. By \eqref{eq:DefOfrhonum}, it is immediate that
	\[\dd\varrho_{k-1}(\upE_1)=\upE_{k-1},\;\dd\varrho_{k-1}(\upH_1)=\upH_{k-1}, \dd\varrho_{k-1}\upS_1=\upS_{k-1}.\]
	Thus \((\upE_{k-1},\upH_{k-1},\upS_{k-1})\) is a principal \(\lie{sl}_2(\R)\)-triple, and we denote 
	\begin{equation}\label{eq:DefPrincipalLiesl2}
		\varrho_{k-1}\lie{sl}_2(\R):=\dd\varrho_{k-1}(\lie{sl}_2(\R))=\gen{\upE_{k-1},\upH_{k-1},\upS_{k-1}}_{\R}.
	\end{equation}

    In the following, the matrix Lie algebras act on \(\R^k\) by matrix multiplication from the right.

\begin{definition}\label{def:standard}
	\hypertarget{lab:StandardRealisation}{\textbf{(Standard Admissible subalgebras and subgroups)}} 
    We call a Lie algebra \(\f\subset \lie{pgl}_k(\R)\)  \emph{standard admissible} if it is one of the following simple subalgebras.
	\begin{enumerate}
		\item[(1)] the Lie algebra
		$\varrho_{k-1}\lie{sl}_2(\R)$. 
		\item[(2)] the Lie algebra containing $\varrho_{k-1}\lie{sl}_2(\R)$ and preserving the bilinear forms
		\[q_k(x,y) = \sum_{i=1}^{k} (-1)^{i-1}x_i y_{k+1-i}.\]
		It is isomorphic to \(\lie{sp}_{2m}(\R)\) for  \(k=2m \), and to \(\lie{so}_{m,m+1}(\R)\) for \(k=2m+1\; (m \geq 2)\).
		\item[(3)] the Lie subalgebra of \(\lie{so}_{3,4}\) containing \(\varrho_{6}  \lie{sl}_2(\R)\) and preserving the alternating 3-form 
		\begin{equation}\label{eq:G22form}
			\Phi=dx_{147}-2dx_{156}-dx_{237}+dx_{246}-dx_{345} \in\bigwedge^3(\R^7)^*,
		\end{equation}
		where $dx_{ijk}$ denotes $dx_i \wedge dx_j \wedge dx_k.$ This is the rank-\(2\) exceptional Lie algebra whose associated Lie group preserves the split Octonions. We will denote \(\f:=\g_{2(2)}\). 
        \item[(4)] \(\f = \lie{pgl}_k(\R)\).
	\end{enumerate}  
	A closed connected subgroup \(F\subset\PGL_k(\R)\) is called \emph{standard admissible} if \(\Lie(F)\) is standard admissible.
\end{definition} 
		Note that we consider only \(m \geq 2\) in part (2) as $\lie{sp}_{2}(\R)=\lie{sl}_2(\R),$ and $ \lie{so}_{2,1}(\R)=\varrho_2\lie{sl}_2(\R).$

	\subsection{Wronskian curves and standard bases}
	Given a linearly independent \(k\)-tuple of homogeneous polynomials  \(\boldP(x,y):=(P_1(x,y),\dots,P_k(x,y))\in \VV_n^k\), we define the associated \emph{Wronskian-curve} which is \(k\times k\) matrix
	\begin{equation}\label{eq:DefOfWronskian}
		\Wp(s):=(P_j^{(i-1)}(1,s))_{1\le i,j\le k},\;s\in\R,
	\end{equation}
	where \(P_j^{(\ell)}(1,s)\) denotes the \(\ell\)-th derivative of \(P_j(1,s)\) with respect to the variable \(s\).
	
	The following lemma is well-known. It follows from \cite[Theorem 2]{Bocher1900LinearDependence} by induction for example.  
	\begin{lemma}
		\label{lem:WronskianNonSing}
		$\Wp(s)$ is a singular matrix for finitely many $s$. 
	\end{lemma}
	
	Notice that if \(\widetilde\boldP\) is another ordered basis for \(\V=\gen{P_1,\dots,P_k}_\R\), then there exists \(g\in \GL_k(\R)\) such that \(\boldP=\widetilde\boldP g\), and in particular, 
    \begin{equation}\label{eq:LinTransOfWronsk}
        \Wp(s)=W_{\widetilde \boldP g}(s)=W_{\widetilde\boldP}(s)g,\;\forall s\in \R,
    \end{equation}
    as the Wronskian is linear in the components of $\boldP$.
    Here, \(\boldP(x,y) g\) denotes the usual row times matrix multiplication.
    In particular, the set of singularities \(\set{s\in\R:\det \Wp(s) = 0}\) is invariant under choice of basis \(\boldP\) for \(\V\), and will be denoted by
	\begin{equation}\label{eq:DefOfSubspaceSingularityPts}
		\Sigma_{\V}:=\set{s\in\R:\det \Wp(s) = 0},\;\;\V=\gen{P_1,\dots,P_k}_\R.
	\end{equation}
	
	For \(g\in \GL_k(\R)\), denote by \(\bkt{g}:=\R^\times g\) its image in \(\PGL_k(\R)\), and we let
	\[\bkt{\Wp}(s):=\R^\times \Wp(s).\]

\begin{definition}\label{def:StandardBasis}
		Let \(B^-\) be the Borel subgroup of lower-triangular matrices in \(\PGL_k(\R)\). Given a\stdlist{} subgroup \(F\subset \PGL_k(\R)\), we say that \(\V\in\Gr_k(\VV_n)\) is an \emph{\(F\)-subspace} if there exists an ordered basis \(\boldP\) of \(\V\) such that \(\bkt{\Wp}(s)\in  B^-F\) for all \(s\in\R\smallsetminus\Sigma_{\V}\). We shall call such a basis a \emph{standard basis of $\V$ for \(F\)}.    
	\end{definition}
   
      By \eqref{eq:LinTransOfWronsk}, \(\V\) is an \(F\)-subspace if and only if for any ordered basis \(\tilde\boldP\) of \(\V\) there exists \(g\in\GL_k(\R)\) such that \(\bkt{W_{\tilde\boldP}}(s)\in B^- Fg\) for all \(s\in \R\smallsetminus\Sigma_{\V}\).
	
    The condition \(\bkt{\Wp}(s)\in  B^-F\) defining \(F\)-subspaces is interpreted geometrically via the trajectories of the Wronskian curves in the space of full flags of \(\R^k\) as follows. 
	Let \(\calF\) be the space of full flags in $\R^k$, i.e. the space of sequences of nested subspaces
	\[\boldsymbol{U}_\bullet:=\{0\}\subseteq \boldsymbol{U}_1\subseteq\dots \subseteq \boldsymbol{U}_k=\R^k,\quad \dim \boldsymbol{U}_i=i.\]
	Denote by \(\be_1,\dots,\be_k\) the standard basis vectors of \(\R^k\), and let 
	\[E_\bullet:=0\subseteq \gen{\be_1}_{\R}\subseteq \gen{\be_1,\be_2}\subseteq\cdots\subseteq \R^k.\]
	Consider the right multiplication action of \(\GL_k(\R)\) on the space \(\R^k\) of row vectors. Then, \(\PGL_k(\R)\) acts on \(\calF\) with the stabilizer of \(E_\bullet\) equal to \(B^-\), so that \(\calF\cong B^-\backslash \PGL_k(\R)\) via the map \(B^-\bkt g\mapsto E_\bullet g\). Thus, the basis \(\boldP\) for \(\V\) is standard for \(F\), if the trajectory \(E_{\bullet} \Wp(s),\;s\in\R\smallsetminus \Sigma_{\V}\) lies in an orbit \(E_\bullet F\).
    
    For\stdlist{} \(F\), trajectories \(E_\bullet F\) are well understood, and known to be Zariski-closed (see the proof of Proposition~\ref{prop:ClassificationOfFInBFg}). In particular, there exist concrete homogeneous polynomial equations which determine these flag orbits, as well as a concrete description of the standard bases. This is the content of the following proposition.

\begin{proposition}\label{prop:ExistanceOfStandardBasis} Given a\stdlist{} subgroup \(F\) with \(\f:=\Lie(F)\), the following conditions determine the \(F\)-subspaces and their standard bases.
		\begin{enumerate}[label=(\arabic*)]
			\item\label{enu:Principalsl2Obstruction} Let \(\f=\varrho_{k-1}\lie{sl}_2(\R)\).
             A subspace \(\V\in \Gr_k(\VV_n)\) is an \(F\)-subspace if and only if it is generated by polynomials of the form
			\[  P_i(x,y)=H(x,y) \frac{A(x,y)^{k-i}B(x,y)^{i-1}}{(i-1)!},\;i=1,\dots,k,\]
             where \(H(x,y),\,A(x,y),\,B(x,y)\) are homogeneous polynomials. Moreover, every standard basis for \(F\) is of the form \(\boldP(x,y)=(P_1(x,y),\dots,P_k(x,y))\) with the polynomials of the above form.
			\item\label{enu:G22Obstruction} Let \(\f=\g_{2(2)}\).  A subspace \(\V\in \Gr_7(\VV_n)\) is an \(F\)-subspace if and only if \(\V\) has an ordered basis \(\boldP(x,y)\) such that 
			\[\Phi(\boldP(x,y), \partial_y\boldP(x,y), \mathbf{v})=0,\;\forall\mathbf{v}\in\R^7,\]
			where \(\Phi\) is given by \eqref{eq:G22form}, and where \(\partial_y\boldP(x,y):=(\partial_y P_1(x,y),\dots,\partial_y P_k(x,y))\). All standard bases for \(F\) are of this form.
			\item\label{enu:BiLinearObstruction} Let \(\f=\lie{so}_{m,m+1}\) (or \(\f=\lie{sp}_{2m}(\R))\). A subspace \(\V\in \Gr_{2m+1}(\VV_n)\) (or \(\V\in \Gr_{2m}(\VV_n)\)) is an \(F\)-subspace if and only if \(\V\) is generated by polynomials \(P_1(x,y),\dots,P_k(x,y)\) with \(k=2m+1\) (or \(k=2m\)) such that
            \[q (\partial_y^i\boldP(x,y),\partial_y^j\boldP(x,y))=0,\;\forall 0\leq i,j< \floor{k/2},\]
			where \(q(\boldsymbol{x},\boldsymbol{y})=\sum_{i=1}^{k}(-1)^{i+1}x_i y_{k+1-i}\).  All standard bases for \(F\) are of this form.
			\item Let \(\f=\lie{pgl}_k(\R)\). Every \(\V\in \Gr_k(\VV_n)\) is an \(F\)-subspace and every ordered basis of \(\V\) is standard.
		\end{enumerate}
	\end{proposition}
	\begin{remark}
	The above criteria can also be deduced from \cite[Sections~4.1, 4.4, and~5.1--5.2]{DoubrovZelenko2013}. In Section~\ref{subsec:Explicit Conditions for the Obstruction}, we establish a more refined result on the decomposition of Wronskian curves that we use later. The argument also gives a self-contained proof of the above criteria.
	\end{remark}
	Given $\V\in \Gr_k(\VV_n)$, define $F_{\V}$ to be  the smallest dimensional\stdlist{} subgroup $F$ such that $\V$ is an $F$-subspace. It is unique since the dimensions of groups in the list differ.

	Let us define the Wronskian curve associated with  \([\Lambda]\in X_k(\VV_n)\). For \(\alpha\in\R^\times\), let \(\boldP_1(x,y),\;\boldP_2(x,y)\) be ordered \(\Z\)-bases for \(\Lambda\) and \(\alpha \Lambda\), respectively. Then, there exists \(\gamma \in \GL_k(\Z)\) such that \(\boldP_2=\alpha \boldP_1\gamma\), so that \(W_{\boldP_2}=\alpha W_{\boldP_1}\gamma\). 
    
    Thus, the following trajectory is well defined
	\begin{equation}\label{eq:WronskianOfLambda}
		W_{[\Lambda]}:\R\smallsetminus \Sigma_{\V}\to\PGL_k(\R)/\PGL_k(\Z),\quad W_{[\Lambda]}(s):=\bkt{W_{\boldP}}(s)y_0,\;s\in\R\smallsetminus\Sigma_{\V},
	\end{equation}
	where \(y_0\in \PGL_k(\R)/\PGL_k(\Z)\) denotes the identity coset, \(\V:=\gen{\Lambda}_{\R},\) and where \(\boldP\) is an arbitrary ordered \(\Z\)-basis for \(\Lambda\).
	
	\begin{definition}\label{def:FLambdaSmallestSTDSubgroupContWronsk}
		Given \([\Lambda]\in X_k(\VV_n),\) consider the list of\stdlist{} subgroups  \(F\subseteq\PGL_k(\R)\) containing \(F_{\gen{\Lambda}_\R}\).
        
        We say that \(F\) is \emph{\([\Lambda]\)-rational}, if there exists an \(F\)-standard basis \(\boldP\)
        of \(\gen{\Lambda}_\R\) and there exists \(g\in\GL_k(\R)\) such that \(F[g]y_0\) is homogeneous and such that \(\boldP g\) is a \(\Z\)-basis for \(\Lambda\).  We define \(F_{[\Lambda]}\) to be the \([\Lambda]\)-rational  subgroup of smallest dimension.
	\end{definition}
    \begin{remark}\label{rem:BorelCharacterisationOfFLambda} For our\stdlist{} subgroups, in view of \cite[Theorem~7.8]{BorelHarishChandra1962} and the Borel density theorem, see e.g. 
\cite[Chap.~V]{Raghunathan1972}, given a\stdlist{} subgroup \(F\) and \(\bkt\Lambda\in X_k(\VV_n)\), we have that \(F_{[\Lambda]}\subset F\) if and only if there exists an \(F\)-standard basis \(\boldP\) for \(\gen{\Lambda}_{\R}\) and \(g\in\GL_k(\R)\), such that \(g^{-1}Fg\) is defined over \(\Q\) and $\boldP g$ is an ordered $\Z$-basis of $\Lambda$.

Moreover,\stdlist{} subgroups are defined over \(\Q\), so \(g^{-1}Fg\) is defined over \(\Q\) whenever \(g\in \GL_k(\Q)\). In particular, if \(\V\in \Gr_k(\VV_n)\), \(\boldP\) is an \(F_{\V}\) standard basis,  \(g\in \GL_k(\Q)\) and \(\Lambda\) is generated by \(\boldP g\) then \(F_{\V}=F_{\bkt\Lambda}\).
\end{remark}
\begin{remark}\label{rem:GenericSubspaces}
    Let \(\V\in \Gr_k(\VV_n)\) such that for a given ordered basis  \(\boldP(x,y)\) for \(\V\), and every bilinear form \(q\) on \(\R^k\), which is symmetric of signature \(m,m+1\) if \(k=2m+1\) and  symplectic if \(k=2m\), there exists a \emph{nonzero} polynomial in the list \[q (\partial_y^i\boldP(x,y),\partial_y^j\boldP(x,y)),\;\; 0\leq i,j< \floor{k/2}.\]
   Then \(F_{\V}=F_{\bkt{\Lambda}}=\PGL_k^\circ(\R)\).
\end{remark}
    \subsection{Classification of the fiber-measures}
	In this subsection, we give a precise description of the limits of \eqref{eq:DefNormBallAverages}, depending on \(\Lie(F_{[\Lambda]})\). 

    \subsubsection{The Multisection Case: \(\Lie(F_{[\Lambda]})=\varrho_{k-1}\lie{sl}_2(\R)\) with \(\Gamma\)-special \(\bkt{\Lambda}\)}\label{sec:multisec}
	We begin with a special arithmetic situation occurring only when \(\V\) is an \(\bkt{\varrho_{k-1}(\SL_2(\R))}\)-space and where \(H(x,y)\) appearing in a standard basis, see Proposition~\ref{prop:ExistanceOfStandardBasis}, has maximal degree.
    
	We say that a subspace \(\V\in \Gr_k(\VV_n)\) is \emph{Maximally Divisible (MD)} if 
	\[\V = H(x,y) \gen{x^{k-1},\,x^{k-2}y,\dots,\,y^{k-1}}_{\R}\cong \gen{H(x,y)}_\R\otimes \VV_{k-1},\]
	for some homogeneous polynomial $H(x,y)$ of degree \(n-(k-1)\).  For example, $\Vk$ 
	is MD, where \(H(x,y)=x^{n-(k-1)}\). 
	
	If \(\V\) is MD, 
    then \(g \in \SL_2(\R)\) maps elements of $X(\V)=\set{[H(x,y)\Lambda']:\;[\Lambda']\in X(\VV_{k-1})}  $ as follows:
	\[g\cdot H(x,y)[\Lambda']=(g\cdot H(x,y))[g\cdot \Lambda'].\]
    
	For \([\Lambda']\in X(\VV_{k-1})\), denote
its stabilizer by $\stab[\Lambda']:=\set{h\in\SL_2(\R):\, h\cdot[\Lambda']=[\Lambda']}.$
	\begin{definition}\label{def:GammaSpecialPoints}
		Let \(\V\in\Gr_k(\VV_n)\) be MD.
		\begin{enumerate}
			\item[(1)] We say that \([\Lambda]=[H(x,y)\Lambda']\in X(\V)\) is \emph{\(\Gamma\)-special} if \(\stab[\Lambda']\) and $\Gamma$ are commensurable. This is equivalent to \(\stab[\Lambda']\cap \Gamma\) being a lattice in \(\SL_2(\R)\).
		
			\item[(2)] If \([\Lambda]\in X(\V)\) is \(\Gamma\)-special, then 
			\(\Gamma\cdot[\Lambda']\) is finite, which we call the \emph{\((\Gamma,\Lambda)\)-packet} and denote
			\[\Gamma\cdot [\Lambda']=\set{[\Lambda'_1],\dots,[\Lambda'_m]}.\]
			We define the \emph{\((\Gamma,\Lambda)\)-principal multisection} \(\calM_{\Gamma,\Lambda}\subseteq X(\CC_k)\) by
			\begin{equation}\label{eq:principalMultiSection}
				\calM_{\Gamma,\Lambda}:=\set{\kappa(\theta)\cdot\left(x^{n-(k-1)}\,[\kappa(-\theta)\cdot\Lambda'_i]\right):\;\theta\in \R/\pi\Z,\,i=1,\dots,m}.  
			\end{equation}
		\end{enumerate}
	\end{definition}
	The principal multisection $\calM_{\Gamma,\Lambda}$ is the union of images of continuous sections \(\sigma_{[\Lambda']} : \CC_k\to X(\CC_k)\) defined for  \([\Lambda']\in X(\VV_{k-1})\) as
	\[\sigma_{[\Lambda']}(\theta):= \kappa(\theta)\cdot\left(x^{n-(k-1)}\,[\kappa(-\theta)\cdot\Lambda']\right),\;\theta\in \R/\pi \Z, \]
	which is the map on the fiber \(X\paren{\Vk}\) sending \(x^{n-(k-1)}[\Lambda'] \mapsto x^{n-(k-1)}[\kappa(-\theta)\cdot \Lambda'] \)  composed with the action of \(\kappa(\theta)\), mapping the fiber \(X\paren{\Vk}\)  to \(\kappa(\theta)\cdot X\paren{\Vk}=X\paren{\kappa(\theta)\cdot\Vk}\).
	
	If the \((\Gamma,\Lambda)\)-packet consists of one element \(\set{[\Lambda']}\), then \eqref{eq:principalMultiSection} is the image of the section \(\sigma_{[\Lambda']}\). As the following proposition shows, this section is \(\Gamma\)-equivariant, which generalizes \cite[Theorem~1.10]{SargentShapira2019}. It is stated for general discrete subgroups \(\Gamma\),  of possibly infinite volume.
	
	Since \(\SL_2(\R)\) preserves \(\CC_k\cong \R/\pi\Z\),  \(\SL_2(\R)\) acts on \(\R/\pi\Z\) defined via  \eqref{eq:CircleIdentification}, namely
	\[g\cdot\theta:=\sigma^{-1}\paren{g \kappa(\theta)\cdot \Vk},\;g\in\SL_2(\R),\,\theta\in \R/\pi\Z,\]
	so that \(\kappa(g\cdot \theta)\cdot \Vk=g\kappa(\theta)\cdot \Vk.\) 
	\begin{proposition}\label{prop:EquivariantSection}
		If $\gamma \in \stab[\Lambda'],$ then
		\[\gamma \cdot \sigma_{[\Lambda']}(\theta)= \sigma_{[\Lambda']}(\gamma\cdot \theta).\]
		If \(\Gamma \subset \stab[\Lambda']  \) is a lattice, and \(\theta_0\in \R/\pi\Z\), then
		\[\overline {\Gamma\cdot\sigma_{[\Lambda']}(\theta_0)}=\set{\kappa(\theta)\cdot\left(x^{n-(k-1)}\,[\kappa(-\theta)\cdot\Lambda']\right):\;\theta\in \R/\pi\Z}. \]
	\end{proposition}
	\begin{proof}
		Notice that for \(g\in\SL_2(\R)\) and \(\theta\in\R/\pi \Z\), there exists \(\alpha\in\R^\times\) such that
		\[g\kappa(\theta)\cdot x^{n-(k-1)}=\alpha \kappa(g\cdot \theta)\cdot x^{n-(k-1)}.\]
		
		Therefore, we obtain the equivariance
		\begin{align*}
			\gamma \cdot \sigma_{[\Lambda']}(\theta)&=\gamma\kappa(\theta)\cdot\left(x^{n-(k-1)}\,[\kappa(-\theta)\cdot\Lambda']\right)
			=\paren{\gamma \kappa(\theta)\cdot x^{n-(k-1)}}\,[\gamma\cdot\Lambda']\\
			&=\paren{ \kappa(\gamma\cdot\theta)\cdot x^{n-(k-1)}}\,[\Lambda']
=\kappa(\gamma\cdot\theta)\cdot\left(x^{n-(k-1)}\,[\kappa(-\gamma\cdot\theta)\cdot\Lambda']\right)=\sigma_{[\Lambda']}(\gamma\cdot\theta).
		\end{align*}
		
		If \(\Gamma \subset \stab[\Lambda'] \) is a lattice, then it is well known that every \(\Gamma\)-orbit on \(\SL_2(\R)/B\cong \R/\pi\Z\) is dense. 
        Since the section \(\sigma_{[\Lambda']}\) is continuous, it follows by the equivariance just proved that 
		\[\overline {\Gamma\cdot\sigma_{[\Lambda']}(\theta_0)}
        =\sigma_{[\Lambda']}(\overline{\Gamma\cdot\theta_0})
        =\sigma_{[\Lambda']}(\R/\pi \Z)
        =\set{\kappa(\theta)\cdot\left(x^{n-(k-1)}\,[\kappa(-\theta)\cdot\Lambda']\right):\;\theta\in \R/\pi\Z}. \]
	\end{proof}

	\begin{theorem}\label{thm:MainTheoremGammaSpecial}
		Let \(\V\in \Gr_k(\VV_n)\) be  a Maximally Divisible subspace.  If \([\Lambda]\in X(\V)\) is \(\Gamma\)-special, then the fiber-measures \(\mu^{\bkt\Lambda}_{\theta}\) for \(\theta\in \R/\pi\Z\) given in Theorem~\ref{thm:GeneralForm} are given by
		\[\mu^{\bkt\Lambda}_{\theta}(f)=\frac{1}{m}\sum_{i=1}^m f(\kappa(\theta)\cdot[x^{n-(k-1)}\kappa(-\theta)\cdot \Lambda'_i]),\;f\in C_c(X_k(\VV_n)),\]
		where \(\set{[\Lambda'_1],\dots,[\Lambda'_m]}\)
		is the \((\Gamma,\Lambda)\)-packet.
	\end{theorem}
	\bigskip
	
	\noindent\textbf{Example.} Let \(k<n+1\) be two natural numbers. It holds that \(\Lambda:=x^{n-(k-1)}\Sym^{k-1}(\Z^2):=x^{n-(k-1)}\gen{x^{k-1},x^{k-2}y,\dots,y^{k-1}}_\Z\) is \(\SL_2(\Z)\)-special. Indeed, for every \(\gamma \in \SL_2(\Z)\),
	\[\gamma\cdot\gen{x^{k-1},x^{k-2}y,\dots,y^{k-1}}_\Z\subset \gen{x^{k-1},x^{k-2}y,\dots,y^{k-1}}_\Z.\]
	In particular, the \((\SL_2(\Z),\Lambda)\)-packet is a singleton \(\{[\Lambda']\},\) where \(\Lambda':=\Sym^{k-1}(\Z^2)\).  Thus, if \(\|\cdot\|\) is \(\SO_2(\R)\)-invariant, we get that
	\[\lim_{T \to \infty}
	\frac{1}{\#\Gamma_T}\sum_{\gamma \in \Gamma_T} f(\gamma \cdot [\Lambda])=\frac1\pi \int_0^{\pi}f(\kappa(\theta)\cdot[x^{n-(k-1)}\kappa(-\theta)\cdot \Lambda'])\,\dd\theta,\;\forall f\in C_c(X_k(\VV_n)).\]

    \subsubsection{The Case \(\Lie(F_{\bkt{\Lambda}})=\varrho_{k-1}\lie{sl}_2(\R)\) with not \(\Gamma\)-special \(\bkt{\Lambda}\), or \(\Lie(F_{\bkt{\Lambda}})=\lie{so}_{m,m+1}(\R)\)}
Since these two cases have similar behavior, we introduce the notions separately and then state a unified result.

Let \(\V\) be an \(F\)-subspace, where \(\f\in \set{\rho_{k-1}\lie{sl}_2(\R),\lie{so}_{m,m+1}(\R)}\) and fix  an \(F\)-standard basis \(\boldP(x,y)=(P_1(x,y),\dots,P_k(x,y))\).

If \(\f=\rho_{k-1}\lie{sl}_2(\R)\),  we denote the finite set \(\widetilde \Sigma  := \Sigma \cup \set{s\in \R:P_1(1,s)=0},\) and we let \(\phi_{\boldP,\lie{f}}(s)=P_2(1,s)/P_1(1,s),\; s\in \R\smallsetminus \widetilde\Sigma\). Notice that, with the notations of the standard basis in Proposition~\ref{prop:ExistanceOfStandardBasis}\ref{enu:Principalsl2Obstruction}, we have  \(B(1,s)/A(1,s)=P_2(1,s)/P_1(1,s)\).  

We define 
	\[I_{\boldP,\lie{f}}^\pm:=\set{s\in \R\smallsetminus \widetilde\Sigma:\,\pm\phi_{\boldP,\lie{f}}'(s)>0}.\]Also, denote \(J_{\lie{f}}:=\diag(1,-1,1,\dots,(-1)^{k-1})\). 

    If \(k=2m+1\)  and \(\f=\lie{so}_{m,m+1}(\R)\), define \(\phi_{\f}(h):=\det\left(h_{ij}+(-1)^{m+j}h_{i,\,2m+2-j}\right)_{1\leq i,j\leq m}\) for \(h\in \GL_{2m+1}(\R)\), and let
        \[\phi_{\boldP,\f}(s):=(\det W_{\boldP}(s))^m\phi_{\f}(W_{\boldP}(s)),\;s\in\R\smallsetminus\Sigma,\]
        where \(\Wp\) is the Wronskian curve \eqref{eq:DefOfWronskian}.
		Similarly, let
		\[I_{\boldP,\lie{f}}^\pm:=\set{s\in \R\smallsetminus\Sigma:\pm \phi_{\boldP,\f}(s) >0}.\]
        We denote here $J_{\f}=\diag(-1,I_{2m-1},-1)$.

\begin{remark}
    Let \(U^-\subseteq \PGL_k(\R)\) be the lower-triangular unipotent subgroup.
By Lemmas~\ref{lem:WronskCritForSL2} and~\ref{lem:CurveDecompAndCondForOrthogonal},
\(I_{\boldP,\f}^+\) and \(I_{\boldP,\f}^-\) are finite unions of intervals
and are the inverse images under \(s\mapsto\bkt{W_{\boldP}}(s)\) of the
disjoint sets \(U^-F\) and \(U^-\bkt{J_\f}F\), respectively.

We include the following observation only to clarify the dependence on
the standard basis and omit its proof. A basis \(\boldP h\)
is \(F\)-standard precisely when
\(\bkt{h}\in N_{\PGL_k(\R)}(F)=F\sqcup\bkt{J_\f}F\).
For \(f\in F\) and two \(F\)-standard bases \(\boldP_1,\boldP_2\) such that \(\boldP_2=\boldP_1  J_{\f}\), it holds that \(\phi_{\boldP_1,\f}=-\phi_{\boldP_2,\f}\).
\end{remark}

       We define the following densities
	\begin{equation}\label{eq:PrinSLtheta}
c_{\boldP,\f}^\pm(\theta_0):=\frac{\int_{I_{\boldP,\f}^\pm}\,\|\kappa(\theta_0)\delta(\infty)\nu(s)\|^{-2}\,\dd s}{\int_{-\infty}^\infty\,\|\kappa(\theta_0)\delta(\infty)\nu(s)\|^{-2}\,\dd s},\;\theta_0\in[0,\pi].
	\end{equation} 

   Now let \(\Lambda\subset \V\) such that \(F_{[\Lambda]}=F\), with the above \(F\). In case that \(\f=\rho_{k-1}\lie{sl}_2(\R)\), assume further that \(\Lambda\) is not \(\Gamma\)-special. 
   
   In particular, \(\Lambda\) has an ordered basis \(\boldP(x,y) g\), with \(g\in \GL_k(\R)\) and the fixed ordered basis \(\boldP\) such that the orbits \(\calO^\pm\subset \PGL_k(\R)/\PGL_k(\Z)\) given by
		\[\calO^+:=F_{[\Lambda]} \bkt g y_0,\quad \calO^-:=F_{[\Lambda]}  \bkt{ J_{\f} \,g }y_0,\]
        are homogeneous.

With the above conventions, we have the following theorem.
\begin{theorem}\label{thm:PrincipalSL2Equid}
		The fiber-measures \(\mu^{\bkt\Lambda}_\theta\), for \(\theta\in\R/\pi\Z\) of Theorem~\ref{thm:GeneralForm}, are given by
        \[\mu^{[\Lambda]}_{\theta}= c_{\boldP,\f}^+(\theta) \left(\kappa(\theta)_*R^\infty_* \mu_{\calO^+}\right)
        + c_{\boldP,\f}^-(\theta) \left(\kappa(\theta)_*R^\infty_* \mu_{\calO^-}\right).\]
		In particular, if \(\calO^+ = \calO^-\), then 
		\[\mu^{[\Lambda]}_{\theta}=\kappa(\theta)_*(\Rep_*\mu_{\calO}).\]
	\end{theorem}

	Observe that \(\calO^+ = \calO^- \) if and only if there exists a \(\gamma \in \PGL_k(\Z)\) such that \[\gamma \in \bkt {g^{-1} J_{\f} F_{\bkt{\Lambda}} g}.\]
    For example, if \(g\) is a rational diagonal matrix, then   \(gJ_{\f}g^{-1}=J_{\f}\in \GL_k(\Z)\). Since \(g\) is rational, \(F\bkt g y_0\) is homogeneous, see Remark~\ref{rem:BorelCharacterisationOfFLambda}. Thus, in every such case, \(\calO^+=\calO^-\).
	
	It is more difficult to construct \(\bkt g\) such that for all \(\gamma \in \PGL_k(\Z)\) it holds that \(\gamma \notin \bkt {g^{-1} J_{\f} F_{\bkt{\Lambda}} g},\) and such that \(F[g]y_0\) is homogeneous.
    
    Construction of examples for such \(\bkt g\) for \(F_{[\Lambda]}=\bkt{\varrho_{k-1}(\SL_2(\R))}\) case will be given in Lemma~\ref{lem:DisjointOLambdaOrbsExample}. See also Remark~\ref{rem:ExamplesDisjointOrbOrthCase} regarding such examples for the \(F_{[\Lambda]}=\bkt{\SO^\circ_{m,m+1}(\R)}\) case.

    We would like now to present several concrete examples. We take the norm \(\|\cdot\|\) to be \(\SO_2(\R)\)-invariant.

	\medskip
	
	\noindent\textbf{Example.} If \(k=2\), then \(F_{[\Lambda]}=\bkt{\varrho_1(\SL_2(\R))}=\PGL^\circ_2(\R)\) for every \(2\)-lattice \(\Lambda\subset \VV_n\). In particular, there exists a unique \(\PGL^\circ_2(\R)\)-invariant probability measure \(\mu_{\Vk}\) on \(X(\Vk)\cong \PGL_2(\R)/\PGL_2(\Z)\cong \PGL^\circ_2(\R)/\PGL^\circ_2(\Z)\), where \(\Vk:=\gen{x^n,x^{n-1}y}_\R\).
Therefore, for the following \(\Gamma\) and \(\Lambda\), we have 
	\[\lim_{T\to\infty}\frac{1}{\#\Gamma_T}\sum_{\gamma \in \Gamma_T}f(\gamma\cdot [\Lambda])=\frac 1\pi\int_0^{\pi}\kappa(\theta)_*\mu_{\Vk}(f)\,\dd\theta,\;\forall f\in C_c(X_2(\VV_n)).\]
	\begin{enumerate}
		\item  Consider \(\Lambda= x^{n-1}\gen{\sqrt{2}x,y}_\Z\subset \VV_n,\) and \(\Gamma = \SL_2(\Z)\). Then, \(\Lambda\) lies in the MD subspace \(\Vk\), but \(\Lambda\) is not \(\SL_2(\Z)\)-special, i.e. varying \(\begin{psmallmatrix}
			a&b\\c&d
		\end{psmallmatrix}\in\SL_2(\Z)\) produces an infinite family of homothety classes of \(2\)-lattices
		\(\begin{psmallmatrix}
			a&b\\c&d
		\end{psmallmatrix}\cdot[\gen{\sqrt{2}x,y}_\Z]=[\gen{\sqrt{2}(ax+cy),bx+dy}_\Z].\)
		Indeed, for \(n_1, n_2 \in \N\), if \[\begin{psmallmatrix}
			1&n_1\\0&1
		\end{psmallmatrix}\cdot[\gen{\sqrt{2}x,y}_\Z]=\begin{psmallmatrix}
			1&n_2\\0&1
		\end{psmallmatrix}\cdot[\gen{\sqrt{2}x,y}_\Z],\]
		then there exist \(q_1,q_2\in \Z\) and \(\alpha \in \R^\times \) such that:
		\[q_1\sqrt{2}x+q_2(n_1x+y) = \alpha (n_2x+y)\iff \begin{cases}
			q_1\sqrt{2}+q_2n_1=\alpha n_2\\
			q_2=\alpha
		\end{cases}\]
		Plugging in the second condition in the first shows that \(q_1\sqrt 2 = q_2(n_2-n_1).\)
		So, if \(n_1\neq n_2\), there is no solution.
		
		\item Let \(\Lambda := \gen{x^n,y^n}_\Z\subset \VV_n\). Recall that throughout the paper we assume \(n\geq 2\).  Then \(\gen{\Lambda}_\R\) is not MD. The above limit is obtained for every lattice subgroup \(\Gamma\subseteq \SL_2(\R)\).		
	\end{enumerate}
    
    \noindent\textbf{Example.} Let \(\V:=\gen{x^6,x^5y,x^3y^3,xy^5,y^6}_{\R}\). Then, \(\V\) is an \(\bkt{\SO^\circ_{2,3}(\R)}\)-subspace, but not \(\bkt{\rho_4(\SL_2(\R))}\)-subspace.
        The ordered basis \(\boldP(x,y):=\paren{x^6,\frac{3}{2}x^5y,-\sqrt{\frac{5}{2}}x^3y^3,\frac{3}{2}xy^5,y^6}\) is \(\bkt{\SO^\circ_{2,3}(\R)}\)-standard. Thus, for every $\Lambda\subset \V$ which is generated by an ordered basis of the form \(\boldP(x,y)g\) with \(g\in\GL_5(\Q)\), we have that \(F_{[\Lambda]}=\bkt{\SO^\circ_{2,3}(\R)}\). 
        
        In particular, for $\Lambda=\gen{x^6,\frac{3}{2}x^5y,-\sqrt{\frac{5}{2}}x^3y^3,\frac{3}{2}xy^5,y^6}_{\Z},$ we get that
        \[\lim_{T\to\infty}\frac{1}{\#\Gamma_T}\sum_{\gamma \in \Gamma_T}f(\gamma\cdot [\Lambda])=\frac 1\pi\int_0^{\pi}\kappa(\theta)_*(\Rep_*\mu_{\bkt{\SO^\circ_{2,3}(\R)}y_0})(f)\,\dd\theta,\;\forall f\in C_c(X_5(\VV_6)).\]

		\subsubsection{The Remaining Cases}
		Finally, we  present the results for \(\Lie(F_{[\Lambda]})\in \set{\g_{2(2)},\lie{sp}_{2m}(\R),\lie{pgl}_k(\R)}\) where $m \geq 2$ and $k \geq 3.$ 
		
		\begin{theorem}\label{thm:mainthm_g2sppgl}
			Let \(F\)\stdlist{} with \(\Lie(F)\in \set{\g_{2(2)},\lie{sp}_{2m}(\R),\lie{pgl}_k(\R)}\).  
            Let \(\V\) be an \(F\)-subspace, and fix an \(F\)-standard basis \(\boldP(x,y)\) for \(\V\). 
            
            Assume that \(\Lambda\subset \V\) such that \(F=F_{[\Lambda]}\), and let \(g\in \GL_k(\R)\) such that \(\boldP(x,y)g\) is an ordered \(\Z\)-basis for \(\Lambda\) and  \(\calO:=F_{[\Lambda]} \bkt g y_0\) is homogeneous. Denote by \(\mu_{\calO}\) be the \(F_{[\Lambda]}\)-invariant probability measure on \(\calO\). 
            
            Then, the fiber-measures of Theorem~\ref{thm:GeneralForm} are given by
		\(\mu^{[\Lambda]}_{\theta}=\kappa(\theta)_*(\Rep_*\mu_{\calO})\) for \(\theta\in \R/\pi\Z\).
		\end{theorem}

        \noindent\textbf{Example.} Let $\Lambda = \gen{x^3,xy^2,y^3}_{\Z}.$ The existence of a standard basis satisfying the condition in Proposition~\ref{prop:ExistanceOfStandardBasis}\ref{enu:BiLinearObstruction} is equivalent to the existence of a nondegenerate quadratic form $\tilde q$ of signature (1,2) satisfying \(\tilde q(x^3,xy^2,y^3)=0.\) If such $\tilde q$ exists, then
            \[\tilde q(x^3,xy^2,y^3) = \tilde q_{11} x^6 + 2\tilde q_{12} x^4y^2 + 2\tilde q_{13}x^3y^3 + \tilde q_{22}x^2y^4 + 2\tilde q_{23}xy^5 + \tilde q_{33}y^6=0,\]
            where $\tilde q_{ij}=\tilde q(\be_i,\be_j).$ Since the above six monomials are distinct, all coefficients $\tilde q_{ij}$ must vanish. Hence, $\tilde q=0,$ contradicting the nondegeneracy of $\tilde q$. Therefore, $\Lie(F_{\gen{\Lambda}_{\R}})= \Lie(F_{\bkt{\Lambda}})=\lie{pgl}_{3}(\R).$

		\section{Sums to integrals}\label{sec:3}
		
		Throughout this section, \(G:=\SL_2(\R)\). We fix a norm \(\|\cdot\|\) on $M_2(\R)$, and for a subset \(\calS\subset G\), we denote 
		\[\calS_T:=\set{g\in \calS:\,\|g\|\leq T}.\]

        We now present a special case of \cite[Proposition~2.1]{BersudskyShah} (see also references therein), which is suitable for our needs. We shall state the result in a slightly different form.
        
		Let \(\Gamma\subset G\)  be a lattice, and  fix natural numbers \(2\leq k\leq n\). We stress however that the following holds true for all natural numbers \(1\leq k\leq n+1\). 
        
         Let \(\dd g\) be a  Haar measure on \(G\), and for $T>0,$ for sufficiently large $T$ such that \(G_T\neq\emptyset\), define a probability measure $\mu_T$ on $G/\Gamma \times X_k(\VV_n)$ by
		\begin{equation}\label{eq:DefJointIntegral}
			\mu_T(\varphi \otimes f):=\frac{1}{\vol(G_T)} \int_{G_T} \varphi(g\Gamma)f(g\cdot[\Lambda])\dd g,    \end{equation}
		for $\varphi \in C_c(G/\Gamma)$ and $f \in C_c(X_k(\VV_n)).$ 
        
       Let $\{O_\delta\}_{\delta \in (0,1]}$ be a nested sequence  of neighborhoods  of  the identity coset in $G/\Gamma$ such that $O_\delta \searrow  \{\Gamma\}$ as $\delta \to0.$ For each $\delta \in (0,1]$, choose a nonnegative bump function $\varphi_\delta\in C_c(G/\Gamma)$ supported in $O_\delta$ and $\mu_{G/\Gamma}(\varphi_\delta) =1,$ where \(\mu_{G/\Gamma}\) is the \(G\)-invariant probability on \(G/\Gamma\).

		\begin{proposition}[Proposition~2.1~\cite{BersudskyShah}]\label{prop:sumtointegral}
			Suppose that $\mu_T$ converges in the weak-$\ast$ topology to a measure $\nu$ on $G/\Gamma \times X_k(\VV_n).$ 
            
		Then, for every $f \in C_c(X_k(\VV_n)),$
			$$\lim_{T \to \infty} \frac{1}{\#\Gamma_T} \sum_{\gamma \in \Gamma_T}f(\gamma \cdot [\Lambda]) = \lim_{\delta \to 0} \nu(\varphi_\delta \otimes f).$$
		\end{proposition}
			
    We remark that \cite[Proposition~2.1]{BersudskyShah}  has two conditions \[
\lim_{\delta\to0}
\limsup_{T\to\infty}
\frac{\operatorname{vol}(G_{(1+\delta)T})}
{\operatorname{vol}(G_T)}
=1, \;\; \mathrm{and}\;\;
\frac{\#\Gamma_T}{\operatorname{vol}(G_T)}
\rightarrow \frac{1}{\vol(G/\Gamma)}.
\] These assumptions hold in our case, see for example \cite{GorodnikNevo}.
		\subsection{Iwasawa decomposition and asymptotics for \texorpdfstring{$\SL_2(\R)$}{SL2(R)}-action on \texorpdfstring{$X_k(\VV_n)$}{Xk(VVn)}}
		
		In this section, we prove the following asymptotic expression of the $\SL_2(\R)$-action on 
        $X_k(\VV_n)$ in the Iwasawa decomposition coordinates \eqref{eq:DefIwasawaCoordinates}.
		Roughly speaking, we show that given $[\Lambda]\in X_k(\VV_n)$ and large \(t\), apart from finitely many \(s\), the trajectory \(s\mapsto \delta(t) \nu(s)\cdot[\Lambda]\) is approximated by the Wronskian \eqref{eq:WronskianOfLambda} trajectory in \(\Vk\) \eqref{eq:DefDiagHighestWeights}, expanded by the diagonal element \(\delta_{k-1}\) given by \eqref{eq:DefOfrhonum}. 
		\begin{proposition}\label{prop:AsympAction}
			Let $[\Lambda]\in X_k(\VV_n)$, and let \(W_{[\Lambda]}:\R\smallsetminus \Sigma \to \PGL_k(\R)/\PGL_k(\Z)\) be the Wronskian curve \eqref{eq:WronskianOfLambda}, where \(\Sigma\) is the finite set \eqref{eq:DefOfSubspaceSingularityPts} of singularities of the Wronskian curve \(\Wp(s)\) \eqref{eq:DefOfWronskian} which corresponds to a \(\Z\)-basis \(\boldP\) of \(\Lambda\). 
			
			Fix \(\epsilon>0\), and fix  $f\in C_c(X_k(\VV_n)).$ Then, 
			\[f((\kappa(\theta)\delta(t)\nu(s))\cdot[\Lambda])=f(\kappa(\theta)\cdot\left[ \Rep \left(\delta_{k-1}(t)W_{[\Lambda]}(s)\right)\right]) +o_{t\to \infty}(1),\]
			uniformly in $\theta\in[0,2\pi]$ and $s \in \R\smallsetminus \bigcup_{\xi\in\Sigma}(\xi-\epsilon,\xi+\epsilon)$. Here, \(\Rep:\PGL_k(\R)/\PGL_k(\Z)\to X(\Vk)\) is defined in \eqref{eq:DefRepOrderedBasis}.
		\end{proposition}

        We need several observations before the proof. With the basis elements \( \{\hat e_i(x,y) \}_{i=1,\dots,n+1}\) given by \(\eqref{eq:DefOfOrderedBasisForRinfty}\), for  \(P(x,y)\in \VV_n\) and \(s\in\R\),
				\begin{equation}\label{eq:CoordinatesOfUnipActingOnPloly}
					\nu(s) \cdot P(x,y)= \sum_{i=1}^{n+1}P^{(i-1)}(1,s)\hat e_i(x,y),
				\end{equation}
				where \(P^{(j)}(1,s):=\frac{\dd ^j}{\dd s^j}P(1,s)\) for \(j=0,\dots,n\). Indeed, the equation \eqref{eq:CoordinatesOfUnipActingOnPloly} is immediate on the basis elements, and the conclusion follows by linearity.

				Recall
				\(
				\Vk = \gen{\hat e_1(x,y),\ldots,\hat e_{k}(x,y)}_{\R},
				\)
				and consider the complementary subspace \[\Uk = \gen{\hat e_{k+1}(x,y),\dots,\hat e_{n+1}(x,y)}_{\R}.\]
                We define
				$\PiK:\VV_n\to \Vk$
				to be the projection  \(\sum_{i=1}^{n+1} a_i \hat e_i(x,y)\mapsto  \sum_{i=1}^{k} a_i \hat e_i(x,y).\)
				
				Notice that for \(P(x,y)\in \VV_n\) and \(s\in \R\), by \eqref{eq:CoordinatesOfUnipActingOnPloly}, \[\PiK(\nu(s)\cdot P(x,y))=\sum_{i=1}^{k}P^{(i-1)}(1,s)\hat e_i(x,y).\]

                For the rest of this subsection, fix \[\V=\gen{P_1(x,y),\dots,P_k(x,y)}_
				\R\in \Gr_k(\VV_n).\]
				Let  \(\Wp(s),\;s\in \R\) be the Wronskian \eqref{eq:DefOfWronskian} corresponding to \(\boldP:=(P_1,\dots,P_k)\) and denote by \( \Sigma :=\Sigma_{\Wp}\subset \R\) the finite set \eqref{eq:DefOfSubspaceSingularityPts}.
                
				\begin{corollary}\label{cor:RestrictedProjIsIso}
					For \(s\in \R\),  the  representing matrix of \(\PiK\mid_{\nu(s)\cdot \V}:\nu(s)\cdot \V\to \Vk\)  with respect to the ordered bases \((\nu(s)\cdot P_1,\dots,\nu(s)\cdot P_k)\) and \((\hat e_1,\dots,\hat e_{k})\) is \(\Wp(s)\). 
					
					In particular, the map \(  \PiK\mid_{\nu(s)\cdot\V}\) is an isomorphism onto \(\Vk\),  \(\forall s\notin \Sigma\). 
				\end{corollary}
				
				For $s \notin \Sigma$, we obtain the direct sum decomposition
				\(
				\VV_n=\nu(s)\cdot\V\oplus \Uk,\;\forall s\in \R\smallsetminus \Sigma.
				\)
				Thus, we may define for \(t\in \R\) an isomorphism
				\(
				g(t,s)\in \GL(\VV_n)
				\) by
				\begin{equation}\label{eq:DefOfg}
					g(t,s)\big|_{\nu(s)\cdot\V}:=\delta(t)\cdot \PiK\big|_{\nu(s)\cdot\V},
					\quad
					g(t,s)\big|_{\Uk}:=\delta(t)\cdot \big|_{\Uk}.
				\end{equation}
				Here, we denote by \(\delta(t)\cdot \in \GL(\VV_n)\) the image of \(\delta(t)=\diag(e^{t/2},e^{-t/2})\) via the representation \eqref{eq:DefOfIrredAction}.

				The following is the key observation needed for the proof of Proposition~\ref{prop:AsympAction}.
                \begin{proposition}\label{prop:gtsAssymptotics}
					Let \(\epsilon>0\), and denote \(\Sigma_\epsilon:=\bigcup_{\xi\in\Sigma} (\xi -\epsilon,\xi+\epsilon)\), where \(\Sigma\) is the finite set where \(\Wp(s)\) is singular. Then, for \(s\in \R\smallsetminus \Sigma_\epsilon\),
					\[
					\norm{\delta(t)\cdot g(t,s)^{-1}-\Id}=O_{\epsilon}((1+|s|)^{-1}e^{-t}),\;\text{ as }t\to\infty.
					\]
					In particular, as $t$ tends to infinity, uniformly in \(s\in \R\smallsetminus  \Sigma_\epsilon\),
					\[
					\delta(t)\cdot g(t,s)^{-1}\xrightarrow[t \to \infty]{} \Id,
					\]
				\end{proposition}
				
				\begin{proof}
					For $s\notin\Sigma$, define the linear map \(B_s:=\left(\PiK|_{\nu(s)\cdot \V}\right)^{-1}-\Id|_{\Vk}:\Vk\longrightarrow\Uk\), sending
					\begin{equation}\label{eq:ExpressionBs}
						\sum_{i=1}^{k} P^{(i-1)}_j(1,s)\hat e_i(x,y)\mapsto \sum_{i=k+1}^{n+1} P^{(i-1)}_j(1,s)\hat e_i(x,y),\;\;j=1,\dots,k.
					\end{equation}
					
					We claim that \[\delta(t)\cdot g(t,s)^{-1}(x)-x=\delta(t)\cdot B_s \delta(t)^{-1}\cdot v,\]
					where \(v:=\PiK(x)\). Indeed, by definition of \(g(t,s)\) \eqref{eq:DefOfg}, for \(x =v + w\) with \(v\in \Vk, w \in \Uk,\) we have
					\[\delta(t)\cdot g(t,s)^{-1}(x)=\delta(t) (\delta(t)\cdot  {\Pi_k^{+}|_{\nu(s)\cdot\V}})^{-1}(v)+ w=\delta(t) B_s \delta(t)^{-1}\cdot v+v+w=\delta(t) B_s \delta(t)^{-1} v+x.\]
					Let $ \| \cdot \|$ be an operator norm. Then,
					\begin{equation}\label{eq:DistanceFromIdFormula}
						\| \delta(t)\cdot g(t,s)^{-1}-\Id\|
						\ll \|\delta(t)\cdot |_{\Uk} \| \|B_s \| \| \delta(t)^{-1}\cdot |_{\Vk}\|,\;\forall s\notin \Sigma.
					\end{equation}
					Note that the largest $\delta(t)$-weight on $\Uk$ is $\frac n2-k$, whereas the largest
					$(\delta(t))^{-1}$-weight on $\Vk$ is $-\frac n2+(k-1)$.
					Thus \(\|\delta(t)\cdot |_{\Uk} \| \|\delta(t)^{-1}\cdot |_{\Vk}\| \ll e^{-t}\). Consequently, by \eqref{eq:DistanceFromIdFormula}, 
					\[\|\delta(t)\cdot g(t,s)^{-1}-\Id\|\ll e^{-t}\|B_s\|, \quad  s\notin \Sigma.\]
					
					Let \(\epsilon>0\). We will now estimate \(\|B_s\|\) in the domain \(s\in\R\smallsetminus\Sigma_\epsilon\). 
					First,  fix \(s_0>0\) large enough such that \([-s_0,s_0]\supset \Sigma_\epsilon\). Due to continuity, 
					\begin{equation}\label{eq:O1BoundOnBs}
						\|B_s\|=O_{\epsilon}(1), \;\; s \in [-s_0,s_0]\smallsetminus\Sigma_\epsilon.
					\end{equation}
					Next, let \(|s|\geq s_0\). 
					To express \(B_s\) more conveniently, define the maps $J_s:\V\to \Vk,\;K_s:\V\to \Uk$ by
					\[J_s(f):=\sum_{i=1}^{k} f^{(i-1)}(1,s)\hat e_i(x,y),\quad
					K_s(f) := \sum_{i=k+1}^{n+1} f^{(i-1)}(1,s)\hat e_i(x,y).\]
					Then,   by \eqref{eq:ExpressionBs}, we get that \(B_s=K_s\circ J_s^{-1}\). 
					Let us estimate \(\|B_s\|\) by estimating 
					\(\|K_s \circ J_s^{-1}\hat e_i(x,y)\|,\ i=1,...,k.\)

					Denote 
					\(\phi_{i,s}(x,y):=J_s^{-1}\hat e_i(x,y) \in \V,\;|s|\geq s_0\) for \(i= 1,\dots, k\). Namely, by the definition of \(J_s\), \(\phi_{i,s}(x,y)\in\V\) is a polynomial such that 
					\begin{equation}\label{eq:JsBasis}
						\phi^{(j-1)}_{i,s}(1,s)=\begin{cases}
							1&\text{if }i=j,\\
							0&\text{if }i\neq j,
						\end{cases}\;\text{ for } j = 1,\dots, k.
					\end{equation}
					Also, note $K_s\circ J_s^{-1}\hat e_i(x,y)=\sum_{\ell=k+1}^{n+1} \phi^{(\ell-1)}_{i,s}(1,s)\hat e_\ell(x,y)$. 
					
					Consider the functionals \(\lambda^{(i)}(s)\in \V^*,\;|s|\geq s_0,\)
					\[\lambda^{(i)}_s:\V\to\R,\quad \lambda^{(i)}_s(f):=f^{(i)}(1,s),\quad i\in \N\cup\{0\}.\]
					For every $|s|\geq s_0$, \(J_s\) is an isomorphism, so that $\set{ \phi_{i,s}}_{i=1}^{k}$ is a basis for \(\V\). Then, by \eqref{eq:JsBasis} the first \(k\) functionals $\lambda_s^{(0)},\dots,\lambda_s^{(k-1)}$ form a basis for \(\V^*\) dual to the basis $\set{ \phi_{i,s}}_{i=1}^{k}$ of \(\V\).  
					
					Then, for every $\ell\geq k$, there exist coefficients \(\paren{c_{\ell,j}(s)}_{j=0}^{k-1}\) such that
					\[\lambda_s^{(\ell)}=\sum_{j=0}^{k-1}c_{\ell,j}(s)\lambda_s^{(j)}.\]
					We observe that 
					\[K_s\circ J_s^{-1}\hat e_i(x,y)=\sum_{\ell=k+1}^{n+1}\lambda_s^{(\ell-1)}(\phi_{i,s}) \hat e_{\ell}(x,y) =\sum_{\ell=k+1}^{n+1}c_{\ell-1,i-1}(s)\hat e_{\ell}(x,y),\;\;\text{ for }i = 1,\dots, k.\]
					Then, to complete the proof, it remains to show that $c_{\ell,i}(s)=O(|s|^{i-\ell})$ with \(0\leq i\leq k-1\) and \(k\leq \ell\leq n\).
					
					Choose a basis of polynomials ${q_1,...,q_k}$ for $\V$ such that  the polynomials \(q_1(1,s),\dots,q_k(1,s)\) are  monic and their degrees \(d_i:=\deg(q_i(1,s))\) are increasing  $d_1<d_2<\cdots<d_k$. Denote the falling factorial as
					\[\quad(x)_i:=x(x-1)\cdots(x-i+1),\;\; i \geq1, \quad (x)_0:=1.\]
					By applying the derivatives directly, we have 
					\(\lambda_s^{(\ell)}(q_i)=(d_i)_\ell s^{d_i-\ell}+O( s^{d_i-\ell-1})\) with the convention that both sides are zero for $\ell >d_i.$
					On the other hand, 
					\[\lambda_s^{(\ell)}(q_i)=\sum_{j=0}^{k-1}c_{\ell,j}(s)\lambda_s^{(j)}(q_i)= \sum_{j=0}^{k-1}c_{\ell,j}(s) \left(s^{d_i-j}(d_i)_j+O(s^{d_i-j-1})\right).\] 
					Dividing by $s^{d_i-\ell}$, we conclude that
					\begin{equation}\label{eq:FallinFactExperession}(d_i)_\ell+O(s^{-1})=\sum_{j=0}^{k-1}c_{\ell,j}(s)s^{\ell-j}(d_i)_j(1+O(s^{-1})),\;\;i=1,...,k.
					\end{equation}
					Set $A(s):=\max\{|c_{\ell,j}(s)s^{\ell-j}|:\;j\leq k-1\},$ and suppose for contradiction that $\sup_{|s|\geq s_0}A(s)=\infty.$
					Then the vectors
					\[\widehat c(s):=\left(\frac{c_{\ell,j}(s)s^{\ell-j}}{A(s)}\right)_{j=0,\dots k-1}\]
					lie in a unit sphere corresponding to a suitable norm. By compactness, there exists a  sequence $\{s_\nu\}_{\nu\in\N}$ such that both $|s_\nu|$ and $A(s_\nu)$ diverge to infinity and 
					\[\lim_{\nu\to\infty}\widehat c(s_\nu)=(\alpha_0,\dots,\alpha_{k-1})\neq 0.\]
					Then, the polynomial
					$ F(x):=\sum_{j=0}^{k-1}\alpha_j(x)_j$ is non-zero and by \eqref{eq:FallinFactExperession},
					\[F(d_i)= \lim_{\nu \to \infty} \frac{(d_i)_\ell}{A(s_\nu)}  = 0,\;\;i=1,\dots,k.\]
					However, $F(x)$ is a polynomial of degree at most $k-1$ having $k$ zeros, which is a contradiction.
					
					We conclude that for $|s| \geq s_0,$
					$\norm{\delta(t)\cdot g(t,s)^{-1}-\Id}=O_{\epsilon}(s^{-1}e^{-t}),\;\text{ as }t\to\infty.$
					Combining with \eqref{eq:O1BoundOnBs}, we have \[
					\norm{\delta(t)\cdot g(t,s)^{-1}-\Id}=O_{\epsilon}((1+|s|)^{-1}e^{-t}),\;\text{ as }t\to\infty.
					\]
				\end{proof}

				\begin{proof}[Proof of Proposition~\ref{prop:AsympAction}]
					Let \([\Lambda]\in X_k(\VV_n)\), \(\Lambda=\gen{P_1(x,y),\dots,P_k(x,y)}_{\Z}\), and denote $\V:=\gen{\Lambda}_\R$. 
					
					By
					\eqref{eq:CoordinatesOfUnipActingOnPloly}, the \(k\)-lattice $\nu(s)\cdot \Lambda$ is generated
					by the vectors
					\[
					\nu(s)\cdot P_j(x,y)=\sum_{i=1}^{n+1} P_j^{(i-1)}(1,s) \hat e_i(x,y),
					\quad j=1,\dots,k.
					\] 
					By the definition \eqref{eq:DefOfg} of $g(t,s)$ on $\nu(s)\cdot \V$, for every \(1\leq j\leq k\),
					\begin{align*}
						g(t,s)\,\nu(s)\cdot P_j(x,y)
						&=\delta(t)\cdot\PiK\,(\nu(s)\cdot P_j(x,y))
						\\&=\delta(t)\cdot\paren{\sum_{i=1}^{k} P^{(i-1)}_j(1,s)\hat e_i(x,y)}\\
						&=e^{\frac{n-(k-1)}{2}t}\paren{\sum_{i=1}^{k} P^{(i-1)}_j(1,s)e^{\frac{k-1-2(i-1)}{2}t}\hat e_i(x,y)}\\
						&=e^{\frac{n-(k-1)}{2}t}\Rep(\delta_{k-1}(t)\Wp(s) \be_j),
					\end{align*}
					where \(\Rep\) is given in \eqref{eq:DefRepForVector}, where \(\Wp\) is the Wronskian \eqref{eq:DefOfWronskian} and where \(\delta_{k-1}(t)\) is given in \eqref{eq:DefOfrhonum}.
					
					Therefore,
					\begin{equation}\label{eq:gReducesToWronskian}
						g(t,s)\,\nu(s)\cdot [\Lambda]=\bigl[e^{\frac{n-(k-1)}{2}t}\Rep\bigl(\delta_{k-1}(t)\,\Wp(s)y_0\bigr)\bigr]=\bigl[\Rep\bigl(\delta_{k-1}(t)W_{[\Lambda]}(s)\bigr)\bigr].
					\end{equation}
					
					Denote  \(\kappa(\theta)\cdot\,\in\GL(\VV_n)\), the image of \(\kappa(\theta)\) obtained by the action \eqref{eq:DefOfIrredAction}. We conclude by Proposition~\ref{prop:gtsAssymptotics} and by compactness that
					\[\|\kappa(\theta)\cdot (\delta(t)\cdot g(t,s)^{-1})(\kappa(\theta)^{-1}\cdot)\,-\Id\|=O_{\epsilon}((1+|s|)^{-1}e^{-t}),\]
					uniformly in \(\theta\in\R\),  $s\in\R\smallsetminus \bigcup_{\xi\in\Sigma}(\xi -\epsilon,\xi+\epsilon)$ and \(t\geq 0\).
					
					Let \(f\in C_c(X_k(\VV_n))\). As the left regular action $(L_hf)(x)=f(hx)$ is uniformly
					continuous in the sup norm, i.e. $\|L_h f-f\|_\infty\to0$ as $h\to\Id$ in $\GL(\VV_n)$, we conclude that
					\begin{align*}
						f((\kappa(\theta)\delta(t)\nu(s))\cdot[\Lambda])&=L_{\kappa(\theta)\cdot (\delta(t)\cdot g(t,s)^{-1})(\kappa(\theta)^{-1}\cdot)}f(\kappa(\theta)\cdot\left[ g(t,s)\nu(s)\cdot \Lambda\right])\\
						&=f(\kappa(\theta)\cdot\left[ \Rep \left(\delta_{k-1}(t)\Wp(s)y_0\right)\right])+o_{t\to\infty}(1),  
					\end{align*}
					uniformly in $\theta$ and $s\in\R\smallsetminus \bigcup_{\xi\in\Sigma}(\xi -\epsilon,\xi+\epsilon)$.  
				\end{proof}

				\subsection{Volume estimates and joint integral asymptotics}

				Using the Iwasawa decomposition, a \(G\)-Haar measure $dg$ is given by
				\begin{equation}\label{eq:iwasawa-decomp}
					\int_G\,f(g)\,\dd g= \int_0^{2\pi} \int_{-\infty}^{\infty}\int_{-\infty}^{\infty}\,f(\kappa(\theta)\delta(t)\nu(s))e^{t}\,\dd t\,\dd s\,\dd \theta,
				\end{equation}
               
				\begin{lemma}\label{lem:AsympDomConvLemma}
					Let \(f:G\to \R\) be a bounded Borel measurable function. Then,
					\[\frac{1}{T^2}\int_{G_T}f(g)\,\dd g = \int_{0}^{2\pi}\int_{-\infty}^\infty \int_{\calE(r,\theta)}f(\kappa(\theta)\delta(2\log T+r)\nu(s))\,\dd s\, e^r \dd r \,\dd\theta+o_{T\to\infty}(1),\]
					where 
					\begin{equation}\label{eq:DefCalEAndDeltaInfty}
						\calE(r,\theta):=\set{s\in \R:\,e^{r/2}\|\kappa(\theta)\de(\infty)\nu(s)\|\leq 1},\quad \de(\infty):=\begin{psmallmatrix}
							1&0\\0&0
						\end{psmallmatrix}.
					\end{equation}
					
				\end{lemma}
				\begin{proof}
					We have
					\[\frac{1}{T^2}\int_{G_T}f(g)\,\dd g =\frac{1}{T^2}\int_{\calB_T}\,f(\kappa(\theta)\delta(t)\nu(s))e^{t}\,\dd s\,\dd t\,\dd \theta,\]
					where \(\calB_T:=\set{(\theta,t,s)\in[0,2\pi]\times\R^2:\,\|\kappa(\theta)\delta(t)\nu(s)\|\leq T}.\) 
					
					For \(T>0\), we change the variable \(t\) as \(t=:2\log T +r,\) so that
					\begin{align}
						\frac{1}{T^2}\int_{G_T}f(g)\,\dd g &=\frac{1}{T^2} \int_{\calE_T}\,f(\kappa(\theta)\delta(2\log T+r)\nu(s))e^{r+2\log T}\,\dd s\,\dd r\,\dd \theta\nonumber\\
						&= \int_{\calE_T}\,f(\kappa(\theta)\delta(2\log T+r)\nu(s))e^{r}\,\dd s\,\dd r\,\dd \theta\label{eq:IntOnBallInIwasChangeCoord},   
					\end{align}
					where 
						$\calE_T :
=\set{(\theta,r,s)\in[0,2\pi]\times\R^2:\,\|\kappa(\theta)\frac{\delta(2\log T+r)}{T}\nu(s)\|\leq 1}. $
			
					For each \(r\in \R\), we have \(\lim_{T\to\infty}\frac{1}{T}\delta(2\log T+r)=e^{r/2}\delta(\infty)\), thus, by norm continuity for each \(\theta, r,s\in \R\), we have \(\lim_{T\to\infty}\|\kappa(\theta)\frac{\delta(2\log T+r)}{T}\nu(s)\|=e^{r/2}\|\kappa(\theta)\delta(\infty)\nu(s)\|\). 
					It follows that \(\underset{T \to \infty}{\lim}\chi_{\calE_T}=\chi_{\calE}\) almost everywhere, where
				 \(\chi_{\calE_T}\) and \(\chi_{\calE}\) are the characteristic functions on \(\calE_T\) and \[\calE:=\set{(\theta,r,s)\in[0,2\pi]\times\R^2:\,e^{r/2}\|\kappa(\theta)\delta(\infty)\nu(s)\|\leq 1}.\] 
					Observe that 
                    \[\calE_T\subset \calD:=\set{(\theta,r,s)\in[0,2\pi]\times\R^2: Ce^{r/2}(1+|s|)\leq 1}, \]
                    since there exists \(C>0\) such that for all \(T>0\), \(Ce^{r/2}(1+|s|)\leq \|\kappa(\theta)\frac{\delta(2\log T+r)}{T}\nu(s)\|\).
					
                    Note that \(\chi_{\calD}\) is integrable, since
					\[\int_0^{2\pi}\int_{\R^2} \chi_{\calD}(\theta,r,s) e^r \dd r\,\dd s\, \dd\theta\ll \int_{-\infty}^\infty (1+|s|)^{-2}\dd s.\]
					By the dominated convergence theorem,
					\[\lim_{T\to \infty}\int_0^{2\pi}\int_{\R^2} |\chi_{\calE_T}(\theta,r,s)-\chi_{\calE}(\theta,r,s)|e^r\,\dd s\,\dd r \,\dd \theta=0.\]
					Therefore, by the triangle inequality and using that \(f\) is bounded, by \eqref{eq:IntOnBallInIwasChangeCoord},  
					\begin{align*}
						\frac{1}{T^2}\int_{G_T}f(g)\,\dd g&=\int_{\calE}\,f(\kappa(\theta)\delta(2\log T+r)\nu(s))e^{r}\,\dd s\,\dd r\,\dd \theta + o_{T\to\infty}(1)\\
						&= \int_{0}^{2\pi}\int_{-\infty}^\infty \int_{\calE(r,\theta)}f(\kappa(\theta)\delta(2\log T+r)\nu(s))\,\dd s\, e^r \dd r \,\dd\theta+o_{T\to\infty}(1).
					\end{align*}
				\end{proof}
			
				\begin{corollary}\label{cor:volume-asymp-SL2}
					For $G_T=\{g\in G:\;\|g\|\leq T\}$, we have
					\[ \lim_{T\to\infty}\frac{\vol(G_T)}{T^2}= 2\int_0^{\pi} \int_{-\infty}^{\infty} \|\kappa(\theta)\delta(\infty)\nu(s)\|^{-2}\,\dd s\,\dd\theta.
					\]
				\end{corollary}
				\begin{proof}
					Applying Lemma~\ref{lem:AsympDomConvLemma} with the constant function \(f\equiv 1\), and using Fubini, we obtain
					\begin{align*}
						\frac{1}{T^2}\vol(G_T) &= \int_0^{2 \pi} \dd \theta\int_{-\infty}^{\infty} \dd s \int_{\set{r\in\R:e^{r/2}\leq \|\kappa(\theta)\de(\infty) \nu(s)\|^{-1}}}\,e^r\dd r+o_{T\to\infty}(1)\\
						&=\int_0^{2\pi}\int_{-\infty}^{\infty}\|\kappa(\theta)\delta(\infty)\nu(s)\|^{-2}\,\dd s\,\dd\theta+o_{T\to\infty}(1)\\
						&=2\int_0^{\pi}\int_{-\infty}^{\infty}
						\|\kappa(\theta)\delta(\infty)\nu(s)\|^{-2}\,\dd s\,\dd\theta+o_{T\to\infty}(1).
					\end{align*}
				\end{proof}
				
				Here is our main conclusion of this section.
				We define
				\begin{equation}\label{eq:DefJointDiag}
					\bdelta_{k-1}(t):=(\delta(t),\delta_{k-1}(t)),\;t\in \R,
				\end{equation}
				where \(\delta_{k-1}(t)=\varrho_{k-1}(\delta(t))\) as in \eqref{eq:DefOfrhonum}, and where \(\delta(t):=\delta_1(t)\).
				
                Let \([\Lambda]\in X_k(\VV_n)\), and let \(\boldP\) be an ordered \(\Z\)-basis for \(\Lambda\). Recall the Wronskian curve \(\Wp:\R\smallsetminus\Sigma\to \GL_k(\R)\)  given by \eqref{eq:DefOfWronskian}, where   \(\Sigma=\set{s\in\R:\,\det W_{\boldP}(s)=0}\). As noted in the introduction, \(W_{[\Lambda]}(s):=\bkt{W_{\boldP}}(s)y_0,\;s\in \R\smallsetminus\Sigma\), where \(y_0\in \PGL_k(\R)/\PGL_k(\Z)\) is the identity coset, is well defined in terms of \([\Lambda]\) only. We define
\begin{equation}\label{eq:DefJointWronskianCurve}
                \varphi_{\boldP}:\R\smallsetminus \Sigma \to \SL_2(\R) \times \PGL_k(\R),\;\;\varphi_{\boldP}(s):=(\nu(s),\bkt{W_{\boldP}}(s)),
                \end{equation}
                and as a consequence, the following curve is well defined
				\begin{equation*}
    \varphi_{[\Lambda]}:\R\smallsetminus \Sigma \to G/\Gamma \times X_k,\;\;\varphi_{[\Lambda]}(s):=\varphi_{\boldP}(s)\bx_0,
				\end{equation*}
				where \(\bx_0:=(x_0,y_0)\in G/\Gamma\times X_k\) denotes the identity coset.
                \begin{proposition}\label{prop:ReducBallAvergesToExpandingCurves}
					Let \(\Phi\in C_c(G/\Gamma\times X_k(\VV_n)),\) and for \(\theta\in [0,2\pi]\) denote \(\Phi_\theta(a,b):=\Phi(\kappa(\theta)a,\kappa(\theta)\cdot b).\) Then
					\begin{align*}
						\frac{1}{T^2}\int_{G_T}\Phi(g\Gamma,g\cdot[\Lambda])\,\dd g &\\= \int_0^{2\pi}\int_{-\infty}^\infty \int_{\calE(r,\theta)\smallsetminus\Sigma}&\Phi_{\theta}\circ(\Id\times\Rep)\bigl(\bdelta_{k-1}(2\log(T)+r)\varphi_{[\Lambda]}(s)\bigr)\,\dd s\, e^r\dd r\,\dd\theta+o_{T\to\infty}(1). 
					\end{align*}
				\end{proposition}
				\begin{proof}
					By Lemma~\ref{lem:AsympDomConvLemma} and the definition of $\Phi_\theta,$
					\begin{align}
						\frac{1}{T^2}& \int_{G_T}\Phi(g\Gamma,g\cdot[\Lambda])\,\dd g \notag \\
						&= \int_0^{2\pi}\int_{-\infty}^\infty \int_{\calE(r,\theta)} \Phi_\theta (\delta(2\log T+r)\nu(s)\Gamma,\delta(2\log T +r)\nu(s) \cdot [\Lambda]) \dd s e^r\dd r \dd \theta +o_{T\to\infty}(1). \label{eq:ReducingBallAverage}
					\end{align}
					Since $\Phi_\theta$ is a compactly supported continuous function, Proposition~\ref{prop:AsympAction} shows that the integrand of \eqref{eq:ReducingBallAverage} can be written as
					\begin{align}
						\Phi_\theta (\delta(2\log T+r)\nu(s)\Gamma, &[\Rep(\delta_{k-1}(2\log T+r) W_{[\Lambda]}(s)]) + o_{(2\log T+r) \to \infty}(1) \notag \\
						&= \Phi_\theta\circ(\Id\times\Rep)(\bdelta_{k-1}(2\log T+r)\varphi_{[\Lambda]}(s)) + o_{(2\log T +r) \to \infty}(1). \label{eq:ReducingIntegrand}
					\end{align}
					If $\lim_{t\to \infty} g(\theta,t,s)=0$ for a bounded function $g$, then a dominated convergence argument gives
					\[\lim_{T\to\infty} \int_0^{2\pi} \int_{-\infty}^\infty \int_{\calE(r,\theta)} g(\theta,2\log T +r,s) \dd s e^r \dd r \dd \theta = 0.\]
					Therefore, the error term in \eqref{eq:ReducingIntegrand} remains negligible after integration in \eqref{eq:ReducingBallAverage}, and the result follows.
				\end{proof}
				
				We note that   Proposition~\ref{prop:ReducBallAvergesToExpandingCurves} reduces the analysis of  the asymptotics of \eqref{eq:DefJointIntegral} with \(T\to\infty\) to the analysis of the following measures.
				
				For $I$ a nonempty bounded interval in \(\R \smallsetminus \Sigma\), we define the measures 
				\begin{equation}\label{eq:DefJointMeasure}
					\mu_{I,\tau}(\Phi):=\frac{1}{|I|}\int_I\Phi\paren{\bdelta_{k-1}(\tau)\varphi_{[\Lambda]}(s)}\,\dd s,\;\;\Phi\in C_c(G/\Gamma\times X_k),\;\tau>0.  
				\end{equation}
				Once the limits of \(\mu_{I,\tau}\) as \(\tau\) tends to \(\infty\) are obtained, together with Proposition~\ref{prop:ReducBallAvergesToExpandingCurves} and the dominated convergence theorem, we will be able to calculate the limits of \eqref{eq:DefJointIntegral} and to obtain our main theorems.

				\section{Lower-dimensional obstructions}\label{sec:4}
				To motivate the work in this section, we note the following observation which is crucial to the analysis of the measures \eqref{eq:DefJointMeasure}. All intervals in this section are assumed to have nonempty interior.
				
				Let \[\bB^-:=B_1^-\times B_2^-\subset \SL_2(\R)\times \PGL_k(\R)\] 
				be the product group of the lower-triangular Borel subgroups. For \(i=1,2\), let \(U_i^-\subset B_i^-\) be the lower triangular unipotent subgroups, let \(T_i\subset B_i^-\) the full diagonal subgroups, and denote
				\[\bU^-:=U_1^-\times U_2^-,\;\bT:=T_1\times T_2 \subset \bB^-,\]
				so that \(\bB^-=\bU^-\bT.\)
				
				Recall the diagonal subgroup \(\bdelta_{k-1} \subset \SL_2(\R)\times \PGL_k(\R)\) given in \eqref{eq:DefJointDiag}. Then \(\bdelta_{k-1}\) contracts \(\bU^-\), namely, for every \(\bu^-\in \bU^-\) we have \[\lim_{t\to\infty}\bdelta_{k-1}(t)\bu^-\bdelta_{k-1}(t)^{-1}=\Id.\] 
				Also, since \(\bdelta_{k-1}\subset\bT\), then clearly \(\bdelta_{k-1}\) commutes with \(\bT\).
				
				Let \([\Lambda]\in X_k(\VV_n)\), and let \(\varphi_{[\Lambda]}(s):=\varphi_{\boldP}(s)\bx_0\) be the curve \eqref{eq:DefJointWronskianCurve}. Suppose that there exists a closed subgroup \(\bF\subset \SL_2(\R)\times \PGL_k(\R)\) and \(p\in \SL_2(\R)/\Gamma \times \PGL_k(\R)/\PGL_k(\Z)\) such that \(\varphi_{[\Lambda]}(s)\in \bB^- \bF p\) for all \(s\in I\). Then, we may decompose the curve as 
                \begin{equation}\label{eq:DesiredDecomp}
                    \varphi_{[\Lambda]}(s)=\bu^-(s) \balpha(s) \boldf(s)p,\,s\in I
                \end{equation}
				where \(\bu^-(s)\in \bU^-,\) and \(\balpha(s)\in \bT\), see e.g. Lemma~\ref{lem:BFDecompLemma} below. Consequently,  the support of the limiting measures of \eqref{eq:DefJointMeasure} as \(\tau\to\infty\) will be contained in the closure of \(\bT \bF p\).
				
				In view of \cite[Theorem~1.5]{Yang2020}, after verifying an avoidance property of the trajectories \(s\mapsto (\nu(s),W(s))^{-1}\) acting on full flag varieties avoiding certain Schubert varieties  (see Section~\ref{sec:5}), the only groups \(\bF\) which will be needed to be considered for the above purpose are algebraic reductive subgroups.

                An additional important goal here will be to determine the precise form of the values of \(\balpha(s)\) in the above decomposition, explicitly. Indeed, in Subsection~\ref{subsec:Explicit Conditions for the Obstruction}, we show that \(\balpha(s)\), modulo \(\bF\), is either \((\Id,\Id)\) or \((\Id,\bkt{J})\), with \(J^2=\Id\).

				Let \( \lie{pgl}_k(\R)\cong\lie{sl}_k(\R)\), which we identify with the Lie algebra of traceless \(k\times k\) matrices. 
				
				We denote
				\begin{equation}\label{eq:DefSL2DiagEmbAlgebra}
					\lie{sl}^{\Delta}_2(\R):=\set{(X,\dd\varrho_{k-1}(X)):\,X\in\lie{sl}_2(\R)}\subset \lie{sl}_2(\R)\times \lie{pgl}_k(\R),
				\end{equation}
				where \(\dd\varrho_{k-1}:\lie{sl}_2(\R)\to \lie{pgl}_k(\R)\) is the differential of \eqref{eq:DefOfIrredRep}.
				\begin{proposition}\label{prop:ClassificationOfFInBFg}
					Let \(\boldP(x,y)\in \VV_n^k\) be a linearly independent ordered tuple. Let \(\varphi_{\boldP}:\R\smallsetminus\Sigma\to \SL_2(\R)\times\PGL_k(\R)\) as in \eqref{eq:DefJointWronskianCurve}, where \(\Sigma:=\set{s\in\R:\,\det(W_{\boldP}(s))=0}\). 
					
					Let \(\bdelta_{k-1}\) be the diagonal subgroup \eqref{eq:DefJointDiag}, and let  \(F\subseteq \SL_2(\R)\times \PGL_k(\R)\) be a connected closed reductive subgroup such that \(\bdelta_{k-1}(\R)\subset F\) and such that  there exists an interval \(I\subset\R\smallsetminus \Sigma\) and a \( \bg=(g_1,g_2)\in \SL_2(\R)\times \PGL_k(\R)\) with \(\varphi_{\boldP}(s)\in \bB^- F \bg\) for all \(s\in I\). 
                    
                    Then, there exists a diagonal element \(\balpha_0  \in \bT\subset \SL_2(\R)\times \PGL_k(\R)\)  such that  either \[\f^{\balpha_0} =\lie{sl}^{\Delta}_2(\R),\;\text{ or }\;\lie{sl}^{\Delta}_2(\R)\subset\f^{\balpha_0}=\lie{sl}_2(\R)\times \lie h,\]
                    where \(\lie{h}\) is\stdlist{}. Here, \(\f^{\balpha_0}:=\Ad(\balpha_0)\f\), with \(\f:=\Lie(F)\). 
                    
                    Moreover, \[\varphi_{\boldP}(s)\in \bB^{-} F^{\balpha_0} (\balpha_0\bg)=\bB^- F \bg,\]
                    for all \(s\in\R\smallsetminus\Sigma\), with \(F^{\balpha_0}=\balpha_0F\balpha_0^{-1}\).
				\end{proposition}
				
				\subsection{Proof of Proposition~\ref{prop:ClassificationOfFInBFg}}
				
				We need the following classification result.
				\begin{proposition}\label{prop:AdmissibleLieAlgebra}
					Let \(\lie m\subset \lie{pgl}_k(\R)\) be a reductive subalgebra that contains \(\upE_{k-1}\) and \(\upH_{k-1}\). 
					
					Then, \(\lie m\) is a\stdlist{} subalgebra, which is in particular a split, simple Lie algebra, with a split Cartan subalgebra \(
					\mathfrak{t}_{\lie m}:=\lie m\cap \mathfrak{t}\).
				\end{proposition}
				This classification is  well-known. Since we were not able to locate a self-contained reference, we include a proof relying on Dynkin's classification in the appendix.
				
				We have the following corollary of Propositions~\ref{prop:AdmissibleLieAlgebra}. 
				Let \(\piL:\SL_2(\R)\times \PGL_k(\R) \to \SL_2(\R)\) and \(\piR:\SL_2(\R)\times \PGL_k(\R) \to \PGL_k(\R)\) be the natural maps. 
				\begin{corollary}\label{cor:ClassificationOfLieAlgInProd}
					Let \(\f\subseteq \lie{sl}_2(\R)\times \lie{pgl}_k(\R)\) be a reductive Lie subalgebra containing the elements $(\upE_1, \upE_{k-1})$ and  $(\upH_1, \upH_{k-1})$.
					Then, \(\lie{sl}_2^\Delta(\R)\subset \f\), and moreover, either \(\f=\lie{sl}_2^\Delta(\R)\) or \(\f=\lie{sl}_2(\R)\times \lie h\), where \(\lie h\) is a\stdlist{} algebra.   
				\end{corollary}
				\begin{proof} Denote
					$\f_1=\dd\piL(\f), \mathfrak f_2=\dd\piR(\f).$ Recall that a homomorphic image of a reductive algebra is reductive, so
					both \(\f_1,\,\f_2\) are reductive. Also, \(\upE_1,\upH_1\in f_1\) and \(\upE_{k-1},\upH_{k-1}\in \f_2\). Then, by Proposition~\ref{prop:AdmissibleLieAlgebra}, both \(\f_1,\,\f_2\) are\stdlist{} simple Lie algebras. For \(\f_1\subset \lie{sl}_2(\R)\), the only possibility is \(\f_1=\lie{sl}_2(\R)\). 
					
					Consider \(\lie{I}:=\ker \dd\piL\cap\f\subset \mathbf{0}\times \f_2\). Since \(\f_2\) is simple, either \(\lie{I}=\mathbf{0}\) or \(\lie{I}=\mathbf 0\times \f_2\). In the first case, we get an onto homomorphism \(\lie{sl}_2(\R)=\f_1\to \f_2\). In fact, since \(\lie{sl}_2(\R)\) is simple, the former map is an isomorphism.
					
					If \(\varrho_{k-1}\lie{sl}_2(\R)\subsetneq \f_2\), then \(\f_2\) is a higher rank simple Lie algebra. In this case, \(\lie{I}=\mathbf 0\times \f_2\) since \(\lie{sl}_2(\R)\) and \(\f_2\) cannot be isomorphic.  In particular, \(\f=\lie{sl}_2(\R)\times \f_2\). 
					
					Suppose that \(\f_2=\varrho_{k-1}\lie{sl}_2(\R)\). By Lemma~\ref{lem:fIsSimpleAndContainsPrincSL2}\eqref{enu:UniqAndContainmentOFrhok-1SL2}, there is a unique homomorphism \(\vartheta:\lie{sl}_2(\R)\to\varrho_{k-1}\lie{sl}_2(\R)\) such that \(\vartheta(\upE_1)=\upE_{k-1},\;\vartheta(\upH_1)=\upH_{k-1}.\) So either \(\f=\lie{sl}_2^\Delta(\R)\) or \(\f=\lie{sl}_2(\R)\times\varrho_{k-1}\lie{sl}_2(\R)\).
				\end{proof}
				
				To prove Proposition~\ref{prop:ClassificationOfFInBFg}, we verify that the assumption needed for the classification Proposition~\ref{prop:AdmissibleLieAlgebra} holds in our situation using the following statement.
				
				\begin{proposition}\label{prop:ExistanceOfPrincipalNilpotInf}
					 Let \(W(s),\;s\in\R\smallsetminus \Sigma\) be a Wronskian curve  \eqref{eq:DefOfWronskian}, and let \(\Sigma:=\set{s\in\R:\,\det(W(s))=0}\). Let  \(F\subseteq \SL_2(\R)\times \PGL_k(\R)\) be a connected closed reductive subgroup such that \(\bdelta_{k-1}(\R)\subset F\) and such that the following properties hold. There exists an interval \(I\subset\R\smallsetminus \Sigma\) and a \(\bg\in \SL_2(\R)\times \PGL_k(\R)\) with \((\nu(s),\bkt{W}(s))\in \bB^- F\bg\) for all \(s\in I\). 
                    
                    Then, there exists a diagonal element \(\balpha_0\in \bT\subset \SL_2(\R)\times \PGL_k(\R)\)  such that \((\upE_1,\upE_{k-1}),\) \((\upH_1, \upH_{k-1})\in \Ad(\balpha_0)\f\).
				\end{proposition}
				
				To prove Proposition~\ref{prop:ExistanceOfPrincipalNilpotInf}, we need the following useful elementary lifting lemma in Lie groups.  
				\begin{lemma}\label{lem:BFDecompLemma}
					Let \(G\) be a Lie group, and let \(B,F\subset G\) be two closed subgroups. Suppose that \(\gamma:I\to BF\) is a smooth curve. Then, there exists a subinterval \(J\subset I\) and smooth curves \(b:J\to B\) and \(f:J\to F\) such that \(\gamma(s)=b(s)f(s)\) for all \(s\in J\). 
				\end{lemma}
				\begin{proof}
					We fix \(s_0\in I\), and we assume without loss of generality that  \(\gamma(s_0)=\Id\).
					
					Choose complements \(\lie m\subset\lie b\) and
					\(\lie c\subset\lie g\) such that
					\[
					\lie b=\lie m\oplus(\lie b\cap\lie f),
					\quad
					\lie g=\lie m\oplus\lie c\oplus\lie f = \lie c\oplus(\lie b+\lie f).
					\]
					Then, there exist open neighborhoods of the origin \(\calM\subset \lie m,\,\CC\subset\lie c,\,\calF\subset\f\) such that
					\[
					\calM \times\CC \times \calF \to G,\;\;(X,Z,Y)\longmapsto\exp X\,\exp Z\,\exp Y,
					\]
					is a diffeomorphism onto its image. Consequently,  there is a subinterval \(s_0\in J\subset I\) such that
					\[
					\gamma(s)
					=\exp X(s)\,\exp Z(s)\,\exp Y(s),\,s\in J,
					\]
					with smooth \(X,Z,Y\) taking values in \(\calM,\,\CC,\,\calF\),
					respectively. Since \(\gamma(s)\in BF\) for all \(s\in J\),
					\[
					\exp Z(s)
					=\exp(-X(s))\,\gamma(s)\,\exp(-Y(s))
					\in BF\cap \exp\CC.
					\]
					To complete the proof, it is sufficient to verify that the image  \(Z(J)\) is countable, after possibly shrinking \(\CC\). Indeed, then the connectedness of \(J\) and
					\(Z(s_0)=0\)  imply that \(Z\equiv0\).
					
					Consider \(\mu:B\times F\to G,\;\mu(b,f):=bf\). Since \(\mu\) has  constant rank, for each \((b,f)\in B\times F\), there exists a neighborhood \(W=W_{b,f}\subset B\times F\) such that \(\mu(W)\) is a submanifold of \(G\). Thus, to prove that \(Z(J)\) is countable, it is sufficient to verify that, after possibly shrinking \(\CC\), for each \(bf\in BF\cap \exp\CC\), the submanifolds \(\mu(W)\) and \(\exp\CC\) are transverse, that is
					\begin{equation}\label{eq:transversaility}
						T_{bf}G=T_{bf}(\exp \CC)\oplus T_{bf}(\mu(W)).
					\end{equation}
					Indeed in this situation, after possibly shrinking \(W\) around \((b,f)\), \(\mu(W)\cap\exp\CC\) is a singleton, and countability of \(BF\cap\exp\CC\) follows by second countability of \(B\times F\).
					
					Denote \(g:=bf\) and  \(W:=W_{b,f}\). At \(g=\Id\), \eqref{eq:transversaility} holds by our choice of \(\lie c\), \(\g=\lie c\oplus (\lie b+\lie f)\). 
					By a standard argument using the smooth dependence of the tangent spaces under the \(G\)-right and left translations,  \eqref{eq:transversaility} holds for \(g\) sufficiently close to the identity. Then,  after shrinking \(\CC\), we obtain \eqref{eq:transversaility}.
				\end{proof}

				\begin{proof}[Proof of Proposition~\ref{prop:ExistanceOfPrincipalNilpotInf}]
					It will be convenient to identify the image of \(W(s)\) in \(\PGL_k(\R)\) with the normalized Wronskian \(\widehat{W}(s):=|\det W(s)|^{-1/k}W(s),\;s\in \R\smallsetminus \Sigma\). Also, we identify \(\lie{pgl}_k(\R)\) with the algebra of \(k\times k\) traceless matrices.
					
					Denote \(\varphi(s):=(\nu(s),\widehat{W}(s))\). Let  \(F\subseteq \SL_2(\R)\times \PGL_k(\R)\) be a connected closed reductive subgroup such that \(\bdelta_{k-1}(\R)\subset F\) and such that  there exists an interval \(I\subset\R\smallsetminus \Sigma\) and a \(\bg\in \SL_2(\R)\times \PGL_k(\R)\) with \(\varphi(s)\in (B_1^-\times B_2^-)F\bg\) for all \(s\in I\). 
					
					Denote  \(\varphi'(s):=\frac{\dd\varphi}{\dd s}(s)\), and consider the logarithmic derivative \(\varphi'(s)\varphi(s)^{-1}\in \lie{sl}_2(\R)\times\lie{sl}_k(\R)\). We claim that 
					\begin{equation}\label{eq:LogDerOfvarphi}
						\varphi'(s)\varphi(s)^{-1}=(\upE_1,\upE_{k-1})+Z(s),
						\quad Z(s):=(Z_1(s),Z_2(s))\in\Lie(B_1^-\times B_2^-), \;\forall s\in I.
					\end{equation}
					Indeed, the first \(k-1\) rows of
					\(W'(s)\) are the last \(k-1\) rows of \(W(s)\). Hence, the matrix form of \(W'(s)W(s)^{-1}\) is given by
					\[\begin{pmatrix}
						0&1&0&\dots&0\\
						0&0&1&\ddots&\vdots\\
						\vdots&\vdots&\ddots&\ddots&0\\
						0&0&\dots&0&1\\
						z_1(s)&z_2(s)&\dots&\dots&z_k(s)
					\end{pmatrix}.\]
					The normalization by \((|\det W(s)|)^{-1/k}\) only adds a diagonal scalar term to the logarithmic derivative. Both entries of \(\varphi(s)\) are Wronskian curves, which verifies \eqref{eq:LogDerOfvarphi}. 
					
					Next, denote by \(U_i^-\subset B_i^-\) the lower triangular unipotent subgroups, and denote \(T_i\subset B_i^-\) the diagonal subgroups, so that \(B_i^-=U_i^-T_i\). By Lemma~\ref{lem:BFDecompLemma}, there is a subinterval \(J\subset I\) and smooth curves \(\bu^-:J\to U_1^-\times U_2^-\),  \(\balpha:J\to T_1\times T_2\) and \(\boldf: J\to F\) such that 
					\[\varphi(s)=\bu^-(s)\balpha(s)\boldf(s)\bg,\;s\in J.\]
					As a consequence, the logarithmic derivative is also expressed as
					\begin{equation}\label{eq:LogDerOfBFg}
						\varphi'(s)\varphi(s)^{-1}
						=
						\bb'(s)\bb(s)^{-1}
						+
						\Ad(\bb(s))\bigl(\boldf'(s)\boldf(s)^{-1}\bigr),\;s\in J,
					\end{equation}
					where \(\bb(s)=\bu^-(s)\balpha(s)\). 
					
					Notice that since \(\bdelta_{k-1}(\R)\subset F\), by taking derivative, we have \((\upH_1,\upH_{k-1})\in \f\), where we recall that \(\upH_{j-1}=\diag(j-1,j-3,\dots,-(j-1))\) with \(j\in \N_{\geq 2} \). 
                    
                    We consider the weight-space decompositions of \(\lie{sl}_2(\R)\times\lie{sl}_k(\R)\) and of the induced weight decomposition \(\f\) corresponding to the \(\ad\)-action of \((\upH_1,\upH_{k-1})\). Namely, we write
					\[\boldf'(s)\boldf(s)^{-1}=:\sum_{\bf r}Y_{\bf r}(s),\] 
					where \(\bf r\) varies on the different weights of \(\ad(\upH_1,\upH_{k-1})\) and \(Y_{\bf r}(s)\) are the corresponding weight vectors.  
					To complete the proof, we compare the weight decompositions of \eqref{eq:LogDerOfvarphi} and
					\eqref{eq:LogDerOfBFg}.
					
					In general, for \(j\in \N\), the weights of \(\upH_{j-1}\) are  \(r=2(j-1),2(j-2),\dots,-2(j-1)\), the weight-vectors corresponding to \(r\neq0\) in \(\lie{sl}_j(\R)\) are given by \(\be_{a,b},\;2(b-a)=r,\) and the zero-weights correspond to the diagonal sub-algebra. The diagonal subgroup preserves each \(\upH_{j-1}\) weight spaces. So \(\upE_1,\upE_{k-1}\) are each weight vectors corresponding to weight 2. Also, note the general fact that if \(v\in \SL_j(\R)\) is a lower triangular unipotent, and \(X\in\lie{sl}_j(\R)\) is a weight vector corresponding to weight \(r\), then \(\Ad(v)X=X+L\), where \(L\) is a sum of weight vectors corresponding to lower weight. 
					
					Fix \(s_0\in J\). Since \(\bb'(s_0)\bb(s_0)^{-1}\) has only nonpositive weights, \(\Ad(\bb(s_0))(\boldf'(s_0)\boldf(s_0)^{-1})\) has highest weight 2, with component \((\upE_1,\upE_{k-1})\).
					Diagonal conjugation preserves weights, while lower unipotent conjugation adds only lower-weight terms. Hence \(\boldf'(s_0)\boldf(s_0)^{-1}\) also has highest weight \(2\), and \((\upE_1,\upE_{k-1}) =\Ad(\balpha(s_0))Y_2(s_0)\in\Ad(\balpha(s_0))\f\). Finally, since \((\upH_1,\upH_{k-1})\in \f\) which commutes with the diagonal element \(\balpha(s_0)\), we also have  \((\upH_1,\upH_{k-1})\in \Ad(\balpha(s_0))\f\)
			      \end{proof}
				
				\begin{proof}[Proof of Proposition~\ref{prop:ClassificationOfFInBFg}]
					Write $G=\SL_2(\R)\times\PGL_k(\R)$,
					$B=B_1^-\times B_2^-$, and let
					$\mathbf G,\mathbf B$ denote the corresponding algebraic groups.

                    Let  \(F\subseteq \SL_2(\R)\times \PGL_k(\R)\) be a connected closed reductive subgroup such that \(\bdelta_{k-1}(\R)\subset F\) and such that  there exists an interval \(I\subset\R\smallsetminus \Sigma\) and a \( \bg=(g_1,g_2)\in \SL_2(\R)\times \PGL_k(\R)\) with \(\varphi_{\boldP}(s)\in \bB^- F \bg\) for all \(s\in I\). 
                    
					By Proposition~\ref{prop:ExistanceOfPrincipalNilpotInf} and Corollary~\ref{cor:ClassificationOfLieAlgInProd}, there exists a diagonal element \(\balpha_0  \in \bT\subset G\)  such that  either \[\f^{\balpha_0} =\lie{sl}^{\Delta}_2(\R),\;\text{ or }\;\lie{sl}^{\Delta}_2(\R)\subset\f^{\balpha_0}=\lie{sl}_2(\R)\times \lie h,\]
                    where \(\lie{h}\) is\stdlist{}.  
					
					Let $\FF$ be the connected split algebraic subgroup corresponding to the\stdlist{} subgroup \(F^{\balpha_0}\),
					so that $F^{\balpha_0}=\FF(\R)^\circ$. For these standard embeddings,
					$\mathbf B_{\FF}:=\mathbf B\cap\FF$ is the standard split Borel
					subgroup of $\FF$. Hence
					$\mathbf B_{\FF}\backslash\FF \cong \mathbf B\backslash\mathbf B\FF $
					is a projective Zariski-closed
					subvariety $Y_0$ 
					\cite[Lemma~6.2.2 and Proposition~6.1.2(iii)]{springer1998}.
					
                    We claim that, using the split Bruhat decomposition, the real points
					of $\mathbf B_{\FF}\backslash\FF$ have representatives in the real connected component $F^{\balpha_0}$. Indeed,  root subgroups are connected, and Weyl-group
					representatives are products of root-group elements
					\cite[Theorem~16.1.3, Lemma~8.1.4(i), and
					Theorem~8.2.8(i)]{springer1998}.
					Thus $Y_0(\R)=B\backslash BF^{\balpha_0}$ is Zariski closed in
					$B\backslash G$. In particular, \(BF^{\balpha_0}\) is Zariski closed in \(G\).
					
					Since $\balpha_0$ is diagonal, $\balpha_0\in B$, so that 
					$BF=BF^{\balpha_0}\balpha_0$. This shows that $BF\bg=BF^{\balpha_0}\balpha_0\bg$ is Zariski closed in $G$, for every \(\bg\in G\).
					Since $W(s)$ is polynomial, the map
					$s\mapsto(\nu(s),\bkt{W}(s))$ is regular on
					$\R\smallsetminus\Sigma$. Its inverse image of $BF\bg$ is therefore
					Zariski closed and contains the infinite interval $I$,
					hence equals $\R\smallsetminus\Sigma$. Namely, \[(\nu(s),\bkt{W}(s))\in BF\bg=BF^{\balpha_0}\balpha_0\bg,\;\;\forall s\in \R\smallsetminus\Sigma.\]
				\end{proof}
				
				\subsection{Wronskian curve decompositions and standard basis criteria}\label{subsec:Explicit Conditions for the Obstruction}
				Our goal of this subsection is to prove Proposition~\ref{prop:ExistanceOfStandardBasis}, and to obtain the explicit decompositions having the form \eqref{eq:DesiredDecomp}. Proposition~\ref{prop:ExistanceOfStandardBasis} will be a direct corollary of Lemmata~\ref{lem:WronskCritForSL2}, \ref{lem:CurveDecompAndCondForOrthogonal}, \ref{lem:CurveDecompAndCondForSymplectic}, \ref{lem:CurveDecompAndCondForExecptionalG22} and Lemma~\ref{lem:s parameter to xy}.
                
                Throughout, we fix a tuple of linearly independent homogeneous polynomials \(\boldP(x,y)=(P_1(x,y),\dots,P_k(x,y))\in \VV_n^k\),  we let \(\bkt{W}(s):=\bkt{W_{\boldP}}(s)\) be the image in \(\PGL_k(\R)\) of the Wronskian \eqref{eq:DefOfWronskian} defined on \(\R\smallsetminus\Sigma\), where  \(\Sigma=\set{s\in\R:\,\det\Wp(s)=0}\).  
				
				We will denote  by \(B_2^-\subset \PGL_k(\R)\) the Borel subgroup of lower triangular matrices, and we let \(U_2^-\subset B_2^-\) be the unipotent radical of \(B_2^-\).
				
				For the following lemma, recall the representation \(\varrho_{k-1}:\SL_2(\R)\to \GL_k(\R)\), which is given by the representing matrices of the irreducible representation \eqref{eq:DefOfIrredRep} with respect to the basis \eqref{eq:DefOfOrderedBasisForRinfty}.
				\begin{lemma}\label{lem:WronskCritForSL2}
					Let \(\f=\varrho_{k-1}\lie{sl}_2(\R),\) and let \(F\subset\PGL_k(\R)\) be the associated connected subgroup. It holds that \(\bkt{W}(s)\in B_2^- F\) on an interval \(I\subset \R\smallsetminus \Sigma\) if and only if \(\bkt{W}(s)\in B_2^- F,\;\forall s\in\R\smallsetminus\Sigma\) if and only if there exist homogeneous polynomials \(H(x,y),\,A(x,y),\,B(x,y)\) such that 
					\begin{equation}\label{eq:PolyOrderedBasisPrincSL2}
						P_i(x,y)=H(x,y) \frac{A(x,y)^{k-i}B(x,y)^{i-1}}{(i-1)!},\;i=1,\dots,k.
					\end{equation}
					Moreover, let 
					\(\widetilde \Sigma := \Sigma\cup \set{s\in \R:\;P_1(1,s)= 0}\)
					and let \(\phi(s):=\frac{P_2(1,s)}{P_1(1,s)}=\frac{B(1,s)}{A(1,s)},\;s\in \R\smallsetminus \widetilde \Sigma.\) Then,  \(\R\smallsetminus\tilde \Sigma=I_1\cup\dots \cup I_m\), where  \(I_j,\;j=1,\dots,m\) denote the connected components of \(\R\smallsetminus\tilde\Sigma\), and there exists an analytic curve \(u^-:\R\smallsetminus\tilde\Sigma\to U_2^-\) such that for each \(j\in\set{1,\dots,m},\) 
					\begin{equation}\label{eq:DecompPrincSL2Case}
						\bkt{W}(s)= u^-(s) \bkt{J^\sigma} f(s),\;\; f(s)=\varrho_{k-1}\paren{\begin{psmallmatrix}
								|\phi'(s)|^{-1/2}&0\\0&|\phi'(s)|^{1/2}
							\end{psmallmatrix}\begin{psmallmatrix}
								1&\phi(s)\\0&1
						\end{psmallmatrix}},\;s\in I_j,    
					\end{equation}
					where
					\[J=\diag(1,-1,1,\dots,(-1)^{k-1}),\;\sigma=\begin{cases}
						1&\textup{ if }\phi'(s)<0,\,\forall s\in I_j,\\
						0&\textup{ if }\phi'(s)>0,\,\forall s\in I_j.
					\end{cases}\]
				\end{lemma}
				\begin{proof}
					Suppose first that \(\bkt W(s)\in B_2^- F\) for all \(s\in I\). Thus, we may write
					\[\bkt{W}(s)=\R^\times b(s)\varrho_{k-1}(M(s)),\;s\in I.\]
					for \(M(s)=\begin{psmallmatrix} \alpha(s) & \beta(s) \\ \gamma(s) & \delta(s) \end{psmallmatrix} \in \SL_2(\R)\), where \(b(s)\in B_2^-\). Observe that the first row of \(\varrho_{k-1}(M(s))\) is of the form
					\[\paren{\alpha(s)^{k-1}, \frac{\alpha(s)^{k-2}\beta(s)}{1!}, \dots, \frac{\beta(s)^{k-1}}{(k-1)!}}.\]
					Indeed, \(\varrho_{k-1}\) is the matrix representation of the irreducible \(\SL_2(\R)\) action \eqref{eq:DefOfIrredRep} with respect to the ordered basis \(e_i(x,y):=\frac{x^{k-i}y^{i-1}}{(i-1)!},\;i=1,\dots,k,\) so that
					\[M(s)\cdot e_i(x,y)=e_i(\alpha(s)x+\gamma(s)y,\beta(s)x+\delta(s)y)=\frac{\alpha(s)^{k-i}\beta^{i-1}(s)}{(i-1)!}e_1(x,y)+\cdots.\]
					
					Thus the first row of \(W(s)\) is expressed as
					$$\paren{P_1(1,s), \dots, P_k(1,s)} = c(s) b_{1,1}(s)\paren{\alpha(s)^{k-1}, \textstyle \frac{\alpha(s)^{k-2}\beta(s)}{1!}, \dots, \frac{\beta(s)^{k-1}}{(k-1)!}},\;c(s)\in\R^\times.$$
					We note that \(b_{1,1}(s)\), \(\alpha(s)\) and \(\beta(s)\) are nonzero for all but finitely many \(s\in I\), since otherwise \(P_1(1,s)\) or \(P_2(1,s)\) is the zero polynomial. That in turn forces a column of \(W(s)\) to be zero, a contradiction to \(\det W(s)\neq 0\).
					
					The former expression implies that  \(\frac{P_{i}(1,s)}{P_{i+1}(1,s)}=\frac{a(s)}{b(s)}i\) for all \(i=1,\dots,k-1\) such that \(a(s),b(s)\) are polynomials with \(\gcd(a,b)=1\). In particular, \(P_{i}(1,s)=P_1(1,s)\frac{b^{i-1}(s)}{a^{i-1}(s)(i-1)!}\). Since the polynomials \(a(s),b(s)\) are coprime,  \(a^{k-1}(s)\) divides \(P_1(1,s)\). Then, \(h(s)=\frac{P_1(1,s)}{a(s)^{k-1}}\) is a polynomial, and \(P_i(1,s)= h(s) \frac{a(s)^{k-i}b(s)^{i-1}}{(i-1)!}\). 
					
					Let \(\kappa_1,\kappa_2,\kappa_3\) be the degrees of \(h,a,b\) correspondingly, and note that
					\begin{align*}
						n \geq \deg(P_1(1,s)) = \kappa_1+(k-1)\kappa_2,\;n \geq \deg(P_k(1,s)) = \kappa_1+(k-1)\kappa_3.
					\end{align*}
					Let \( d :=\max \{\kappa_2, \kappa_3\}\), and let 
					\[H(x,y):=x^{n-(k-1)d}h(y/x),\;\;A(x,y):=x^{d}a(y/x),\;B(x,y):=x^{d}b(y/x).\]
					In general, for every polynomial \(F(s)\) such that \(M\geq\deg(F)\), we have that \(x^M F(y/x)\) is a homogeneous polynomial. Then, \(H(x,y),\,A(x,y),\,B(x,y)\) are homogeneous polynomials. Using $P_i(x,y)=x^{n}P_i(1,\frac yx),$  we have
					\begin{align*}
						P_i(x,y)=x^{n}h\left(\frac{y}{x}\right)\frac{a(y/x)^{k-i}b(y/x)^{i-1}}{(i-1)!}
						=H(x,y)\frac{A(x,y)^{k-i}B(x,y)^{i-1}}{(i-1)!}.
					\end{align*}
					
					Conversely, assume  \eqref{eq:PolyOrderedBasisPrincSL2}. Since \(P_1(1,s)\) is a polynomial, it vanishes only for finitely many \(s\). Then, with \(\widetilde\Sigma = \Sigma \cup \set{s\in \R: \,P_1(1,s)=0}\) and \(\phi(s)=B(1,s)/A(1,s)\) for \(s\in\R\smallsetminus\tilde\Sigma\), we may write
					\[P_i(1,s)=H(1,s)\frac{A(1,s)^{k-i}B(1,s)^{i-1}}{(i-1)!}=h(s)\frac{\phi(s)^{i-1}}{(i-1)!},\;\forall s\in \R\smallsetminus\tilde\Sigma,\;i=1,\dots,k,\]
					where \(h(s)=H(1,s)A(1,s)^{k-1}\). 
					
					Denote
					\[v(t):=\left(1,t,\frac{t^2}{2!},\ldots,\frac{t^{k-1}}{(k-1)!}\right), \quad p(s):=h(s)v(\phi(s))\in\R^k.\]
					By induction, for each \(i\geq 0\) there exist smooth functions \(a_{ir}\), \(0\leq r\leq i\), such that
					\[p^{(i)}(s)=\sum_{r=0}^i a_{ir}(s)\,v^{(r)}(\phi(s)), \; a_{ii}(s)=h(s)(\phi'(s))^i.\]
					
					Recall \(\varrho_{k-1}\begin{psmallmatrix}
						1&s\\0&1 
					\end{psmallmatrix}=(v_j^{(i-1)}(s))_{i,j=1,\dots,k}\); see \eqref{eq:DefOfrhonum}. Consequently,  we get that
					\begin{equation}\label{eq:InProofDecompWPrincSl2}
						W(s)=b(s)\varrho_{k-1}\begin{psmallmatrix}
							1&\phi(s)\\0&1
						\end{psmallmatrix},\, s\in \R\smallsetminus\tilde \Sigma 
					\end{equation}
					where \(b(s)\) is lower triangular. In particular, this shows that \((\nu(s),\bkt W(s))\in \bB^-(\SL_2(\R)\times F)\) on an interval, and by Proposition~\ref{prop:ClassificationOfFInBFg}, we get that \(\bkt W(s)\in B_2^- F\) for all \(s\in \R\smallsetminus\Sigma\). 
					
					Finally, we show \eqref{eq:DecompPrincSL2Case}. By \eqref{eq:InProofDecompWPrincSl2}, the entries on the diagonal of \(b(s)\) are
					\[(b_{1,1}(s),\dots,b_{k,k}(s))=\bigl(h,\;h\phi',\;h(\phi')^2,\ldots,h(\phi')^{k-1}\bigr).\]
                    Thus, we may write $b_{i,i}(s)=h(s)(-1)^{\sigma (i-1)} |\phi'(s)|^{i-1}$ and
                    \[
                    b(s)=
                    h(s)|\phi'(s)|^{(k-1)/2}\,
                    u^-(s)J^\sigma
                    \varrho_{k-1}
                    \begin{pmatrix}
                    |\phi'(s)|^{-1/2}&0\\
                    0&|\phi'(s)|^{1/2}
                    \end{pmatrix}.
                    \]
					Consequently, we obtain \eqref{eq:DecompPrincSL2Case} after multiplying \(W(s)\) by the scalar matrix  \(\frac{1}{h(s)|\phi'(s)|^{(k-1)/2}}\Id\).
				\end{proof}
				
                \begin{corollary}\label{cor:CurveDecompAndCondFordiagonallyEmbeddedPrincipalSL2}
					Let \(\f=\lie{sl}^{\Delta}_2(\R)\subset\lie{sl}_2(\R)\times\lie{pgl}_k(\R)\) be as in \eqref{eq:DefSL2DiagEmbAlgebra}, and denote the associated connected subgroup by \(F:= \{ ( g, [\varrho_{k-1}(g)]): g \in SL_2(\R)\}\). Then the following are equivalent.
                    \begin{enumerate}
                        \item 
                   \((\nu(s),\bkt{W}(s))\in (B_1^-\times B_2^-) F \) on an interval \(I\subset \R\smallsetminus \Sigma,\) 
                   \item \( (\nu(s),\bkt{W}(s))\in (B_1^-\times B_2^-) F,\;\forall s\in \R\smallsetminus\Sigma,\)
                   \item there exists a homogeneous polynomial \(H(x,y)\) of degree \(n-(k-1)\) such that
		\begin{equation}\label{eq:PolyOrderedBasisDiag}
						P_i(x,y)=H(x,y)\frac{x^{k-i}y^{i-1}}{(i-1)!},\;i=1,\dots,k.
					\end{equation}
                     \end{enumerate}
					Moreover, 
					\[(\nu(s),\bkt{W}(s))=(\nu(s),\R^\times u^-(s)\varrho_{k-1}(\nu(s))),\;\forall s\in\R\smallsetminus\Sigma.\]
				\end{corollary}
				\begin{proof}
					Recall the representation \(\varrho_{k-1}:\SL_2(\R)\to \GL_k(\R)\) given by the representing matrices of the irreducible representation \eqref{eq:DefOfIrredRep} with respect to the basis \eqref{eq:DefOfOrderedBasisForRinfty}. If \((\nu(s),\bkt{W}(s))\in (B_1^-\times B_2^-)F\) on \(I\), where \(\f=\lie{sl}_2^{\Delta}(\R),\)  then 
					\[(\nu(s),\bkt{W}(s))=(b_1(s)f(s),\R^\times b_2(s)\varrho_{k-1}(f(s))),\quad \forall s\in I\]
					for $b_1(s)\in B_1^-, b_2(s) \in B_2^-$ and $f(s) \in \SL_2(\R).$ The first coordinate gives $f(s)=b_1(s)^{-1}\nu(s),$ thus $\bkt{W}(s)=\R^\times b_2(s)\varrho_{k-1}(b_1(s)^{-1})\varrho_{k-1}(\nu(s)).$ Since $b_2(s)\varrho_{k-1}(b_1(s)^{-1}) \in B_2^-,$ in view of the first row of \(\R^\times b_2(s)\varrho_{k-1}(b_1(s)^{-1})\varrho_{k-1}(\nu(s))\), we get that the first row entries of \(W(s)\) are of the form
					\[P_i(1,s)=h(s)\frac{s^{i-1}}{(i-1)!},\;\forall s\in I,\;i=1,\dots,k.\]
					Therefore
					\[P_i(x,y)=x^n P_i(1,y/x)=x^{n-k+1}h(y/x)\frac{x^{k-1}(y/x)^{i-1}}{(i-1)!}=x^{n-k+1}h(y/x)\frac{x^{k-i}y^{i-1}}{(i-1)!}.\]
					Conversely, by applying Lemma~\ref{lem:WronskCritForSL2} to the equation \eqref{eq:PolyOrderedBasisDiag}, we have $\phi(s)=s$ and thus
					\[\bkt{W}(s)=\R^\times u^-(s)\varrho_{k-1}(\nu(s)),\quad \forall s \in I.\]
					
					Finally, by Proposition~\ref{prop:ClassificationOfFInBFg},  \((\nu(s),\bkt{W}(s))\in (B_1^-\times B_2^-) F\) on an interval \(I\subset \R\smallsetminus \Sigma\) if and only if \((\nu(s),\bkt{W}(s))\in (B_1^-\times B_2^-) F,\;\forall s\in \R\smallsetminus\Sigma\).
				\end{proof}
				
				We proceed to the situation that \(\bkt{W}(s)\in B_2^-F\), with \(\f=\lie{so}_{m,m+1}(\R)\). 
				We need the following lemma, which gives a polynomial criterion on \(g\in \SO_{m,m+1}(\R)\) determining if \(g\) belongs to identity component.
				\begin{lemma}\label{lem:ConnCompsofSOmm}
					Let \(J=\diag(-1,\Id_{2m-1},-1)\in \SO_{m,m+1}(\R)\). We have \[\SO_{m,m+1}(\R)=\SO_{m,m+1}^{\circ}(\R)\sqcup J\SO_{m,m+1}^{\circ}(\R),\] where \(\SO_{m,m+1}^{\circ}(\R)\) is the identity component. Moreover, for 
					\begin{equation}\label{eq:DefPhiSO}
						\phi_{\lie{so}_{m,m+1}}(g):=\det\left(
						g_{ij}+(-1)^{m+j}g_{i,\,2m+2-j}
						\right)_{1\leq i,j\leq m},
					\end{equation}
					the element \(g\in \SO_{m,m+1}(\R)\) belongs to the identity component if and only if \(\phi(g)>0\).
				\end{lemma}
				\begin{proof}
					We recall that \(\SO_{m,m+1}(\R)\)  has two connected components, see \cite[Proposition~1.145]{Knapp2002}.
					
					We proceed to describe the criterion determining to which component a given element belongs to. Recall that \(\SO_{m,m+1}(\R)=\set{g\in \SL_{2m+1}(\R):\,q(\bx g,\by g)=q(\bx,\by)},\) where \(q(\bx,\by):=\sum_{i=1}^{2m+1}(-1)^{i+1}x_i y_{2m+2-i}\). Since \(q(\be_i,\be_j)=0\) for all \(1\leq i,j\leq m\), the first \(m\)-rows of \(g\in \SO_{m,m+1}(\R)\) span a totally isotropic subspace. We note that totally isotropic subspaces for a quadratic form of signature \(m,m+1\) of dimension \(m\) are maximal totally isotropic subspaces, and we now define an orientation on those subspaces. On a connected component, the orientation of the first \(m\)-rows cannot change. 
					
					Consider the following two complementary subspaces of \(\R^{2m+1}\),
					\[V_1:=\gen{\be_j:1\leq j\leq m}_\R,\;\;V_2:=\gen{\be_{m+1},\,\be_j-(-1)^{m+j}\be_{2m+2-j}:1\leq j\leq m}_\R.
					\]
					Notice that \(q\mid_{V_2}\) is definite. Let \(\Pi:\R^{2m+1}\to V_1\) be the projection onto \(V_1\) with kernel \(V_2\). Given a totally isotropic subspace \(U\subset\R^{2m+1}\), notice that \(\Pi\mid_U: U\to V_1\) is injective. Indeed, for \(\boldsymbol{u}\in U\), we have \(q(\boldsymbol{u})=0\), and if  in addition \(\boldsymbol{u}\in\ker(\Pi\mid_U)=U\cap V_2\subset V_2\), then \(\boldsymbol{u}=0\) since \(V_2\) is a \(q\)-definite subspace. In particular, for (maximal) totally isotropic subspace \(U\)  such that \(\dim U=m\), we get that \(\Pi\mid_U\) is an isomorphism. Thus, by fixing an orientation on \(V_1\), we obtain an orientation on each \(m\)-dimensional totally isotropic subspace. 
					
					The projection \(\Pi\colon\mathbb{R}^{2m+1}\to V_1\) with kernel \(V_2\) is given concretely by
					\[\Pi(x)=\sum_{j=1}^m
					\left(x_j+(-1)^{m+j}x_{2m+2-j}\right)\be_j.
					\]
					Therefore, if \(g=(g_{ij})\in \SO_{m,m+1}(\R)\), then the projections of its first \(m\)-rows are positively oriented if and only if
					\[
					\phi(g):=\det\left(
					g_{ij}+(-1)^{m+j}g_{i,\,2m+2-j}
					\right)_{1\leq i,j\leq m}>0.
					\]
					In particular, since \(\phi(I_{2m+1})>0\),  \(g\) is in the identity component if and only if \(\phi(g)>0\). Finally, for \(J=\diag(-1,I_{2m-1},-1)\in \SO_{m,m+1}(\R)\), we have \(\phi(J)<0\), so that the non-identity component is \(J\SO_{m,m+1}^\circ(\R).\)
				\end{proof}
				
				\begin{lemma}\label{lem:CurveDecompAndCondForOrthogonal}
					Let \(\f=\lie{so}_{m,m+1}(\R)\), and let
					\(F\subset\PGL_{2m+1}(\R)\) be the associated connected subgroup.
					Let \(I\subset\R\smallsetminus\Sigma\) be a nonempty open interval.
					Then \(\bkt{W}(s)\in B_2^-F\) for all \(s\in I\) if and only if
					\(\bkt{W}(s)\in B_2^-F\) for all \(s\in\R\smallsetminus\Sigma\),
					if and only if, for all \(s\in I\), the \(m\)-dimensional subspace
					\[
					\gen{
						\paren{P_1^{(i-1)}(1,s),\dots,P_{2m+1}^{(i-1)}(1,s)}
						:i=1,\dots,m
					}_\R\subset\R^{2m+1}
					\]
					is totally isotropic for the form \(q(\bx,\by):=\sum_{i=1}^{2m+1}(-1)^{i+1}x_i y_{2m+2-i}.\)
					
					Moreover, suppose these equivalent conditions hold, and let
					\(I_1\sqcup\dots\sqcup I_L=\R\smallsetminus\Sigma\) be the
					connected components. Then for each \(1\leq i\leq L\), there
					exist analytic curves \(u^-:I_i\to U_2^-\) and \(f:I_i\to F\)
					such that
					\[
					\bkt{W}(s)=\R^\times u^-(s)J^{\sigma}f(s), \quad
					\sigma=
					\begin{cases}
						0,&(\det W(s))^m\phi_{\lie{so}_{m,m+1}}(W(s))>0,
						\quad\forall s\in I_i,\\
						1,&(\det W(s))^m\phi_{\lie{so}_{m,m+1}}(W(s))<0,
						\quad\forall s\in I_i.
					\end{cases}
					\]
				\end{lemma}
               
				\begin{proof}
					We have
					\(\bkt{W}(s)\in B_2^-F\) for all \(s\in\R\smallsetminus\Sigma\)
					if and only if \(\bkt{W}(s)\in B_2^-F\) for all \(s\in I\) by Proposition~\ref{prop:ClassificationOfFInBFg} applied to \(s\mapsto (\nu(s),\bkt W(s))\in \bB^-(\SL_2(\R)\times F)\).
					
					Consider the subspace \(E_m:=\gen{\be_1,\dots,\be_m}_\R\),
					which is totally isotropic for \(q\).
					Since \(F\) preserves \(\R^\times q\) and \(B_2^-\) preserves
					\(E_m\), then \(\bkt{W}(s)\in B_2^-F\) implies that
					\(E_mW(s)\) is totally isotropic. This proves the only if part.
					
					Conversely, assume that \(E_mW(s)\) is totally isotropic for
					all \(s\in I\). Let \(Q\) be the symmetric matrix defining \(q\), namely
					\[
					Q_{i,j}=(-1)^{i+1}\delta_{i+j,\,2m+2},\;i,j=1,\dots,2m+1
					\]
					We denote \(A(s):=W(s)QW(s)^\top\). Since \(E_m W(s)\) is totally isotropic, 
					\[A(s)_{i,j}=q(\be_i W(s),\be_j W(s))=0 ,\]
					for all \(i,j\leq m\).
					Since \(\frac{\dd}{\dd s}\be_i W(s)=\be_{i+1}W(s)\) for \(i<2m+1,\) the derivative product rule gives
					\[
					\frac{\dd}{\dd s}A(s)_{i,j}
					=A(s)_{i+1,j}+A(s)_{i,j+1},
					\; 1\leq i,j<2m+1.
					\]
					An induction along the anti-diagonals, using this identity
					and the symmetry of \(A(s)\), yields
					\[
					A(s)_{i,j}=0\quad(i+j\leq2m+1),\;\;
					A(s)_{i,\,2m+2-i}=(-1)^{i+1}\lambda(s),
					\]
					where \(\lambda(s):=A(s)_{1,\,2m+1}\).
					
					It follows that \(A(s)Q\) is lower triangular,
					with every diagonal entry equal to \(\lambda(s)\).
					Taking determinants gives
					\[
					\lambda(s)^{2m+1}
					=\det(A(s)Q)
					=\det(W(s))^2\det(Q)^2>0,
					\; s\in I.
					\]
					Hence \(\lambda(s)>0\), and
					\(\widetilde u(s):=\lambda(s)^{-1}A(s)Q\)
					is lower-triangular unipotent.
					Its unique lower-triangular unipotent square root \(u(s):=\exp\paren{\frac12\log\widetilde u(s)}\) is analytic on \(I\).
					
					Since \(Q\widetilde u(s)^\top Q=\widetilde u(s)\),
					the matrix \(Qu(s)^\top Q\) is also a lower-triangular
					unipotent square root of \(\widetilde u(s)\).
					By uniqueness, \(Qu(s)^\top Q=u(s)\), and consequently
					\[
					u(s)Qu(s)^\top
					=u(s)^2Q
					=\lambda(s)^{-1}A(s)Q^2=\lambda(s)^{-1}A(s).
					\]
					Now set
					\[
					h(s):=\sgn(\det W(s))
					\paren{\sqrt{\lambda(s)}u(s)}^{-1}W(s).
					\]
					The preceding identity gives \(h(s)Qh(s)^\top=Q\), while
					the determinant identity gives \(\det h(s)=1\).
					Thus \(h(s)\in\SO_{m,m+1}(\R)\) is analytic on \(I\), and
					\[
					\bkt{W}(s)=\R^\times u(s)h(s).
					\]
					
					Finally, Lemma~\ref{lem:ConnCompsofSOmm} completes the proof. In fact,  \[
					\phi_{\lie{so}_{m,m+1}}(h(s))
					=\sgn(\det W(s))^m\lambda(s)^{-m/2}
					\phi_{\lie{so}_{m,m+1}}(W(s)),
					\]
					its sign is thus the sign of
					\((\det W(s))^m\phi_{\lie{so}_{m,m+1}}(W(s))\). 
				\end{proof}

				We have the following similar statement for \(\f=\lie{sp}_{2m}(\R)\).  The proof is analogous to the proof of Lemma~\ref{lem:CurveDecompAndCondForOrthogonal}, hence omitted.
				
				\begin{lemma}\label{lem:CurveDecompAndCondForSymplectic}
					Let \(\f=\lie{sp}_{2m}(\R),\) and let \(F\subset\PGL_{2m}(\R)\) be the associated connected subgroup. It holds that \(\bkt{W}(s)\in B_2^- F\) on an interval \(I\subset \R\smallsetminus \Sigma\) if and only if \(\bkt{W}(s)\in B_2^- F,\;\forall s\in\R\smallsetminus\Sigma\) if and only if  for all \(s\in I\), the \(m\)-dimensional subspace
					\[\gen{\paren{P^{(i-1)}_1(1,s),\dots,P^{(i-1)}_{2m}(1,s)}:\,i=1,\dots,m}_\R\subset \R^{2m},\]
					is totally isotropic for the  form \(q(\bx,\by)=\sum_{i=1}^{2m}(-1)^{i+1}x_i y_{2m+1
						-i}.\)
					Moreover, let \(I_1\sqcup \dots \sqcup I_L = \R\smallsetminus\Sigma\) be the connected components of \(\R\smallsetminus\Sigma\). Then for each \(1\leq i\leq L\), there exist analytic curves \(u^-:I_i\to U_2^-\) and \(f:I_i\to F\) such that \[\bkt{W}(s)=\R^\times u^-(s)J^\sigma f(s),\,\forall s\in I_i,\]
					where \(\sigma \in \set{0,1}\) and \(J=\diag(I_m,-I_m).\)
				\end{lemma}
				\begin{remark}
					There are similar polynomial conditions determining \(\sigma\), but we omit them since it is not required for our purposes, see Lemma~\ref{lem:6.6}.
				\end{remark}

				\begin{lemma}\label{lem:CurveDecompAndCondForExecptionalG22}
					Let \(\f=\g_{2(2)}\), and denote by \(F\subseteq \PGL_7(\R)\) the associated connected subgroup. Let \(\Phi=dx_{147}-2dx_{156}-dx_{237}+dx_{246}-dx_{345}\in \bigwedge^3(\R^7)^*\) be the alternating \(3\)-form defining \(\f\) in \eqref{eq:G22form}. It holds that $\bkt W(s) \in B_2^- F$ on an interval $I\subset \R \smallsetminus \Sigma$ if and only if the 1-form $\Phi(P(1,s), P'(1,s),-)$ is the zero form for all \(s\in I\). Moreover,
					\[\bkt{W}(s)=\R^\times u^-(s)f(s)\]
					where $u^-(s)\in U^-_2$ and $f(s) \in F$ are analytic.
				\end{lemma}
				\begin{proof}
					The class \([\Phi]\) is $F$-invariant, that is, for all \(\bv,\bu,\bw\in\R^7\), 
					\[\R^\times \Phi(\bv,\bu,\bw)= \R^\times \Phi(\bv g,\bu g, \bw g),\;\;\forall g\in F.\]
					Now, for \(\bu,\bv\in \gen{\be_1,\be_2}_{\R}\), the one form \(\Phi(\bu,\bv,-)\) is zero. Notice also that \(\gen{\be_1,\be_2}_{\R}B_2^-=\gen{\be_1,\be_2}_{\R}.\) Therefore, if \(\bkt{W}(s)\in B_2^-F\),  then $\Phi(\be_1W(s),\be_2W(s),-)$ is  trivial.
					
					Conversely, let \(I\subseteq\R\smallsetminus\Sigma\) be an interval, and suppose that $\Phi(\be_1W(s), \be_2W(s),-),\;\forall s\in I$ is the zero form. Denote 
					\[R_{ijk}(s):=\Phi(\be_iW(s), \be_jW(s), \be_kW(s))\]
					for $i,j,k\in\{1,\dots,7\}.$ We first show that $R_{ijk}(s)=0,\;\forall s\in I$ for every $i+j+k<12$. 
					
					Since \(\Phi\) is alternating, it suffices to consider \(i<j<k\).
					For the following, recall \(\frac{\dd}{\dd s}\be_i W(s)=\be_{i+1} W(s),\,i=1,\dots,6,\) which, together with the derivative product rule, gives \[\frac{\dd}{\dd s}R_{ijk}(s)=R_{i+1\,jk}(s)+R_{i j+1\,k}(s)+R_{ijk+1}(s),\;\forall i,j,k<7.\] 
					For \(i+j+k\leq 8\), the only
					such triples are \((1,2,3),\,(1,2,4),\, (1,2,5),\,(1,3,4).\)  We already have \(R_{12k}=0\) on \(I\), for all \(k\) since $\Phi(\be_1W(s), \be_2W(s),-),\;\forall s\in I$ is the zero form. Taking derivatives of \(R_{12k}\)  gives
					\(R_{13k}=0\).  
					
					We proceed by increasing \(i+j+k\). Differentiating \(R_{134},R_{234},R_{135}\),
					in this order, gives \(R_{234}=R_{235}=R_{145}=0.\) 
					Finally, differentiating \(R_{136},R_{145},R_{235}\) gives
					\[
					R_{236}+R_{146}=0,\;
					R_{245}+R_{146}=0,\;
					R_{245}+R_{236}=0,
					\]
					so \(R_{146}=R_{236}=R_{245}=0\). This proves that
					\(R_{ijk}=0\) on \(I\) whenever \(i+j+k<12\).
					
					Next, we determine the coefficients with \(i+j+k=12\).
					Differentiating \(R_{137}=R_{146}=R_{236}=R_{245}=0\), and using
					also \(\Phi(\be_1W,\be_3W,-)=0\), gives
					\(R_{237}+R_{147}=0\),
					\(R_{246}+R_{156}+R_{147}=0\),
					\(R_{246}+R_{237}=0\), and
					\(R_{345}+R_{246}=0\).
					Thus, setting \(r(s):=R_{147}(s)\), we obtain
					\((R_{147},R_{156},R_{237},R_{246},R_{345})
					=(r,-2r,-r,r,-r)\) on \(I\).
					
                    We claim that $r(s)$ is a polynomial that does not vanish on $I.$ For a vector $\mathbf a =(a_1,\dots, a_7),$ consider the interior product
						\begin{align*}
							\iota_\mathbf{a} \Phi &= \Phi(\mathbf a,-,-) = \;a_1dx_{47} - 2a_1 dx_{56} -a_2dx_{37} +a_2dx_{46} + a_3dx_{27} -a_3dx_{45}-a_4dx_{17}\\ & -a_4dx_{26} 
							+a_4dx_{35} +2a_5dx_{16} -a_5dx_{34}  -2a_6dx_{15} +a_6dx_{24} +a_7dx_{14} -a_7dx_{23}.
						\end{align*}
						If $\iota_\mathbf a \Phi \wedge \iota_\mathbf a \Phi = 0,$ by comparing the coefficients of $dx_{4567},dx_{3467},dx_{2457},dx_{1346},dx_{1245},dx_{1234},$ we obtain
						\[a_1=a_2=a_3=a_5=a_6=a_7=0.\]  
						The coefficient of $dx_{1267}$ is $2a_4^2-4a_3a_5=2a_4^2,$ hence $a_4=0$ and $\mathbf a =\mathbf 0.$ However, if $r(s_0)=0,$ then $\Phi(\be_1W(s_0),-,-)=R_{157}\theta^5\wedge \theta^7 +R_{167}\theta^6 \wedge \theta^7$ in the dual basis $\theta^1,\dots,\theta^7$ of $\be_1W(s_0),\dots, \be_7W(s_0).$ Hence 
						\[\iota_{\be_1 W(s_0)} \Phi \wedge \iota_{\be_1 W(s_0)} \Phi =0,\]
						so the above observation implies that $\be_1W(s_0)=\mathbf 0,$ contradicting the non-singularity of $W(s_0)$. This proves the claim. 
                        
                        Therefore, by multiplying $\frac1{r(s)^{1/3}}$ on the left, the modified $$R^{[1]}_{ijk}:=\Phi(\be_i \frac1{r(s)^{1/3}}W(s), \be_j\frac1{r(s)^{1/3}}W(s), \be_k\frac1{r(s)^{1/3}}W(s))$$ 
                        have the values
						$(R^{[1]}_{147}, R^{[1]}_{156}, R^{[1]}_{237}, R^{[1]}_{246}, R^{[1]}_{345})=(1,-2,-1,1,-1).$
					
					Finally, we  show that by multiplying the elements of $U^-_2$ on the left, which act as the operation \(\be_i \mapsto \be_i + c(s)\be_j\) for \(i>j\), we can make $R_{ijk}=0$ for every $i+j+k=13,14,\dots,18.$ By the operation $u^{[1]}(s):\be_7 \mapsto \be_7+ \frac12R^{[1]}_{157} \be_6,$ the new $$R^{[2]}_{ijk}:=\Phi(\be_iu^{[1]}(s)\frac{W(s)}{r(s)^{1/3}},\be_ju^{[1]}(s)\frac{W(s)}{r(s)^{1/3}},\be_ku^{[1]}(s)\frac{W(s)}{r(s)^{1/3}})$$ 
                    becomes
					\[\textstyle R^{[2]}_{157}=\Phi\left(\be_1 \frac{W(s)}{r(s)^{1/3}}, \be_5\frac{W(s)}{r(s)^{1/3}},(\be_7+\frac{R^{[1]} _{157}(s)}{2}\be_6)\frac{W(s)}{r(s)^{1/3}}\right) = R^{[1]}_{157}(s)+\frac{R^{[1]}_{157}(s)}{2}R^{[1]}_{156}(s)=0,\]
					while $R^{[2]}_{ijk}=R^{[1]}_{ijk}$ for $i+j+k\leq 12.$ By repeating this process, the operations 
					$$\textstyle u^{[2]}(s):\be_2 \mapsto \be_2+ \frac12R^{[2]}_{256} \be_1,\ u^{[3]}(s):\be_4 \mapsto \be_4+ R^{[3]}_{247} \be_3,\ u^{[4]}(s): \be_3 \mapsto \be_3+ -R^{[4]}_{346} \be_2$$ in order give $R^{[5]}_{ijk}=0$ for $i+j+k=13,$ while they do not affect the previous terms. Next, by the operations 
					$$\textstyle u^{[5]}(s): \be_7 \mapsto \be_7-\frac12 R^{[5]}_{167} \be_5,\ u^{[6]}(s): \be_5 \mapsto \be_5 + R^{[6]}_{257}\be_3,\dots, u^{[15]}(s):\be_7 \mapsto \be_7+\frac12R^{[15]}_{567}\be_1$$
					in order, we can conclude that $R^{[16]}_{ijk}=0$ for $i+j+k\neq 12,$ and $(R^{[16]}_{147}, R^{[16]}_{156}, R^{[16]}_{237}, R^{[16]}_{246}, R^{[16]}_{345})=(1,-2,-1,1,-1),$ where
					\[\textstyle R^{[16]}_{ijk}=\Phi\left(\be_iu^{[15]}(s)\cdots u^{[1]}(s) \frac{W(s)}{r(s)^{1/3}}, \be_ju^{[15]}(s)\cdots u^{[1]}(s) \frac{W(s)}{r(s)^{1/3}}, \be_ku^{[15]}(s)\cdots u^{[1]}(s) \frac{W(s)}{r(s)^{1/3}}\right).\]
					By denoting $u^-(s):=\left(u^{[15]}(s)\cdots u^{[1]}(s)\right)^{-1}$ and $f(s)=\frac{1}{r(s)^{1/3}}(u^-(s))^{-1}W(s),$ we conclude that 
					$$\Phi(\be_if(s), \be_jf(s), \be_kf(s))=\Phi(\be_i,\be_j,\be_k)$$
					for every $i,j,k \in \{1,\dots,7\}$, which completes the proof.
				\end{proof}

To conclude Proposition~\ref{prop:ExistanceOfStandardBasis}, we have the following lemma.
\begin{lemma}\label{lem:s parameter to xy}
    Let \(\boldP(x,y)=(P_1(x,y),\dots,P_k(x,y))\in \VV_n^k\) denote a linearly independent tuple, and let \(I\subset \R\) be an open interval.
    \begin{enumerate}
        \item  Let \(q: \R^k\times\R^k\to\R\) be a bilinear form. Then, 
            \[q (\partial_y^i\boldP(x,y),\partial_y^j\boldP(x,y))=0,\;\forall 0\leq i,j< \floor{k/2},\]
			if and only if, the  subspace
					\[
					\gen{
						\paren{P_1^{(i-1)}(1,s),\dots,P_{k}^{(i-1)}(1,s)}
						:i=1,\dots,\floor{k/2}
					}_\R\subset\R^{k},
					\] 
					is totally isotropic for \(q\) for all \(s\in I\).
        \item  Take \(k=7\). Let \(\Phi\in \bigwedge^3(\R^7)^*\) be an alternating \(3\)-form. Then,
        \[\Phi(\boldP(x,y), \partial_y\boldP(x,y), \mathbf{v})=0,\;\forall\mathbf{v}\in\R^7,\]
        if and only if the \(1\)-form $\Phi(P(1,s), P'(1,s),-)$ is the zero form for all \(s\in I\).
    \end{enumerate}
\end{lemma}
\begin{proof}
For each indicated \(i,j\) in (1) or fixed \(\mathbf v\) in (2), set
\[
H(x,y)=
\begin{cases}
q(\partial_y^i\boldP(x,y),\partial_y^j\boldP(x,y)) & \text{in (1)},\\
\Phi(\boldP(x,y),\partial_y\boldP(x,y),\mathbf v) & \text{in (2)}.
\end{cases}
\]
In (1), \(H(1,s)\) is the pairing of two basis elements, and
bilinearity extends the vanishing to all
pairs of vectors in the spanned subspace. 

Clearly, if \(H(x,y)\equiv 0\), then in particular \(H(1,s) = 0\) for all \(s\). Also, since \(H(1,s)\) is a polynomial, if \(H(1,s)\) vanishes on an open interval, it vanishes for all
\(s\).

Next, notice that in each case, \(H(x,y)\) is a homogeneous polynomial, and denote its degree by \(d\in \N\). Then,
\[
H(x,y)=x^dH(1,y/x)=0,\;\;\forall x\neq 0.
\]
Since \(H(x,y)\) is a polynomial, \(H(x,y)\equiv 0.\)
\end{proof}

                Finally, we record the following simple observation, which will be used in Proposition~\ref{prop:DecompAndMinimalityOfJointCurve}.
\begin{lemma}\label{lem:CurveDecompforPGL}
                    Let $\f=\lie{pgl}_k(\R),$ and let $F\subset\PGL_k(\R)$ be the associated connected subgroup. Let $I_1\sqcup \cdots \sqcup I_L = \R\smallsetminus \Sigma$ be the connected components of $\R \smallsetminus \Sigma$. Then for each $1\leq i \leq L,$ there exists an analytic curve $f:I_i \to F$ such that 
                    \[\bkt{W}(s)=\bkt{J^\sigma} f(s), \quad \sigma=\begin{cases} 0, & \text{ if } k \text{ odd, or } k \text{ even and } \det(W(s))>0, \quad \forall s\in I_i, \\ 1, & \text{ if } k \text{ even and } \det(W(s))<0,\quad \forall s\in I_i,\end{cases}\]
                    where $J=\diag(-1,\Id_{k-1}).$
                \end{lemma}
                \begin{proof}
                    If $k$ is odd, then $\PGL_k(\R)$ is connected, so we can choose $\bkt{W}(s)=\R^\times f(s).$ If $k$ is even, then $\PGL_k(\R) = F \sqcup JF$ and $\bkt{W}(s) \in F$ if and only if $\det(W(s))>0.$
                \end{proof}

				\section{Linear stability}\label{sec:5}
				We recall the general notions required for us to apply \cite{Yang2020}, and note a couple of additional observations which will be useful later.
				
				Let \(\GG\) be a connected algebraic semisimple group defined over \(\R\). We denote \(G:=\GG(\R)^\circ\)  the identity component of the real points of
				\(\GG\) and $\mathbb G_m(\R) = \R^\times$. 
				Fix a maximal \(\R\)-split torus \(\TT\subseteq \GG\). 
				Denote the set of characters and cocharacters on $\TT$ by
				$$X^*(\TT):= \{ \chi:\TT\to\mathbb G_m  \},\quad X_*(\TT):= \{ \delta :\mathbb G_m\to \TT  \}, $$
				where $\chi, \delta$ are algebraic group homomorphisms. To prevent confusion with our notation of \(\delta\) in \eqref{eq:DefIwasawaCoordinates}, we stress that in this section, \(\delta\) will denote a general cocharacter.
                
				The natural pairing
				\[
				\langle\cdot,\cdot\rangle:X^*(\TT)\times X_*(\TT)\to \Z,
				\quad
				\chi(\delta(t))=t^{\langle\chi,\delta\rangle},\;t\in \mathbb G_m
				\]
				is perfect, i.e. its induced inclusions \(X_*(\TT)\hookrightarrow \Hom(X^*(\TT),\Z)\) and \(X^*(\TT)\hookrightarrow\Hom(X_*(\TT),\Z)\) are
				isomorphisms.
				
				Let \(\Phi:={}_\R\Phi(\GG,\TT)\subseteq X^*(\TT)\) be the set of non-trivial \(\R\)-roots for the adjoint action of \(\TT\) on the Lie algebra of \(\GG\). Fix a set of positive roots \(\Phi^+\subseteq \Phi\). We have that \(\Phi^+\) is closed under addition and \(\Phi=\Phi^+\sqcup (-\Phi^+)\).  We denote the set of dominant cocharacters by
				$$X_*(\TT)^+:= \{\delta\in X_*(\TT): \gen{\alpha,\delta}\geq 0, \, \forall \alpha\in \Phi^+\}. $$
				Note that $X_*(\TT)^+$ contains the zero cocharacter. 
				
				Since \(\TT\) is split over \(\R\), every cocharacter \(\delta \in X_*(\TT)\) induces a homomorphism on real points \(\delta_{\R} \colon \R^\times := \mathbb G_m(\R) \to T := \TT(\R)\). By abuse of notation, we denote this induced homomorphism by \(\delta\) as well.
				
				Denote by
				\(H_\delta :=\delta_*\frac{d}{dt}\in\mathfrak t:=\Lie(T)\) the derivative 
				of \(\delta(t)\) at \(t=1\). The Killing inner product and the Killing length on $X_*(\TT)$ are defined by
				\begin{equation}\label{def:Killing}
					(\delta_1, \delta_2)
					:=
					\frac{1}{2}\Tr\left(
					\ad H_{\delta_1} \ad H_{\delta_2}
					\right),
					\quad
					\|\delta\|^2
					:=(\delta, \delta)= \sum_{\alpha \in \Phi^+}m_\alpha\gen{\alpha,\delta}^2, 
				\end{equation}
				where \(m_\alpha\) is the dimension of the root-space corresponding to \(\alpha\).
				
				For a  cocharacter $a:\mathbb G_m\to \GG$, define the parabolic subgroup 
                \begin{equation}\label{eq:PaStableGp}
                    P(a):=\{g\in G:\lim_{t\to\infty}a(t)ga(t)^{-1}
				\text{ exists}\}.
                \end{equation}
				Assume in addition that \(a\in X_*(\TT)^+\), and choose a minimal parabolic subgroup 
				\(B\subseteq P(a)\) whose Lie algebra of \(B\) is spanned by the \(\TT\)-negative root spaces together with the zero weights.
				
				We note that for a general cocharacter \(a: \mathbb G_m\to \GG \), we may always  choose a maximal $\R$-split torus $\TT$ containing $a$, i.e. \(a: \mathbb G_m\to \TT \), and we may choose $\Phi^+$ so that $a$ is dominant.

				Denote the Weyl group by \(W_G:=N_G(T)/Z_G(T)\).

				\begin{definition}\label{def:WeakUNstab} For \(w\in W=W_G\), denote a representative in \(N_G(T)\) again by $w$, and define \(X_w\) to be the Zariski closure in \(G/P(a)\) of \(BwP(a)/P(a)\subseteq G/P(a)\).
					\begin{enumerate} 
						\item[(1)] A \emph{Schubert variety} in \(G/P(a)\)
						is a subset of the form \(gX_w\subseteq G/P(a)\), where \(g\in G\).
						\item[(2)] A Schubert variety \(gX_w \) in \(G/P(a)\) is called \emph{weakly unstable with
							respect to \(a\)}, if there exists 
						\(0\neq\delta\in X_*(T)^+\) such that \((\delta,a^w)\ge 0\),
						where \(a^w:=waw^{-1}\).
					\end{enumerate}
				\end{definition}
				We have the following lemma, which follows from definitions.
				
				\begin{lemma}\label{lem:productWeakSchubertStability}
					For \(i=1,2,\) let \(\GG_i\)  be a semi-simple group over \(\R\), and let \(\TT_i\subseteq \GG_i\) be a maximal \(\R\)-split torus. Choose positive roots, and let \(d=(d_1,d_2):\mathbb G_m\to\TT_1\times \TT_2\) be a dominant cocharacter over \(\R\).
					
					Denote \(G_i:=\GG_i(\R)\).  Then,
					\[
					(G_1\times G_2)/P_G(d)\simeq G_1/P_{G_1}(d_1)\times G_2/P_{G_2}(d_2),
					\quad
					W_G
					=
					W_{G_1}\times W_{G_2}.
					\]
					Consequently, Schubert varieties are of the form \(g X_w=g_1 X_{w_1}\times g_2 X_{w_2}\) for \(w=(w_1,w_2)\in W_G\), and a Schubert variety is weakly unstable with respect to \(d\) if
					and only if at least one of \(g_iX_{w_i}\) is weakly unstable with respect to
					\(d_i\).
				\end{lemma}
				
				\begin{proof}
					The identity \(P_G(d)=P_{G_1}(d_1)\times P_{G_2}(d_2)\) immediately follows
					from the definition.
					The decompositions of the \((G_1\times G_2)/P(d)\), Weyl group, and Schubert varieties
					then follow from the product decompositions of \(G:=G_1\times G_2\), \(P(d)\), \(B\), and
					\(W\). Thus \(X_{(w_1,w_2)}=X_{w_1}\times X_{w_2}\).
					
					Every dominant cocharacter of $G$ has the form
					\(\delta=(\delta_1,\delta_2)\), where each \(\delta_i\) is dominant for \(G_i\),
					and not both \(\delta_i\) are zero. The Killing inner product splits as
					\[
					(\delta,d^w)
					=
					(\delta_1,d_1^{w_1})+(\delta_2,d_2^{w_2}).
					\]
					Hence, if \((\delta,d^w)\ge 0\) for some nonzero dominant \(\delta\), then
					one of the two summands is nonnegative for some nonzero dominant \(\delta_i\). Conversely if one of the summands, say $(\delta_1,d_1^{w_1}) \geq 0,$ then $\delta =( \delta_1, 0)$ is dominant for which  \((\delta,d^w)\ge 0.\)
				\end{proof}
				
				Our goal now is to prove the following lemma. 
				\begin{lemma}\label{lem:NotSchubInGImpliesNonSchubInF}
					Let \(\GG\) be a semisimple group over \(\R\), let \(\mathbf F\subseteq\GG\) be a semi-simple subgroup defined
					over \(\R\) which is split over \(\R\).  Let \(I\subset\R\) be an interval and let $\phi:I\to F$ be an \emph{analytic} curve, and let \(a:\mathbb G_m\to\mathbf F\) be a cocharacter. 
					
					Suppose that the curve \(\{\phi^{-1}(s)P_G(a)/P_G(a)\}_{\;s\in I} \subseteq G/P_G(a)\) is not contained in any weakly-unstable Schubert variety in \(G/P_G(a)\). Then the curve \(\{\phi^{-1}(s)P_F(a)/P_F(a)\}_{s\in I}\) is not contained in any weakly-unstable Schubert variety in \(F/P_F(a)\).
				\end{lemma}

				The following is needed for the proof of Lemma~\ref{lem:NotSchubInGImpliesNonSchubInF}. 
				
				\begin{lemma}\label{lem:SchubRep}
					Let \(\FF\) be a connected semisimple algebraic group over \(\R\) which
					is \(\R\)-split, and denote \(F:=\FF(\R)\). Fix a cocharacter \(a:\mathbb G_m\to\FF\), and denote
					\(P:=P_F(a)\). Let \(\TT\subset\FF\) be a maximal \(\R\)-split torus
					containing \(a\), and choose the positive roots so that \(a\) is dominant.
					Let \(\Phi:=\Phi(\FF,\TT)\), and let \(B\subseteq P\) be the connected real
					Borel subgroup whose Lie algebra is
					\[
					\Lie(B)
					=
					\lie{t}\oplus\bigoplus_{\alpha\in\Phi^-}\lie{f}_\alpha,
					\quad
					\lie{t}:=\Lie(T),\;T:=\TT(\R).
					\]
					Suppose that for \(w \in N_F(T)/T\),
					\(
					X_w:=\overline{BwP/P}\subseteq F/P,
					\)
					is weakly unstable with respect to \(a\).
					
					Then, there exist a finite-dimensional rational representation
					\(V\) over \(\R\) of \(\FF\) and a vector \(v\in V\) such that \(F\cdot v\) is not
					closed and \(\pi_P^{-1}(X_w)
					\subset
					\{h\in F:h\cdot v\in V^{0-}(a)\},\) where
					\[
					\pi_P(h):=h^{-1}P,
					\quad
					V^{0-}(a)
					:=
					\left\{
					u\in V:
					\lim_{t\to\infty}a(t)\cdot u
					\text{ exists}
					\right\}.
					\]
				\end{lemma}
				
				\begin{proof}
					Since \(X_w\) is weakly unstable, there exists a nonzero dominant
					cocharacter \(\delta\in X_*(\TT)^+\) such that \((\delta,a^w)\geq 0.\) By the perfect character--cocharacter pairing, there exists a character
					\(\widehat{\delta}\in X^*(\TT)\) characterized by
					\[
					\langle\widehat{\delta},\eta\rangle=(\delta,\eta)
					\quad
					\text{for every }\eta\in X_*(\TT).
					\]
					The character \(\widehat{\delta}\) is dominant, see e.g. \cite[\S 2.1]{ShahYang2024}. 
					
					Then, by the highest-weight theorem, see e.g. \cite[Théorème~7.2]{Tits1971}, there exists an irreducible
					finite-dimensional rational representation \(W\) over \(\R\) of \(\FF\),
					with highest weight \(\widehat{\delta}\). Let \(V:=W^*\) be the dual
					representation, that is,
					\[
					(h\cdot\varphi)(x):=\varphi(h^{-1}\cdot x),
					\quad
					h\in F,\; \varphi\in V,\; x\in W.
					\]
					Then \(-\widehat{\delta}\) is a lowest weight of \(V\), and the
					corresponding weight space \(V_{-\widehat{\delta}}\) is a one-dimensional
					\(B\)-stable subspace.
					
					Fix \(0\neq v\in V_{-\widehat{\delta}}\). Since \(v\) has
					\(\TT\)-weight \(-\widehat{\delta}\), we have
					\[
					\delta(t)\cdot v
					=
					t^{\langle-\widehat{\delta},\delta\rangle}v
					=
					t^{-(\delta,\delta)}v,\;t>0.
					\]
					Since \(\delta\neq0\), we have \((\delta,\delta)>0\) and thus
					\(
					\lim_{t\to\infty}\delta(t)\cdot v=0.
					\)
					Thus \(0\in\overline{F\cdot v} \smallsetminus F\cdot v\).
					Consequently, \(F\cdot v\) is not closed.
					
					Now let \(h\in\pi_P^{-1}(BwP/P)\). Since
					\(\pi_P(h)=h^{-1}P\), we have \(h\in Pw^{-1}B\), and we may write
					\[
					h=pw^{-1}b,
					\quad
					p\in P,\quad b\in B.
					\]
					Since \(V_{-\widehat{\delta}}\) is \(B\)-stable,
					we have \(b\cdot v\in V_{-\widehat{\delta}}\). Hence
					\[
					a(t)w^{-1}(b\cdot v)
					=
					w^{-1}\bigl(wa(t)w^{-1}\bigr)\cdot (b\cdot v)
					=
					w^{-1}a^w(t)\cdot(b\cdot v) 
					=
					t^{\langle-\widehat{\delta},a^w\rangle}w^{-1}\cdot(b\cdot v).
					\]
					Moreover, \(\langle-\widehat{\delta},a^w\rangle
					=
					-(\delta,a^w)
					\leq0.\) It follows that \(\lim_{t\to\infty}a(t)w^{-1}\cdot(b\cdot v)\) exists, and
					therefore \(w^{-1}b\cdot v\in V^{0-}(a)\).
					
					Note that \(V^{0-}(a)\) is \(P\)-invariant, thus
					\(
					h\cdot v
					=
					pw^{-1}b\cdot v
					\in V^{0-}(a).
					\)
					We have thus shown that
					\[
					Pw^{-1}B
					\subset
					\{h\in F:h\cdot v\in V^{0-}(a)\}.
					\]
					The set on the right is closed because \(V^{0-}(a)\) is a linear
					subspace of \(V\). Taking closures gives
					\[
					\pi_P^{-1}(X_w)
					\subset
					\{h\in F:h\cdot v\in V^{0-}(a)\}.
					\]
					This proves the lemma.
				\end{proof}

				\begin{proof}[Proof of Lemma~\ref{lem:NotSchubInGImpliesNonSchubInF}]
					Suppose for contradiction that
					\(\widetilde\phi(s):=\phi(s)^{-1}P_F(a)/P_F(a),\;s\in I\) is contained in a weakly unstable Schubert variety in
					\(F/P_F(a)\). Thus there exist \(f_0\in F\) and a weakly unstable standard
					Schubert variety \(X_w^F\subseteq F/P_F(a)\) such that \(\widetilde\phi(I)\subseteq f_0X_w^F.\)
					
					It will be convenient at this point to replace $\phi$ by $\psi(s)=\phi(s)f_0$, so that  $\widetilde\psi(s)=f_0^{-1}\widetilde \phi(s)\in X_w^F$.
                    
					Thus
					\[
					\psi(I)\subseteq \pi_{P_F(a)}^{-1}(X_w^F),
					\]
					where \(\pi_{P_F(a)}(h)=h^{-1}P_F(a).\)
					
					By  Lemma~\ref{lem:SchubRep},  there exist a finite-dimensional
					rational \(F\)-representation \(V_F\) over \(\R\) and a vector \(v_F\in V_F\) such that \(F\cdot v_F\) is not
					closed and
					\[
					\psi(s)\cdot v_F\in V_F^{0-}(a),
					\quad\text{for all }s\in I.
					\]
					
					Let \(F\) act on \(G\times V_F\) on the right by
					\[
					(g,u)\cdot f=(gf,f^{-1}u),\;g\in G,\,u\in V_F,\,f\in F,
					\]
					and denote the quotient vector bundle with fiber $V_F$ (see \cite[Definition 2.1]{Timashev2011}) by 
					\[
					X:=G\times^F V_F:=(G\times V_F)/F.
					\]
					Since \(F\) is reductive and \(V_F\) is affine \cite{Matsushima1960}, \(X\) is
					affine by  \cite[Chapter~1, \S2, Theorem~2.1]{Timashev2011}. 
					
					We consider the \(G\)-action on \(X\) defined by
					\[g\cdot [x,u]:=[gx,u],\;g\in G,\,[x,u]\in X.\]
                    
					By
					the equivariant embedding theorem for affine \(G\)-varieties, \(X\) admits a
					\(G\)-equivariant closed embedding into a finite-dimensional rational
					\(G\)-representation \(W\) over \(\R\); see e.g. \cite[Proposition~1.9]{Brion2010Actions} for the proof over \(\C\). A proof over \(\R\) is essentially verbatim. 
					
					Now denote \(x:=[e,v_F]\in X\subseteq W.\) Since the orbit \(F\cdot v_F\) is not closed, then the \(G\)-orbit \(G\cdot x\) is not closed; indeed, if \(u\in \overline{F\cdot v_F}\smallsetminus F\cdot v_F\), and there exists \(g\in G\) such that \(g\cdot [e,v_F]=[e,u]\), then \(g\in F\), which is a contradiction. In addition, for every \(s\in I\), since \(\psi(s)\in F\),
					\[
					\psi(s)\cdot x=[e,\psi(s)\cdot v_F],
					\]
					and since \(\psi(s)\cdot v_F\in V_F^{0-}(a)\), the limit
					\[
					\lim_{t\to\infty}a(t)\psi(s)\cdot x
					=
					\lim_{t\to\infty}[e,a(t)\psi(s)\cdot v_F]
					\]
					exists in \(X\), hence also in \(W\). Therefore
					\[
					\psi(s)\cdot x\in W_G^{0-}(a):=\{w\in W: \lim_{t\to\infty}a(t)\cdot w\text{ exists}\}
					\;\text{for all }s\in I.
					\]
					
					By \cite[Proposition 2.7]{Yang2020}  together with
					\cite[Lemma~2.3(c)]{Yang2020}, we conclude that \(\widetilde \psi(s),\;s\in I\) is contained in a finite union of weakly unstable Schubert varieties in
					\(G/P_G(a)\). Since \(\psi\) is analytic,  \(\widetilde\psi(I)\) is contained in a weakly unstable Schubert variety in \(G/P_G(a)\).
                    
                    Consequently, \(\widetilde\phi(I)=f_0\widetilde\psi(I)\) is contained in a weakly unstable Schubert variety in \(G/P_G(a)\). This is a contradiction.
				\end{proof}
				
				\subsection{Stability of Wronskian curves}\label{sec:StabOfWronskians}
				We denote \(\GG:=\SL_k\),  $G:=\SL_k(\R)$. We fix the cocharacter $$a(t):=\diag(t^{k-1},t^{k-3},\dots,t^{-(k-1)}),\;t\in \mathbb G_m.$$ 
                Here, notice that \(a(t)=\delta_{k-1}(2\log t),\;t>0\), where \(\delta_{k-1}\) is given by \eqref{eq:DefOfrhonum}.
                
                The maximal \(\R\)-split torus \(\TT\subseteq \SL_k\) containing \(a\) is the full diagonal subgroup
				\[T:=\TT(\R)=\{\diag(t_1,\dots,t_k):\;t_1\cdots t_k=1\}\subseteq G.\]
				
				Choose the usual positive roots $\Phi^+=\{\alpha_{i,j}(\diag(t_1,\dots,t_k)):=t_i/t_j:1\leq i<j\leq k\},$ so that \(a\) is dominant. We let $B^-$ be the lower triangular Borel subgroup. 
				\begin{lemma} We have
					\(B^-=P(a)\), where \(P(a)\) is defined in \eqref{eq:PaStableGp}
				\end{lemma}

				Let $\Wp(s)\in\GL_k(\R),\;s\in \R\smallsetminus\Sigma$, \(\Sigma:=\{s\in\R:\,\det\Wp(s)=0\}\), be a Wronskian curve corresponding to linearly independent polynomials \(P_1(x,y),\dots,P_k(x,y)\), see \eqref{eq:DefOfWronskian}.
                
                It will be convenient to consider the normalization
                \begin{equation}
                    \widehat \Wp(s):=S(s)\Wp(s),\; S(s):= \diag(\det W_P(s)^{-1},1,\ldots,1),\; s\in \R\smallsetminus\Sigma,
                \end{equation}
        so that \(\widehat \Wp(s)\in \SL_k(\R)\), for all \(s\in \R\smallsetminus\Sigma\).

                The main goal of this section is to prove the following statement.
				\begin{proposition}\label{prop:WronskCurveNotContainedInSchub}
					For every nonempty open interval $I\subset\R\smallsetminus \Sigma$, the curve $\widehat \Wp(s)^{-1}B^-/B^-\in G/B^-$ is not contained in any weakly unstable Schubert variety with respect to \(a\). 
				\end{proposition}
                \begin{remark}\label{rem:ReductionStabilityFromSLToPGL}
				    Consider the image of \(\bkt{a(t)}:=\R^\times a(t)\in \PGL_k(\R)\). Then, \(P_{\PGL_k(\R)}(\bkt{a})=\operatorname{B}^-\) the lower-triangular Borel subgroup of \(\PGL_k(\R)\). 
                    
                    We have that \(\PGL_k(\R)/\operatorname{B}^-\cong \SL_k(\R)/B^-\), and in particular, the weakly-unstable Schubert varieties in \(\PGL_k(\R)/\operatorname{B}^-\) with respect to \(\bkt a\) are identified with the weakly-unstable Schubert varieties in \(\SL_k(\R)/B^-\) with respect to \(a\). Importantly, for \(\bkt{\Wp}(s)=\R^\times \Wp(s),\;s\in \R\smallsetminus\Sigma\), we get that  \(\bkt{\Wp}(s)^{-1} \operatorname{B}^-/B^-\) is not contained in a weakly unstable Schubert variety in \(\PGL_k(\R)/\operatorname{B}^-\) with respect to \(\bkt a\) if and only if $\widehat \Wp(s)^{-1}B^-/B^-\in G/B^-$ is not contained in any weakly unstable Schubert variety with respect to \(a\).
				\end{remark}
				
				The Schubert varieties in this setup are the (usual) closures of the Bruhat cells, and we will describe them geometrically in terms of the space of full flags in \(\R^k\) following \cite[Section~6]{Fulton1992Flags}. To prove Proposition~\ref{prop:WronskCurveNotContainedInSchub}, we need to derive Proposition~\ref{prop:NecessaryCondForWeakUnstabSchubert} below, which gives a lower bound on the dimensions of the intersections of the subspaces varying in the flags of the weakly-unstable Schubert varieties with certain fixed subspaces.
				
				Denote by $\F$ the space of all full flags in $\R^k$, i.e. sequences of nested subspaces
				\[F_\bullet:=\{0\}\subseteq F_1\subseteq\dots \subseteq F_k=\R^k,\quad \dim F_i=i.\]
				Denote by \(\be_1,\dots,\be_k\in \R^k\) the standard basis of \(\R^k\). We have that $B^-$ is the stabilizer of the flag $E_{\bullet}$, formed by the subspaces $$E_i:= \gen{\be_{k-i+1},\dots,\be_k},\;\;i=1,\dots,k.$$  
                Thus we will identify $\F$ with
				$ G/B^-$ via the map $gE_\bullet \mapsto gB^-.$
				
				The Weyl group \(\Weyl:=N_G(T)/T \) is isomorphic to the group $\upS_{k}$ of permutations on \(\set{1,\dots,k}\). We will express elements $w\in\Weyl$ as permutations, using the notation
				\[w=(w(1),w(2),\dots,w(k)).\]
				By the Bruhat decomposition \cite[Theorem~16.1.3]{springer1998}, 
				$$\mathcal F = G/B^-=\bigsqcup_{w\in \Weyl} C_w,$$
				where $C_w$ is the Bruhat cell
				\[
				C_w:=B^-wB^-/B^- ,\;w\in\Weyl.
				\]
				A Schubert variety associated with $w\in\Weyl$ is of the form $gX_w=g \overline{ B^-wB^-/B^-} =g \overline{C_w}.$
				
				Notice that for \(w\in \Weyl\) and for \(i=1,\dots,k\), we have  \(w \gen{e_i}_{\R}=\gen{e_{w(i)}}\). We  define \(r_w(p,q)\) to be the number of the standard basis lines in \(E_p\) whose images under
				\(w\) lie in \(E_q\)
				\begin{equation}\label{eq:rwpqByNumOfElts}
					r_w(p,q)
					:=
					\#\{\ell\in\{k-p+1,\ldots,k\}: k-q+1 \leq w(\ell)\leq k\}, \quad 1\le p,q\le k. 
				\end{equation}
				Then, $\dim(wE_p\cap E_q)=r_w(p,q).$
				
				We claim that the Bruhat cell $C_w$ satisfies
				\[
				C_w = \{F_\bullet\in\mathcal F: \dim(F_p\cap E_q)=r_w(p,q)\text{ for all }1\le p,q\le k\}.
				\]
				Indeed, if \(F_\bullet\in C_w\), then \(F_\bullet=bwE_\bullet\) for some
				\(b\in B^-\). Since \(B^-\) preserves \(E_\bullet\), we get
				\[
				\dim(F_p\cap E_q)= \dim(bwE_p\cap bE_q)= \dim(wE_p\cap E_q) = r_w(p,q).
				\]
				Conversely, suppose \(F_\bullet\) satisfies these equalities. By the Bruhat
				decomposition, \(F_\bullet\in C_u\) for a unique \(u\in W\). Applying the
				previous paragraph to \(C_u\), we obtain
				\[
				r_u(p,q)=\dim(F_p\cap E_q)=r_w(p,q)
				\quad\forall p,q.
				\]
				Notice that  $r_w(p,q)$ determines a permutation, and therefore determines a Weyl group element. Thus, \(u=w\), and 
				\(F_\bullet\in C_w\).
				
				By \cite[Lemma~6.1(c)]{Fulton1992Flags},
				\begin{equation}\label{eq:SchubVarInTermsOfRankFunc}
					X_w = \{F_\bullet\in\mathcal F: \dim(F_p\cap E_q)\ge r_w(p,q) \text{ for all }1\le p,q\le k\}.
				\end{equation}
				
				For every dominant cocharacter \(\delta\in X_*(\TT)\), there exists an exponent vector $(r_1,\dots,r_k)\in \Z^k$  with $\sum_{i=1}^k r_i=0$ such that
				\[\delta(t)=\diag(t^{r_1},\dots,t^{r_k}), \;t\in\R^\times.\]
				The Weyl group action on $X_*(T)$ given by
				$w\cdot \delta(t):=w \delta(t) w^{-1}=:\delta^{w}(t),\;w\in\Weyl$, is expressed explicitly as
				\[\diag(t^{r_1},\dots,t^{r_k})^w=\diag(t^{r_{w^{-1}(1)}},\dots,t^{r_{w^{-1}(k)}}),\;\;w\in\Weyl.\]
				Given $\delta_1,~\delta_2\in X_*(\TT)$, set
				$\delta_j(t)=\diag(t^{r_{j,1}},\dots,t^{r_{j,k}}),\;j=1,2$.	
				The Killing length and the Killing inner product are expressed, respectively, by
				\[ \|\delta\|^2=\frac{1}{2}\tr(\ad(\diag(r_{1},\dots,r_{k}))^2)=k\sum_{i=1}^k r_i^2, \quad(\delta_1,\delta_2)=k\sum_{i=1}^k r_{1,i}r_{2,i}.\]
				The roots \(\Phi\) of \(\TT\) are given by $\Phi :=\{\alpha_{i,j}\}_{1 \leq i \neq j \leq k}\subseteq X^*(\TT)$, where
				\[\gen{\alpha_{i,j},\delta}= \gen{\alpha_{i,j} ,\diag(t^{r_1},\dots,t^{r_k})}:=r_i-r_j,\]
				and we choose the positive roots 
				$\Phi^+:=\{\alpha_{i,j}:i<j\}.$ 
				Then, the dominant cocharacters are
				$$\delta(t)=\diag(t^{r_1},\dots,t^{r_k}) \;\;{\text{with}} 
				\;\; r_1\geq r_2\geq\dots\geq r_k.$$
				
				The following lemma describes the permutations that determine the weakly unstable Schubert varieties.
				\begin{lemma}\label{lem:UnstablePermConstraint} Given $w\in \Weyl,$ there exists a nontrivial dominant $\delta\in X_*(T)^+$ such that $( \delta,a^w)\geq 0$ if and only if there exists $q \in \{1,\dots,k-1\}$ such that 
					\begin{equation}
						w^{-1}(1)+\dots+w^{-1}(q)\leq \frac{q(k+1)}{2}.
					\end{equation}
				\end{lemma}
				\begin{proof}
					For $q\in \{1,\dots,k-1\}$, consider the  dominant cocharacters $\delta_q$ whose exponent vector is of the form
					\[\textbf{H}_q:=(k-q,\dots,k-q,\,-q,\dots,-q)\in \Z^k,\]
					where \(k-q\) is repeated \(q\) times and \(-q\) is repeated \(k-q\) times. It follows that
					\begin{align*}
						\frac1k (\delta_q,a^w)&=\textbf{H}_q\cdot(k+1-2w^{-1}(1),\dots,k+1-2w^{-1}(k))\\
						&=(k-q)\bigl((k+1)q-2\sum_{1\leq i\leq q}w^{-1}(i)\bigr)-q\bigl((k+1)(k-q)-2\sum_{q<i\leq k}w^{-1}(i)\bigr)\\
						&=qk(k+1)-2k\sum_{1\leq i\leq q}w^{-1}(i),
					\end{align*}
					with the last step following from \(\paren{\sum_{1 \leq i\leq q} +\sum_{q<i\leq k}} w^{-1}(i)=\frac{k(k+1)}{2}\).
					Thus $(\delta_q,a^w)\geq 0$ if and only if 
					\[	w^{-1}(1)+\dots+w^{-1}(q)\leq \frac{q(k+1)}{2}.\]
					
					To complete the proof, we note that every exponent vector of a non-trivial dominant cocharacter 
					\[(r_1,\dots,r_k)\in\Z^k,\quad r_1\geq\dots\geq r_k,\quad\sum_{i=1}^kr_i=0\]
					can  be written as a rational linear combination of $\textbf{H}_1,\dots,\textbf{H}_{k-1}$ with nonnegative but not all zero coefficients. Indeed, the latter claim is obtained using that
					\[(\be_i-\be_{i+1})\cdot \frac{1}{k}H_q = \begin{cases}
						1,&\text{ if } i=q,\\
						0,&\text{ otherwise},
					\end{cases}\]
					for \(i=1,\dots,k-1\), where \(\be_j\in\Z^k\) denotes a standard basis vector for \(j=1,\dots,k\).
					
					Therefore, if $(\delta,a^w)\geq 0$ for some $\delta\in X_*(T)^+$, then there exists $q\in\{1,\dots,k-1\}$ such that $(\delta_q,a^w)\geq0$, which completes the proof.
				\end{proof}
				
				As mentioned above, the following statement will be key in the proof of Proposition~\ref{prop:WronskCurveNotContainedInSchub}.
				\begin{proposition}\label{prop:NecessaryCondForWeakUnstabSchubert}
					Let $gX_w$ be a weakly unstable Schubert variety as in Definition~\ref{def:WeakUNstab}. Then, there exists \(1\leq q\leq k-1\) and \(p_1,\dots,p_q\in\N\), with
					\begin{equation}\label{eq:SeqOfpiWithAvgCond}
						1\leq p_1<\dots<p_q\leq k,\quad \sum_{j=1}^qp_j\leq\frac{q(k+1)}{2},
					\end{equation}
					such that
					\begin{equation}\label{eq:LowerBoundOnDimIntersctionsCond}
						\forall F_\bullet\in g X_w,\quad \dim( F_{p_j}\cap gE_q)\geq j,\;\;\forall 1\leq j\leq q.
					\end{equation}
				\end{proposition}
				\begin{proof} It is enough to consider the case $g=\Id.$ 
					
					To satisfy \eqref{eq:LowerBoundOnDimIntersctionsCond}, consider  for \(1\leq q\leq k\) the natural numbers \(1\leq p_1<\dots<p_q\leq k\), which are chosen such that
					\begin{equation}\label{eq:ChoiceOfPi}
						\{w^{-1}(k),w^{-1}(k-1), \dots , w^{-1}(k+1-q)\}=\{k+1-p_q,\dots, k+1-p_2, k+1-p_1\}.
					\end{equation}
					With this choice, we have 
					\[r_w(p_j,q)=\#\{\ell\in\{k-p_j+1,\dots , k\}:\; w(\ell)\in\{k-q+1,\dots, k\}\}=j.\]
					Then, by \eqref{eq:SchubVarInTermsOfRankFunc}, we get that for all $F_\bullet\in X_w$, 
					$\dim( F_{p_j}\cap E_q)\geq j,\;\;\forall 1\leq j\leq q.$ 
					
					It remains to show that we can choose \(q\) such that \eqref{eq:SeqOfpiWithAvgCond} holds.
					
					By Lemma~\ref{lem:UnstablePermConstraint}, there exists $1\leq \widetilde q\leq k-1$ such that \(	w^{-1}(1)+\dots+w^{-1}(\widetilde q)\leq \frac{\widetilde q(k+1)}{2}.\)
					We let $q=k-\widetilde q$  and let $1\leq p_1<\dots<p_q\leq k$ as in \eqref{eq:ChoiceOfPi}. Then,
					\begin{align*}
						\sum_{j=1}^q (k+1-p_j) &=\sum_{j=1}^q w^{-1}(k+1-j) 
						\\&= \frac{k(k+1)}2 - \sum_{j=1}^{\widetilde q} w^{-1}(j)\\
						&\geq \frac{k(k+1)}2 - \frac{\widetilde q(k+1)}{2}=\frac{q(k+1)}{2},
					\end{align*}
					which in turn shows that \eqref{eq:SeqOfpiWithAvgCond} holds, that is
					\(\sum_{j=1}^q p_j \leq \frac{q(k+1)}2.\) 
				\end{proof}	
				
				\begin{proof}[Proof of Proposition~\ref{prop:WronskCurveNotContainedInSchub}]
					In view of the identification of \(G/B^-\) with the space of flags, the curve  \(\widehat \Wp(s)^{-1}B^-\in G/B^- ,\;s\in I\) is identified with the curve of flags \[F_\bullet(s):=\widehat \Wp(s)^{-1}E_\bullet=\Wp(s)^{-1}E_\bullet,\;s\in I.\]  It will be convenient to identify $\R^k$ with the following space of polynomials  $$V:=\gen{P_1(1,x),\dots,P_k(1,x)}_\R,$$ by identifying the isomorphism \(\iota:\R^k\to V\) defined by $$\iota(\be_i):= P_i(1,x),\;\;i=1,\dots,k.$$ 
                    Here, \(P_1(x,y),\dots,P_k(x,y)\) are the homogeneous polynomials defining \(\Wp\), see \eqref{eq:DefOfWronskian} to recall. 
					
					Fix \(s\in I\). We claim that, under this identification, the flag $\iota(F_\bullet(s))$ is given by the sequence of nested subspaces of $V$ given by:
					\begin{equation}\label{eq:DescriptionOfFisVanDer}
						\iota( F_i(s))=\{f\in V:\;f(s)=f^{(1)}(s)=\cdots=f^{(k-i-1)}(s)=0\},\; i=1,\dots,k-1. 
					\end{equation}
					  To see this, notice that for \(j=1,\dots,k\)
                    \[\iota(\Wp^{-1}(s)\be_j)= \sum_{\sigma=1}^k(\Wp^{-1}(s))_{\sigma,j} P_\sigma(1,x)=(W_{\boldP}(x)W_{\boldP}(s)^{-1}))_{1,j}.\]
					The vectors generating \(F_i(s)\) are \(W_{\boldP}^{-1}(s)\be_{k-(i-1)},\dots,W_{\boldP}^{-1}(s)\be_k\), so that,
					\[\iota(F_i(s))=\gen{ (W_{\boldP}(x)W_{\boldP}(s)^{-1}))_{1,k-(i-1)} \dots, (W_{\boldP}(x)W_{\boldP}(s)^{-1}))_{1,k}}_\R,\]
                    for \(i=1,\dots k\).
                    Evaluating at \(s\in I\), we have 
                    \[\frac{d^{\ell-1}}{dx^{\ell-1}}(W_{\boldP}(x)W_{\boldP}(s)^{-1}))_{1,j}\mid_{x=s}=(W_{\boldP}(s)W_{\boldP}(s)^{-1}))_{\ell,j}=\delta_{\ell,j},\]
                    for \(\ell,j=1,\dots,k.\) Thus, for each \(1\leq i<k\), the generators of \(\iota(F_i(s))\) satisfy the conditions of \eqref{eq:DescriptionOfFisVanDer}. 
					
					Suppose for contradiction that there exists a weakly unstable Schubert variety \(gX_w\), where $g\in G,\; w\in\Weyl,$ such that $F_\bullet(s)\in gX_w,\;\forall s\in I$. Then, according to Proposition~\ref{prop:NecessaryCondForWeakUnstabSchubert}, there exists a $q\in \{1,\dots,k-1\}$, a subspace $U\subseteq \R^k$ with $\dim(U)=q$ and natural numbers \(1\leq p_1<\cdots < p_q\leq k\) with
					\begin{equation}\label{eq:DimSumOfFs}
						\sum_{i=1}^qp_i\leq \frac{q(k+1)}{2},
					\end{equation}
					such that
					\[\dim\left(\iota(F_{p_j}(s))\cap \iota(U)\right)=\dim\left(F_{p_j}(s)\cap U\right)\geq j,\;\forall 1\leq j\leq q,\;\forall s\in I.\]
					
					Consider for \(s\in I\) and \(j\in \{1,\dots,q\}\) the map:
					\[D_{s,j}:\iota(U)\to \R^{k-p_j},\quad D_{s,j}(f):=(f(s),f^{(1)}(s),\dots,f^{(k-p_j-1)}(s))^\top.\]
					Then, by \eqref{eq:DescriptionOfFisVanDer}, $\ker D_{s,j}$ is nontrivial for all $s\in I$ and $j\in\{1,\dots,q\}$. As we will now show, this implies that \(k-p_j<q\). To verify this, choose a basis $\{f_1(x),\dots,f_q(x)\}$ for \(\iota(U)\). Observe that if \(q\leq k-p_j\), then the Wronskian matrix $W_{\textbf{f}}(s)$  associated with those polynomials  is  formed by the first \(q\) rows of a matrix representing   \(D_{s,j}\). 
                    Now $W_{\textbf{f}}(s)$ is nonsingular for all $s \in J\subseteq I$ for a suitable subinterval $J$, see Lemma~\ref{lem:WronskianNonSing}. Then, if \(q\leq k-p_j\), we get that \(\ker D_{s,j}=0\). Namely, we can conclude that if \(\ker D_{s,j}\neq 0\), then indeed, \(k-p_j<q\).
					
					Then, since $p_j$ are natural numbers strictly increasing in $j$, 
					\[k-q+j\leq p_j,\;\;\forall j\in \{1,\dots,q\}.\]
					We obtain that
					\[\sum_{j=1}^q p_j \geq q(k-q)+\frac{q(q+1)}{2}=\frac{q(2k-q+1)}{2}>\frac{q(k+1)}{2},\]
					where in the last inequality we used that \(k>q\). This contradicts \eqref{eq:DimSumOfFs}.
				\end{proof}

				\section{Joint equidistribution of expanding Wronskian curves}\label{sec:6}
				
				We state \cite[Theorem~1.5]{Yang2020} in a weaker form, suitable for our needs. Denote \(\bG:=\SL_2(\R)\times \PGL_k(\R)\), and denote \(\bGamma:=\Gamma\times \PGL_k(\Z)\), where \(\Gamma \subset \SL_2(\R)\) is a lattice.
				
				For the following, recall the one-parameter diagonal subgroup \(\bdelta_{k-1}\) defined in \eqref{eq:DefJointDiag}. Notice that \(\bdelta_{k-1}\) is contained in the subgroup \[\SL_2^\Delta(\R):=\set{(h,\bkt{\varrho_{k-1}(h)}):\,h\in\SL_2(\R)},\]
				whose Lie algebra is \(\lie{sl}_2^{\Delta}(\R)\) given in \eqref{eq:DefSL2DiagEmbAlgebra}.

                We call a subgroup \(\bF\subset \bG\)  a \emph{joint\stdlist{} subgroup}, if \(\bF\)  is either \(\SL_2^\Delta(\R)\) or \(\SL_2^\Delta(\R) \subsetneq \bF=\SL_2(\R)\times H\), where where \(H\) is\stdlist{} subgroup. 

				\begin{theorem}{\emph{\cite[Theorem~1.5]{Yang2020}}}\label{thm:PengYu}
					Let \(\bF\subset \bG\) be a joint\stdlist{} subgroup, and fix an \(\bx\in\bG/\bGamma\) such that \(\bF\bx\) is homogeneous.  Let \(I\) be a closed interval, and let \(\boldf:I\to \bF\) be an analytic curve such that the following two assumptions hold:
					\begin{enumerate}
						\item\label{enu:PengyuMinimalityAsum} Let \(P_{\bF}(\bdelta_{k-1}):=\set{h\in \bF: \lim_{t\to\infty} \bdelta_{k-1}(t) h\bdelta_{k-1}(t)^{-1}\text{ exists}}.\)  Then, for every proper closed  algebraic reductive subgroup \(\tilde \bF\subsetneq \bF\) containing \(\bdelta_{k-1}(\R)\) and every \(h\in \bF\) such that \(\tilde \bF h\bx\) is homogeneous,   \(\boldf \) is not contained in \(P_{\bF}(\bdelta_{k-1}) \tilde\bF h\).
						\item\label{enu:PengyuLinStabAsum} The  curve \(s\mapsto f(s)^{-1}P_{\bF}(\delta_{k-1})\)
                        is not contained in any weakly unstable Schubert variety \(\bF/P_{\bF}(\bdelta_{k-1})\) with respect to \(\bdelta_{k-1}\). 
					\end{enumerate}
					Then, 
					\begin{equation}\label{eq:PengYuEquid}
						\lim_{\tau\to\infty}\int_I \Phi(\bdelta_{k-1}(\tau )\boldf(s)\bx)\,\dd s =|I| \mu_{\bF\bx}(\Phi),\;\;\forall \Phi\in C_c(\bG/\bGamma),
					\end{equation}
					where \(\mu_{\bF\bx}\) is the \(\bF\)-invariant probability on \(\bF\bx\).  \end{theorem}

                    We note that \(\set{\bg\in \bG: \lim_{t\to\infty}\bdelta_{k-1}(t)\bg\bdelta_{k-1}^{-1}(t)\text{ exists}}=\bB^-\).
					Then,
					\(P_{\bF}(\bdelta_{k-1})=\bB^-\cap\bF\).
                    \begin{definition}\label{def:SmallestJointF}
				    Let \([\Lambda]\in X_k(\VV_n)\),  let \(\boldP\in\VV_n^k\) be an ordered \(\Z\)-basis for \(\Lambda\), and let \(\varphi_{\boldP}\) be given by \eqref{eq:DefJointWronskianCurve}. 
                    
                    We define \(\bF_{[\Lambda]}\subset \bG\) to be the smallest joint\stdlist{} subgroup \(\bF\) such that there exists \((\Id,[g_0])\in\bG\) with the properties that \(\varphi_{\boldP}(s)\in \bB^-\bF(\Id,[g_0])\) for all \(s\in\R\smallsetminus\Sigma\),  and such that     \(\bF(\Id,\bkt{g_0})\bx_0\) is homogeneous, where \(\bx_0 \in \bG/\bGamma\) denotes the identity coset.
				\end{definition}
               
				\begin{remark}
				    Recall the\stdlist{} subgroup \(F_{\bkt\Lambda}\subseteq\PGL_k(\R)\) given by Definition~\ref{def:FLambdaSmallestSTDSubgroupContWronsk}.  Let \(\piR:\bG\to \PGL_k(\R)\) be the natural map. Then, \(\piR(\bF_{[\Lambda]})=F_{\bkt \Lambda}\).  
				\end{remark}
				\begin{proposition}\label{prop:DecompAndMinimalityOfJointCurve}
					Let \([\Lambda]\in X_k(\VV_n)\), and let \(\boldP\in \VV_n^k\) be an ordered \(\Z\)-basis of \(\Lambda\). Consider \(\bF_{[\Lambda]}\subset \bG\) and \((\Id,\bkt{g_0})\in \bG\) be as in Definition~\ref{def:SmallestJointF}.
					
					 Then, there exists finitely many open intervals \(I^+_1,\dots, I^+_{m^+},I^-_1,\dots, I^-_{m^-}\subset \R\smallsetminus \Sigma\) such that \(\R\smallsetminus (I^+\cup I^-)\) is finite, where  \(I^+=\cup_{i=1}^{m_+} I_i^+,\;I^-=\cup_{i=1}^{m_-} I_i^-\), and in addition, such that the following holds. 
					
					There exists a diagonal matrix \(\bkt J\in \PGL_k(\R)\) and there exist analytic curves \(\bu^-:I^+\cup I^- \to \bU^-,\; \boldf:I^+\cup I^- \to \bF_{[\Lambda]} \) such that 
					\begin{equation*}
						\varphi_{\boldP}(s)=\bu^-(s) (\Id, \bkt{J^{\sigma_\pm}})\boldf(s)(\Id,\bkt{g_0}),\;\forall s\in I^{\pm},  
					\end{equation*}
					where \(\sigma_+=0,\;\sigma_-=1\). 
                    
                    The matrix \(J\) is given explicitly by
                    \[J=\begin{cases}
                        \Id,&\text{if } \Lie(\bF_{[\Lambda]})\in \set{\lie{sl}_2^{\Delta}(\R),\, \lie{sl}_2(\R)\times \g_{2(2)}},\\
                        \diag(1,-1,1,\dots,(-1)^{k-1})&\text{if }\Lie(\bF_{[\Lambda]})=\lie{sl}_2(\R)\times\varrho_{k-1}\lie{sl}_2(\R),\\
                        \diag(-1,I_{2m-1},-1)&\text{if }\Lie(\bF_{[\Lambda]})=\lie{sl}_2(\R)\times\lie{so}_{m,m+1}(\R),\\
                        \diag(I_m,-I_m)& \text{if }\Lie(\bF_{[\Lambda]})=\lie{sl}_2(\R)\times\lie{sp}_{2m}(\R),\\
                       \diag(-1,I_{k-1}) &\text{if }\Lie(\bF_{[\Lambda]})=\lie{sl}_2(\R)\times\lie{pgl}_k(\R).
                    \end{cases}\]

					Moreover, the Assumptions \eqref{enu:PengyuMinimalityAsum} and \eqref{enu:PengyuLinStabAsum} in Theorem~\ref{thm:PengYu} hold for the curve \(\boldf\) in each interval \(I^\pm_i\).
				\end{proposition}
This proposition combines the preceding decomposition and stability
results to make Theorem~\ref{thm:PengYu} applicable to the minimal
group \(\bF_{[\Lambda]}\).
The appearance of the factor \(J\) in the decomposition need not
produce two distinct homogeneous orbits.
In particular, in the \(\lie{sp}_{2m}(\R)\) case, Lemma~\ref{lem:6.6} shows that the two
associated homogeneous orbits coincide.

				\begin{proof}
					To simplify notation, write \(F:=\bF_{[\Lambda]}\subset \SL_2(\R)\times\PGL_k(\R).\) 
                    
					Let \((\Id,\bkt{g_0})\in \bG
                    \) as in the statement of the proposition, with \(\varphi_{\boldP}(s)\in\bB^-F(\Id,\bkt{g_0})\) for every
					\(s\in\R\smallsetminus\Sigma\). 
					
					Now,
					\(W_{\boldP g_0^{-1}}(s)=W_{\boldP}(s)g_0^{-1}\), and 
					consequently, putting
					\(\varphi(s):=\paren{\nu(s),\bkt{W_{\boldP g_0^{-1}}}(s)}\),
					we have \(\varphi(s)\in\bB^-F\) for every
					\(s\in\R\smallsetminus\Sigma\).
					Apply Corollary~\ref{cor:CurveDecompAndCondFordiagonallyEmbeddedPrincipalSL2}
					if \(\Lie(F)=\lie{sl}_2^\Delta(\R)\). If
					\(F=\SL_2(\R)\times H\) with \(H\subset \PGL_k(\R)\), apply the appropriate one of
					Lemmata~\ref{lem:WronskCritForSL2},
					\ref{lem:CurveDecompAndCondForOrthogonal},
					\ref{lem:CurveDecompAndCondForSymplectic}, 
					\ref{lem:CurveDecompAndCondForExecptionalG22}, and
                    \ref{lem:CurveDecompforPGL}
					to the second component of \(\varphi(s)\) while keeping the first factor unchanged.
					These results provide finitely many intervals
					\(I_1,\dots,I_m\subset\R\smallsetminus\Sigma\) whose complement in \(\R\) is finite,  a diagonal
					matrix \(J\in\GL_k(\R)\), and analytic curves
					\(\bu^-:\bigcup_i I_i\to\bU^-\) and
					\(\boldf:\bigcup_i I_i\to F\) such that
					\[
					\varphi_{\boldP}(s)
					=\bu^-(s)(\Id, \bkt{J^{\sigma_i}})
					\boldf(s)(\Id, \bkt{g_0}),
					\quad s\in I_i,\quad \sigma_i\in\set{0,1}.
					\]
					
					We now verify Assumption~\eqref{enu:PengyuMinimalityAsum}. Denote \(P:=F\cap\bB^-\). Suppose for contradiction that, for some \(i\), there exist a proper closed
					algebraic reductive subgroup \(\bdelta_{k-1}(\R)\subseteq\tilde F\subsetneq F\)
					and \(h\in F\) such that \(\boldf(I_i)\subset P\tilde Fh\) and such that \(\tilde F h (\Id,\bkt{g_0})x_0\) is homogeneous.
					
					We may assume that \(\tilde F\) is connected, after restricting
					to a nonempty subinterval and adjusting \(h\). 
					
					Since \(\bu^-(s)\), \((\Id,\bkt{J^{\sigma_i}}) \in \bB^-\), and since
					 \(P\subset \bB^-\), the decomposition of \(\varphi_{\boldP}\) gives
					\[\varphi_{\boldP}(s)\in\bB^-\tilde Fh(\Id,\bkt{g_0}),\;\forall s\in I_i. \]
                    
					By Proposition~\ref{prop:ClassificationOfFInBFg}, there exists a diagonal
					\(\balpha_0\) such that
					\(\tilde F^{\balpha_0}:=\balpha_0\tilde F\balpha_0^{-1}\) is jointly\stdlist{} subgroup, and
					\(\varphi_{\boldP}(s)\in\bB^-\tilde F^{\balpha_0}(\balpha_0h)(\Id,\bkt{g_0})\) for every
					\(s\in\R\smallsetminus\Sigma\). Also, since \(\SL_2^{\Delta}(\R)\subset \tilde F^{\balpha_0}\), we may choose \(f_0\in \tilde F^{\balpha_0}\) such that
					\[f_0\balpha_0h(\Id,\bkt{g_0})=(\Id,\bkt{g'_0}),\;\;g'_0\in\GL_k(\R),\]so that \(\tilde F^{\balpha_0}(\balpha_0h)(\Id,\bkt{g_0})=\tilde F^{\balpha_0}(\Id,\bkt{g'_0})\), and in particular,
                    \[\varphi_{\boldP}(s)\in\bB^-\tilde F^{\balpha_0}(\Id,\bkt{g'_0}),\;\forall s\in\R\smallsetminus\Sigma.\]
					
					Furthermore,
					\[\tilde F^{\balpha_0}(\Id,\bkt{g'_0})\bx_0=\balpha_0\tilde Fh(\Id,\bkt{g_0})\bx_0\] is homogeneous. 
					Thus \(\tilde F^{\balpha_0}\) satisfies the defining requirements for
					\(\bF_{[\Lambda]}\), whereas
					\(\dim \tilde F^{\balpha_0}=\dim\tilde F<\dim F\).
					This contradicts the choice of \(F\) of minimal dimension.
					Hence Assumption~\eqref{enu:PengyuMinimalityAsum} holds
					on every \(I_i\).
					
					Finally, we verify Assumption~\eqref{enu:PengyuLinStabAsum}.  Since \(\bu^{-}(s)(\Id,J^{\sigma_i})\in \bB^-\), we obtain
					\[
					(\Id,\bkt{g_0})^{-1}\boldf(s)^{-1}\bB^-
					=\varphi_{\boldP}(s)^{-1}\bB^-
					=\bigl(\nu(s)^{-1},
					\bkt{W_{\boldP}}(s)^{-1}\bigr)\bB^-,
					\quad s\in I_i.
					\]
					
					By Proposition~\ref{prop:WronskCurveNotContainedInSchub}, the curves \(\nu(s)^{-1}B_1^-\) and  \(\bkt{W_{\boldP}}(s)^{-1}B_2^-\), for \(s\in I_i\), are not contained in weakly-unstable Schubert varieties in \(\SL_2(\R)/B_1^-\) and \(\PGL_k(\R)/B_2^-\) with respect to \(\delta_1\) and \(\delta_{k-1}\) correspondingly (see also Remark~\ref{rem:ReductionStabilityFromSLToPGL}).
					
					Lemma~\ref{lem:productWeakSchubertStability} therefore shows
					that \(\boldf(s)^{-1}\bB^-\), \(s\in I_i\), is not contained
					in any weakly unstable Schubert variety of
					\(\bG/\bB^-\) with respect to \(\bdelta_{k-1}\).
					Finally, \(F\) is split semisimple and \(\boldf:I_i\to F\)
					is analytic. Applying
					Lemma~\ref{lem:NotSchubInGImpliesNonSchubInF}, with ambient
					group \(\bG\), shows that
					\(\boldf(s)^{-1}P_F(\bdelta_{k-1})\), \(s\in I_i\), is not
					contained in any weakly unstable Schubert variety of
					\(F/P_F(\bdelta_{k-1})\).
					This is Assumption~\eqref{enu:PengyuLinStabAsum}.
				\end{proof}
				\subsection{Proof of Theorem~\ref{thm:GeneralForm}}
				Let \(2\leq k\leq n\) be two natural numbers and let \([\Lambda]\in X_k(\VV_n)\). let \(\boldP\in \VV_n^k\) be an ordered \(\Z\)-basis of \(\Lambda\). As above, let \(\varphi_{\boldP}\) be given by \eqref{eq:DefJointWronskianCurve}, and let  \(\bF_{[\Lambda]}\subset \bG\) be as in Definition~\ref{def:SmallestJointF}.				
				Let \(\emptyset\neq \calI\subset \R\) be a bounded subset composed of finitely many components. 
				Let \(I^\pm\) be the finite union of intervals as in Proposition~\ref{prop:DecompAndMinimalityOfJointCurve}, and denote \(\calI^\pm:=\calI\cap I^\pm\). Then \(\calI^\pm=\cup_{i=1}^{l_\pm} \calI^\pm_{i}\), where  \(\calI^\pm_{i}\) denotes a non-empty bounded interval (not necessarily closed or open). By Proposition~\ref{prop:DecompAndMinimalityOfJointCurve}, for each \(1\leq i\leq l_{\pm}\) there exist analytic curves \(\bu^-:I_1\cup\dots\cup I_m \to \bU^-,\; \boldf:I_1\cup\dots\cup I_m \to \bF_{[\Lambda]} \) such that for every \(i=1,\dots,m\),
				\begin{equation}\label{eq:DecompOfJointCurve}
					\varphi_{\boldP}(s)=\bu^-(s) (\Id, \bkt{J^{\sigma_\pm}})\boldf(s)(\Id,\bkt{g_0}),\;\forall s\in \calI^\pm_i.  
				\end{equation}
				Here \(J\in \GL_k(\R)\) is a diagonal matrix given in Proposition~\ref{prop:DecompAndMinimalityOfJointCurve}. 
				
				Let \(\epsilon>0\), and choose for each \(1\leq i\leq l^\pm\) a closed subinterval \(\calI_{\epsilon,i}^\pm\subset \calI_i^\pm\) such that \(|\calI_i^\pm\smallsetminus\calI_{\epsilon,i}^\pm|< \epsilon.\)  Let \(\Phi \in C_c(\bG/\bGamma)\). Then,  by Proposition~\ref{prop:DecompAndMinimalityOfJointCurve} and by Theorem~\ref{thm:PengYu}, for all \(\Psi\in C_c(\bG/\bGamma)\), and for \(\bx:=(\Id,\bkt{g_0})\bx_0\),
				\begin{align}
					|\calI_{\epsilon,i}^\pm|(\Id,\bkt{J^{\sigma_\pm}})_*\mu_{\bF_{[\Lambda]}\bx}(\Psi)&=\lim_{\tau\to\infty}\int_{\calI^\pm_{\epsilon,i}} \Psi((\Id,\bkt{J^{\sigma_\pm}})\bdelta_{k-1}(\tau)\boldf(s)\bx)\,\dd s   \nonumber\\
					&=\lim_{\tau\to\infty}\int_{\calI^\pm_{\epsilon,i}} \Psi(\bdelta_{k-1}(\tau)(\Id,\bkt{J^{\sigma_\pm}})\boldf(s)\bx)\,\dd s   \nonumber\\
					&=\lim_{\tau\to\infty}\int_{\calI^\pm_{\epsilon,i}} \Psi(\bdelta_{k-1}(\tau)\varphi_{\boldP}(s)\bx_0)\,\dd s,\nonumber
				\end{align}
				where in the last step we have used \eqref{eq:DecompOfJointCurve} with \(\lim_{\tau\to\infty}\bdelta_{k-1}(\tau)
            \bu^-(s)\bdelta_{k-1}(\tau)^{-1}=\Id\) uniformly in \(s\in\calI^\pm_{\epsilon,i}\). 

                For the following, we shall denote
                \begin{equation}\label{eq:boldCalOOrbitsDef}
\bcalO^+:=\bF_{[\Lambda]}\bx,\;\bcalO^-:=(\Id,\bkt{J})\bF_{[\Lambda]}\bx,
                \end{equation}
and we denote by \(\mu_{\bcalO^+}\) the \(\bF_{[\Lambda]}\)-invariant probability on \(\bcalO^+\), and we let \(\mu_{\bcalO^-}:=(\Id,\bkt{J})_*\mu_{\bcalO^+}\).
                
				Thus, by taking \(\epsilon\to 0\), we conclude that
				\begin{equation}\label{eq:Woohoo}
					\lim_{\tau\to\infty}\int_{\calI} \Psi(\bdelta_{k-1}(\tau)\varphi_{\boldP}(s)\bx_0)\,\dd s = |\calI^+| \mu_{\bcalO^+}(\Psi)+|\calI^-| \mu_{\bcalO^-}(\Psi).
				\end{equation}
				
				We now apply Proposition~\ref{prop:ReducBallAvergesToExpandingCurves}. Let \(\Phi\in C_c(\SL_2(\R)/\Gamma \times X_k(\VV_n))\). We have
					\begin{align*}\frac{1}{T^2}&\int_{\SL_2(\R)_T}\Phi(g\Gamma,g\cdot[\Lambda])\,\dd g  \\& =\int_0^{2\pi}\int_{-\infty}^\infty \int_{\calE(r,\theta)\smallsetminus\Sigma}\Phi_{\theta}\circ(\Id\times\Rep) \bigl(\bdelta_{k-1}(2\log(T)+r)\varphi_{\boldP}(s)\bx_0\bigr)\,\dd s\, e^r\dd r\,\dd\theta+o_{T\to\infty}(1). \end{align*}
				Here, recall that \(\Phi_{\theta}(g\Gamma, [\Lambda])=\Phi(\kappa(\theta) g \Gamma, \kappa(\theta)\cdot [\Lambda]).\) 
                For a fixed \(r\in \R\), we shall use \eqref{eq:Woohoo} to calculate
				\[\lim_{T\to\infty}\int_{\calE(r,\theta)\smallsetminus\Sigma}\Phi_{\theta}\circ(\Id\times\Rep)\bigl(\bdelta_{k-1}(2\log(T)+r)\varphi_{\boldP}(s)\bx_0\bigr)\,\dd s.\]
				Denote 
                \(\calE^\pm(r,\theta):=\calE(r,\theta)\cap I^\pm=
                \set{s\in \R:\,e^{r/2}\|\kappa(\theta)\de(\infty)\nu(s)\|\leq 1}\cap I^\pm.\)

    For fixed \(r,\theta\),  \eqref{eq:Woohoo} shows that
				the inner integral converges to
				\[|\calE^+(r,\theta)|\,(\Id,\Rep)_*\mu_{\bcalO^+}(\Phi_\theta)
				+|\calE^-(r,\theta)|\,(\Id,\Rep)_*\mu_{\bcalO^-}(\Phi_\theta).\]
				The functions \(\Phi_\theta\) are uniformly bounded, and
				equivalence of matrix norms gives
				\(e^r|\calE(r,\theta)|
				\ll e^{r/2}\mathbf{1}_{\{r\le r_0\}}\)
				for some fixed \(r_0\), uniformly in \(\theta\).
				Thus dominated convergence yields
                \begin{align}
                &\lim_{T\to\infty}\frac1{T^2}
				\int_{\SL_2(\R)_T}\Phi(g\Gamma,g\cdot[\Lambda])\,\dd g\nonumber\\
                &=\int_0^{2\pi}\int_{\R}\Bigl(|\calE^+(r,\theta)|\,(\Id,\Rep)_*\mu_{\bcalO^+}(\Phi_\theta)
				+|\calE^-(r,\theta)|\,(\Id,\Rep)_*\mu_{\bcalO^-}(\Phi_\theta)\Bigr)e^r\,\dd r\,\dd\theta.\nonumber
                \end{align}
				
				By Fubini,
				\begin{align*}
					\int_{\R}|\calE^\pm(r,\theta)|e^r\,\dd r
					&=\int_{I^\pm}
					\int_{-\infty}^{-2\log\|\kappa(\theta)\de(\infty)\nu(s)\|}
					e^r\,\dd r\,\dd s\\
					&=\int_{I^\pm}
					\|\kappa(\theta)\de(\infty)\nu(s)\|^{-2}\,\dd s
				\end{align*}
				Finally, also by applying Corollary~\ref{cor:volume-asymp-SL2},

                \begin{align}
                  &\lim_{T\to\infty}\frac1{\vol(\SL_2(\R)_T)}
				\int_{\SL_2(\R)_T}\Phi(g\Gamma,g\cdot[\Lambda])\,\dd g\nonumber\\
                &=\int_0^{2\pi}\Bigl(\omega^+(\theta)\,(\Id,\Rep)_*\mu_{\bcalO^+}(\Phi_\theta)
				+\omega^-(\theta)\,(\Id,\Rep)_*\mu_{\bcalO^-}(\Phi_\theta)\Bigr)\dd\theta,\label{eq:NormBallJointContinuousGenFormula} 
                \end{align}
				where \(\omega^\pm(\theta_0)=\frac{\int_{I^\pm}
					\|\kappa(\theta_0)\de(\infty)\nu(s)\|^{-2}\,\dd s}{2\int_0^{\pi}\int_{\R}\|\kappa(\theta)\de(\infty)\nu(s)\|^{-2}\,\dd s\,\dd\theta}.\)
				
				We now proceed to go over the list of the joint\stdlist{} subgroups, and to apply Proposition~\ref{prop:sumtointegral}.
				\subsubsection{The Case of \(\bF_{[\Lambda]}=\SL_2^\Delta(\R)\): Proof of Theorem~\ref{thm:MainTheoremGammaSpecial}}
				
				For the following, recall Definition~\ref{def:GammaSpecialPoints} of the \(\Gamma\)-special points.
				\begin{lemma}\label{lem:JointGammaSpecial}
					Let \([\Lambda]\in X_k(\VV_n)\) and let \(\bF_{[\Lambda]}\subset\bG\) be given by Definition~\ref{def:SmallestJointF}. Then  \(\Lie(\bF_{[\Lambda]})=\lie{sl}_2^\Delta(\R)
					\Longleftrightarrow
					[\Lambda]\text{ is }\Gamma\text{-special}.\)
					\end{lemma}
				
				\begin{proof}
Assume that \(\bF_{[\Lambda]}=\SL_2^{\Delta}(\R)\). Then, there exists an ordered \(\Z\)-basis \(\boldP\) for \(\Lambda\), and \(\bg:=(\Id,\bkt{g_0})\in \bG\), such that \(\varphi_{\boldP}(s)\in \bB^-\SL_2^{\Delta}(\R)\bg\) for all \(s\in \R\smallsetminus\Sigma\), and moreover, \(\SL_2^{\Delta}(\R) \bg\bx_0\) is homogeneous. 
In particular, consider
\[\stab[\R^\times g_0 \PGL_k(\Z)]:=\set{h\in\SL_2(\R):\R^\times\varrho(h)g_0 \PGL_k(\Z)=\R^\times g_0 \PGL_k(\Z)}.\]
Then, \(\stab[\bkt{g_0} \PGL_k(\Z)]\) is commensurable with \(\Gamma\), where \(\varrho:=\varrho_{k-1}\) is given by \eqref{eq:DefOfIrredRep}.

Denote, as in \eqref{eq:DefOfOrderedBasisForRinfty}, \((\hat e_1(x,y),\dots,\hat e_k(x,y)):=\paren{\frac{x^{k-1}}{0!},\frac{x^{k-2}y}{1!},\dots,\frac{y^{k-1}}{(k-1)!}}\). Apply Corollary~\ref{cor:CurveDecompAndCondFordiagonallyEmbeddedPrincipalSL2} to \(\varphi_{\boldP} \bg^{-1} = \varphi_{\boldP g_0^{-1}}\), which gives that there exists a homogeneous polynomial \(H(x,y)\) of degree \(n-(k-1)\) such that 
\[\boldP g_0^{-1} = \paren{H(x,y)\hat e_1(x,y),H(x,y)\hat e_2(x,y),\dots,H(x,y)\hat e_k(x,y)}.\]
So \(\Lambda = H(x,y)\Lambda'\), where \((\hat e_1(x,y),\dots,\hat e_k(x,y))g_0\) is an ordered \(\Z\)-basis for \(\Lambda'\). 
Recall 
\(\stab[\Lambda']:=\set{h\in\SL_2(\R): h\cdot[\Lambda']=[\Lambda']},\) where the \(\SL_2(\R)\) action is given by \eqref{eq:DefOfIrredAction}. We will now show that \(\stab[\R^\times g_0 \PGL_k(\Z)]=\stab[\Lambda']\). 

Indeed, in the coordinates of the ordered basis \((\hat e_1(x,y),\dots,\hat e_k(x,y))\) of \(\VV_{k-1}\), we obtain the identification of \(X(\VV_{k-1})\cong \PGL_k(\R)/\PGL_k(\Z)\) under which \([\Lambda']\) is identified with \(g_0 \PGL_k(\Z)\), and the representation \(\varrho\) is obtained by taking the representing matrices of the action \eqref{eq:DefOfIrredAction}, see \eqref{eq:DefOfIrredRep} to recall. This shows that \(\Lambda\) is \(\Gamma\)-special.

Vice versa, assume that \([\Lambda]\) is \(\Gamma\)-special. Then, since \(\Lambda\) is contained in an MD space, then \(\Lambda= H(x,y)\Lambda'\), where \(H\) is a homogeneous polynomial of degree \(n-(k-1)\), and there exists an ordered \(\Z\)-basis for \(\Lambda'\) of the form \((\hat e_1(x,y),\dots,\hat e_k(x,y))g_0\). Then, by Corollary~\ref{cor:CurveDecompAndCondFordiagonallyEmbeddedPrincipalSL2}, \(\varphi_{\boldP}(s)\in \bB^- \SL_2^\Delta(\R)(\Id,\bkt{g_0})\), where \(\boldP = H(x,y) (\hat e_1(x,y),\dots,\hat e_k(x,y))g_0.\) 
Finally, since \(\stab[\Lambda']\cap \Gamma\) is a lattice in \(\SL_2(\R)\), see Definition~\ref{def:GammaSpecialPoints} to recall, we get that, using the above coordinate identification, that \(\SL^{\Delta}_2(\R)(\Id,\bkt{g_0})\bx_0\) is homogeneous. This completes the proof.
\end{proof}
				
				\begin{proof}[Proof of Theorem~\ref{thm:MainTheoremGammaSpecial}]
					Suppose that \(\bkt\Lambda \in X_k(\VV_n)\) is \(\Gamma\)-special. By Lemma~\ref{lem:JointGammaSpecial}, we have \(\bF_{[\Lambda]}= \SL^\Delta_2(\R)\). According to Proposition~\ref{prop:DecompAndMinimalityOfJointCurve}, the matrix \(J\) in this case is the identity, so that \(\mu_{\bcalO^+}=\mu_{\bcalO^-}\), where \(\bcalO^\pm\) is defined in \eqref{eq:boldCalOOrbitsDef}.  Let \(\Phi\in C_c(\SL_2(\R)/\Gamma \times X_k(\VV_n))\).  Then, \eqref{eq:NormBallJointContinuousGenFormula} simplifies to
                    \[\lim_{T\to\infty}\frac1{\vol(\SL_2(\R)_T)}
				\int_{\SL_2(\R)_T}\Phi(g\Gamma,g\cdot[\Lambda])\,\dd g
                =\int_0^{2\pi} \,(\Id,\Rep)_*\mu_{\bcalO^+}(\Phi_\theta)\,\frac12 w(\theta)\,\dd\theta,\]
                  where  \(\omega^+(\theta)+\omega^-(\theta)=\frac12 w(\theta),
					\)
                    where \(w(\theta)\) is as in \eqref{eq:def-wtheta}. 

                    Let  $\varphi\in C_c(\SL_2(\R)/\Gamma)$ and let \(f\in C_c(X_k(\VV_n))\). Consider the product function 
                    \[(\varphi\otimes f)(g\Gamma,[\Lambda]):=\varphi(g\Gamma)f([\Lambda]),\; g\in\SL_2(\R),\,[\Lambda]\in X_k(\VV_n).\]
                    Recall by Definition~\ref{def:GammaSpecialPoints} of \(\Gamma\)-special points, that we may   write \(\Lambda=H\Lambda'\) where \(H\) is of degree \(d:=n-(k-1)\), and recall that \(\Gamma':=\Gamma\cap\stab[\Lambda']\) is a lattice in \(\SL_2(\R)\).
					Let \(\mu_{\Gamma'}\) be the invariant probability on \(\SL_2(\R)/\Gamma'\). 
        We claim that 
        \begin{equation}\label{eq:PushForOfJointIntegralDiagCase}
            (\Id,\Rep)_*\mu_{\bcalO^+}(\varphi\otimes f) = \int_{\SL_2(\R)/\Gamma'}
					\varphi(h\Gamma)
					f([x^d(h\cdot\Lambda')])
					\,\dd\mu_{\Gamma'}(h\Gamma').\end{equation}
                 Recall \(\bcalO^+=\{(h\Gamma,\bkt{ \varrho(h)g_0}\PGL(k,\Z)): h \in \SL(2,\R)\} \cong \SL(2,\R)/\Gamma'\). Then,
                 \[(\Id,\Rep)_*\mu_{\bcalO^+}(\varphi\otimes f) = \int_{\SL_2(\R)/\Gamma'} \varphi(h\Gamma) f\left(\Rep(\varrho(h)g_0\PGL(k,\Z))\right) d\mu_{\Gamma'}(h\Gamma').\]
Recall by definition of $\Rep$ \eqref{eq:DefRep},
\[\Rep(\varrho(h)g_0\PGL(k,\Z)) = \left[x^d\left(h \cdot \hatRep(g_0)\gen{\hat{e}_1(x,y),\dots, \hat{e}_k(x,y)}_\Z\right) \right].\]
As argued in the last paragraph of the proof of Lemma~\ref{lem:JointGammaSpecial},
\[\boldP(x,y)=H(x,y)(\hat P_1(x,y), \dots, \hat P_k(x,y))=H(x,y)(\hat e_1(x,y), \dots, \hat e_k(x,y))g_0.\]
is an ordered \(\Z\)-basis for \(\Lambda\), where \(\Lambda'=\gen{\hat P_1(x,y), \dots, \hat P_k(x,y)}_\Z\). By definition of $\hatRep,$
\[\hatRep(g_0)(\hat e_j(x,y)) = \sum_{i=1}^k (g_0)_{i,j}\hat e_i(x,y) = \sum_{i=1}^k \hat e_i(x,y)(g_0)_{i,j} = \hat P_j(x,y).\]
This proves that \(\Rep(\varrho(h)g_0\PGL(k,\Z)) = \left[x^d (h \cdot \Lambda')\right]\). Therefore, \eqref{eq:PushForOfJointIntegralDiagCase} follows.

                    We proceed to apply Proposition~\ref{prop:sumtointegral}. Let $O_\delta $ and $\varphi_\delta$ as in Proposition~\ref{prop:sumtointegral}. 
                    We conclude that
                     \begin{align*}&\lim_{T\to\infty}\frac1{\vol(\SL_2(\R)_T)}\int_{\SL_2(\R)_T}(\varphi_{\delta}\otimes f)(g\Gamma,g\cdot[\Lambda])\,\dd g
                \\
                &=\int_0^{2\pi} \int_{\SL_2(\R)/\Gamma'}
					\varphi_{\delta}(\kappa(\theta)h\Gamma)
					f(\kappa(\theta)\cdot[x^d(h\cdot\Lambda')])
					\,\dd\mu_{\Gamma'}(h\Gamma')\,\frac12w(\theta)\,\dd\theta\\
                    &=\int_0^{2\pi} \int_{\SL_2(\R)/\Gamma'}\varphi_{\delta}(h\Gamma)
					f(\kappa(\theta)\cdot[x^d(\kappa(-\theta)h\cdot\Lambda')])
					\,\dd\mu_{\Gamma'}(h\Gamma')\,\frac12w(\theta)\,\dd\theta,
                     \end{align*}
                     where in the last equality we used left invariance of \(\mu_{\Gamma'}\).
                     
                     Choose representatives \(\gamma_1,\dots,\gamma_m\) for
					\(\Gamma/\Gamma'\). By Definition~\ref{def:GammaSpecialPoints}, the \(\Gamma\)-packet is given by \([\Lambda_i']=\gamma_i\cdot[\Lambda']\) for \(i=1,\dots,m\).
					
					The fibers of the covering \(\SL_2(\R)/\Gamma'\to \SL_2(\R)/\Gamma\) are of size \(m=[\Gamma:\Gamma']\), and for a sufficiently small neighborhood \(U \Gamma / \Gamma \subset \SL_2(\R)/\Gamma\) of \(\Gamma\), the preimage of \(U\) is the disjoint union \(U\gamma_i\Gamma'\) for \(i=1,\dots,m\). 
					Normalized Haar measure assigns weight \(1/m\) to each sheet \(U\gamma_i\Gamma'\), see a more detailed argument in \cite[Section 6.2.2]{BersudskyShah}.
					Consequently, since
					\(\mu_{\SL_2(\R)/\Gamma}(\varphi_\delta)=1\) and the supports shrink
					to \(\Gamma\),
					\[
					\begin{aligned}
						&\int_{\SL_2(\R)/\Gamma'}\varphi_\delta(h\Gamma)
						f\!\left(
						\kappa(\theta)\cdot
						[x^d(\kappa(-\theta)h\cdot\Lambda')]
						\right)\,\dd\mu_\Gamma'(h\Gamma')\\
						&\quad=
						\frac1m\sum_{i=1}^m
						f\!\left(
						\kappa(\theta)\cdot
						[x^d(\kappa(-\theta)\cdot\Lambda_i')]
						\right)+o(1).
					\end{aligned}
					\]
					
					  It follows that for the limiting measure $\eta$ in \eqref{eq:NormBallJointContinuousGenFormula}, 
					\[
					\lim_{\delta\to0}\eta(\varphi_\delta\otimes f)
					=
					\frac1m\sum_{i=1}^m\int_0^{2 \pi}
					f\!\left(
					\kappa(\theta)\cdot
					[x^{n-(k-1)}(\kappa(-\theta)\cdot\Lambda_i')]
					\right)\frac{w(\theta)}{2}\,\dd\theta.
					\]
                    Finally, since  \(\kappa(\theta+\pi)=-\kappa(\theta)\), so that
                    \[\int_0^{2 \pi}
					f\!\left(
					\kappa(\theta)\cdot
					[x^{n-(k-1)}(\kappa(-\theta)\cdot\Lambda_i')]
					\right)\frac{w(\theta)}{2}\,\dd\theta=\int_0^{ \pi}
					f\!\left(
					\kappa(\theta)\cdot
					[x^{n-(k-1)}(\kappa(-\theta)\cdot\Lambda_i')]
					\right)w(\theta)\,\dd\theta.\]
					The conclusion now follows from
					Proposition~\ref{prop:sumtointegral}, namely, 
\[
\lim_{T\to\infty}\frac1{\#\Gamma_T}
\sum_{\gamma\in\Gamma_T}f(\gamma\cdot[\Lambda])
=
\frac1m\sum_{i=1}^{m}
\int_0^\pi
f\!\left(
\kappa(\theta)\cdot
\left[
x^{n-(k-1)}
\bigl(\kappa(-\theta)\cdot\Lambda_i'\bigr)
\right]
\right)
w(\theta)\,d\theta .
\]     
				\end{proof}
\subsubsection{The Case where \(\bF_{[\Lambda]}\) is a product group: Proof of Theorems~\ref{thm:PrincipalSL2Equid},~\ref{thm:mainthm_g2sppgl}}
Suppose that \(\bF_{[\Lambda]}=\SL_2(\R)\times H\), where \(\Lie(H):=\lie{h}\) is\stdlist{}. In particular, with \((\Id,\bkt{g_0}) \in \bG\) as in Definition~\ref{def:SmallestJointF}, we have \(\bF_{[\Lambda]}\bx = \SL_2(\R)/\Gamma \times H[g_0]y_0\), where \(\bx=(\Id,\bkt{g_0})\bx_0, \) with \(\bx_0\in \bG/\bGamma \) the identity coset \(\bGamma=\Gamma\times\PGL_k(\Z)\).
Let \(\calO^+:= H [g_0]y_0\) and \(\calO^-:= J\calO^+\),  with \(J\) as in Proposition~\ref{prop:DecompAndMinimalityOfJointCurve}. As mentioned in the introduction, \(J\) normalizes \(H\), so that \(\calO^-=H[Jg_0]y_0\). Let \(\mu_{\calO^\pm}\) be the \(H\)-invariant probabilities on the corresponding orbits \(\calO^\pm\). 

Then, \eqref{eq:NormBallJointContinuousGenFormula} is  rewritten as
\begin{align}
    &\lim_{T\to\infty}\frac{1}{\vol(G_T)}
    \int_{\SL_2(\R)_T}\Phi(g\Gamma,g\cdot[\Lambda])\,\dd g\nonumber\\
    &=\int_0^{2\pi}
    (\Id\times\Rep)_*\Bigl(
        \mu_{\SL_2(\R)/\Gamma}
        \otimes\mu_{\calO^+}
    \Bigr)(\Phi_\theta)
    \,\omega^+(\theta)\,\dd\theta\nonumber \\&\;\;\;+\int_0^{2\pi}
    (\Id\times\Rep)_*\Bigl(
        \mu_{\SL_2(\R)/\Gamma}
        \otimes\mu_{\calO^-}
    \Bigr)(\Phi_\theta)
    \,\omega^-(\theta)\,\dd\theta,\nonumber
\end{align}
where \(\mu_{\SL_2(\R)/\Gamma}\) is the \(\SL_2(\R)\)-invariant probability on \(\SL_2(\R)/\Gamma\). 

We proceed to apply Proposition~\ref{prop:sumtointegral}. Let \(O_{\delta}\subset \SL_2(\R)/\Gamma\) for \(0<\delta<1\) be a sequence of decreasing identity coset neighborhoods, and let $\varphi_\delta\in C_c(\SL_2(\R)/\Gamma)$ supported in $O_\delta$ such that $\mu_{\SL_2(\R)/\Gamma}(\varphi_\delta) =1.$ Also, let \(f\in C_c(X_k(\VV_n))\), and consider the product function \[\Phi(g\Gamma,[\Lambda]):=\varphi_{\delta}(g\Gamma)f([\Lambda]),\; g\in\SL_2(\R),\,[\Lambda]\in X_k(\VV_n).\]
For \(\theta \in \R\), we have \[\Phi_{\theta}(g\Gamma, [\Lambda])=\varphi_{\delta}(\kappa(\theta) g \Gamma) f(\kappa(\theta)\cdot [\Lambda])=:\varphi_{\delta,\theta}(g \Gamma) f_{\theta}([\Lambda]).\] 
Thus,
\begin{align}
    &\lim_{T\to\infty}\frac{1}{\vol(G_T)}
\int_{\SL_2(\R)_T}\Phi(g\Gamma,g\cdot[\Lambda])\,\dd g\nonumber\\
&= \int_0^{2\pi}\mu_{\SL_2(\R)/\Gamma}(\varphi_{\delta,\theta})\,\Rep_* \mu_{\calO^+}(f_{\theta})\, \omega^+(\theta)\,\dd\theta+\int_0^{2\pi}\mu_{\SL_2(\R)/\Gamma}(\varphi_{\delta,\theta})\,\Rep_* \mu_{\calO^-}(f_{\theta})\, \omega^-(\theta)\,\dd\theta.\nonumber
\end{align}
Due to \(\SL_2(\R)\)-invariance of \(\mu_{\SL_2(\R)/\Gamma}\), we have \(\mu_{\SL_2(\R)/\Gamma}(\varphi_{\delta,\theta})=\mu_{\SL_2(\R)/\Gamma}(\varphi_{\delta})=1\). Let \(\mu^\pm:=\Rep_* \mu_{\calO^\pm}\), so that \(\Rep_* \mu_{\calO^\pm}(f_{\theta})=\kappa(\theta)_*\mu^\pm(f)\), 
where the pushforward is via the irreducible action \eqref{eq:DefOfIrredAction}. 

For all \(0<\delta<1,\) it follows that
\begin{align*}
&\lim_{T\to\infty}\frac{1}{\vol(G_T)}
\int_{\SL_2(\R)_T}\varphi_{\delta}(g\Gamma)f(g\cdot[\Lambda])\,\dd g\\
&=\int_0^{2\pi} \kappa(\theta)_*\mu^+(f)\,\omega^+(\theta)\,\dd\theta+\int_0^{2\pi} \kappa(\theta)_*\mu^-(f)\,\omega^-(\theta)\,\dd\theta\\
&=\int_0^{\pi}  \kappa(\theta)_*\mu^+(f)\,2\omega^+(\theta)\,\dd\theta+\int_0^{\pi} \kappa(\theta)_*\mu^-(f)\,2\omega^-(\theta)\,\dd\theta
\end{align*}
where in the passage to the last line we used  \(\kappa(\theta+\pi)_*\mu^+(f)\,\omega^+(\theta+\pi)=\kappa(\theta)_*\mu^+(f)\,\omega^+(\theta)\), which holds since  \(\kappa(\theta+\pi)=-\kappa(\theta).\) 

By Proposition~\ref{prop:sumtointegral}, we conclude that 
\begin{align}\label{eq:WOHOOConvexCombFormulaAlmostDONE}
\lim_{T \to \infty} \frac{1}{\#\Gamma_T} \sum_{\gamma \in \Gamma_T}f(\gamma \cdot [\Lambda])=\int_0^{\pi}  \kappa(\theta)_*\mu^+(f)\,2\omega^+(\theta)\,\dd\theta+\int_0^{\pi} \kappa(\theta)_*\mu^-(f)\,2\omega^-(\theta)\,\dd\theta.
\end{align}

In many situations, we get that \(\mu:=\mu^+=\mu^-\), and \eqref{eq:WOHOOConvexCombFormulaAlmostDONE} reduces to
\begin{equation}\label{eq:WOHOONOConvexCombFormulaAlmostDONE}
    \lim_{T \to \infty} \frac{1}{\#\Gamma_T} \sum_{\gamma \in \Gamma_T}f(\gamma \cdot [\Lambda])=\int_0^{\pi}  \kappa(\theta)_*\mu(f)\,\omega(\theta)\,\dd\theta,
\end{equation}
where \(\omega(\theta)=2\omega^+(\theta)+2\omega^-(\theta).\)
This happens precisely when \(H[g_0]y_0=H[Jg_0] y_0\), where \(J\) as in Proposition~\ref{prop:DecompAndMinimalityOfJointCurve}. 

When \(\Lie(H)=\g_{2(2)},\) we have \(J=\Id\), so \eqref{eq:WOHOONOConvexCombFormulaAlmostDONE} is automatic. When \(\Lie(H)=\lie{pgl}_k(\R),\) we have $\calO^+=\calO^-,$ since $\PGL^\circ _k(\R)$ acts transitively on $\PGL_k(\R)/\PGL_k(\Z)$. When \(\Lie(H)=\lie{sp}_{2m}(\R)\), we have the following lemma which shows that  \(\mu^+=\mu^-\) when \(\Lie(\bF_{[\Lambda]})=\lie{sl}_2(\R)\times \lie{sp}_{2m}(\R).\)

Denote \(H:= \PSp_{2m}(\R)\), where \( \PSp_{2m}(\R)\) is the connected\stdlist{} subgroup corresponding to \(\lie{sp}_{2m}(\R)\). Then, \(H \bkt{g_0} y_0\) is homogeneous, where \(y_0\) denotes the identity coset in \(\PGL_{2m}(\R)/\PGL_{2m}(\Z)\). 
                \begin{lemma}\label{lem:6.6}
                    Let \(J=\diag(I_m,-I_m)\). Then, \(\bkt{J} H \bkt{g_0} y_0 = H \bkt{g_0} y_0\).
                \end{lemma}
             \begin{proof}
Let $Q$ be the matrix of the symplectic form \(q(\bx,\by)=\sum_{i=1}^{2m} (-1)^{i+1}x_{2m+1-i}y_i\) defining $H$. Namely, \( H=\{[h]\in\PGL_{2m}(\R):hQh^T=Q\}.\) 

We notice that \(JQJ^\top = -Q\), which also shows that \([J]\) normalizes \(H\). We will now show that there exists a \(\gamma \in \GL_{2m}(\Z)\) such that 
\begin{equation}\label{eq:ToBeProvedForSympFlip}
    (g_0 \gamma g_0^{-1}) Q (g_0 \gamma g_0^{-1})^\top = -Q
\end{equation}
Assuming that, this shows that \([J(g_0 \gamma g_0^{-1})]\in H\), and consequently, since  $[J]$ normalizes $H$, 
\[
    [J]H[g_0]y_0
    =H[Jg_0]y_0
    =H[g_0]y_0,
\]
as required.

Denote \( \Omega:=g_0^{-1}Q(g_0^{-1})^\top.\) We claim that $\R\Omega$ is defined over \(\Q\).

Indeed, since \(H[g_0] y_0\) is homogeneous, \(\widetilde\Delta_0 :=g_0^{-1}\Sp_{2m}(\R)g_0\cap \SL_{2m}(\Z)\) is a lattice in $g_0^{-1}\Sp_{2m}(\R)g_0$. By the Borel density theorem, see e.g. 
\cite[Chap.~V]{Raghunathan1972}, 
$\widetilde\Delta_0$ is Zariski dense in $g_0^{-1}\Sp_{2m}(\R)g_0$. Since the only alternating
matrices preserved by $g_0^{-1}\Sp_{2m}(\R)g_0$ are the scalar multiples of $\Omega$, we have
\[
    \R\Omega
    =
    \left\{
        B\in\operatorname{M}_{2m}(\R): B^\top = -B,\;
        \delta B\delta^\top=B
        \text{ for all }\delta\in\widetilde\Delta_0
    \right\}.
\]
The right-hand side is defined by linear equations with integral
coefficients, proving the observation.

Thus, there exists  $c\in\R^\times$ such that $c\Omega$ is an integral matrix.
Let \(W=\begin{psmallmatrix}
    0&1\\-1&0
\end{psmallmatrix}\). By \cite[Chap.~IX, \S5, no.~1, Thm.~1]{BourbakiAlgebreIX},
there exists $U\in\GL_{2m}(\Z)$ such that \( U(c\Omega)U^\top\) is the block diagonal matrix 
\[
    U(c\Omega)U^\top
    =
    \diag(d_1 W,d_2W,\dots,d_mW).
\]
Let \(S:=\diag(1,-1,\ldots,1,-1)\), and note that \[S \diag(d_1 W,d_2W,\dots,d_mW) S^\top = - \diag(d_1 W,d_2W,\dots,d_mW).\]
Taking \(\gamma:=U^{-1}SU\in\GL_{2m}(\Z)\), we obtain \( \gamma\Omega\gamma^\top=-\Omega.\)

Since \(\Omega=g_0^{-1}Q(g_0^{-1})^\top\), we get \eqref{eq:ToBeProvedForSympFlip}.
\end{proof}

For the cases that \(\Lie(H)\in \set{\varrho_{k-1}\lie{sl}_2(\R),\lie{so}_{m,m+1}(\R)}\), there are examples such that \(\calO^+\cap\calO^-=\emptyset.\) 

The following gives a construction of  examples for the case that \(\Lie(H)=\varrho_{k-1}\lie{sl}_2(\R)\). 

\begin{lemma}\label{lem:DisjointOLambdaOrbsExample}
		Let \(A(x,y),\,B(x,y)\) be linearly independent homogeneous polynomials of degree \(d\geq2\). Let \(k\geq 4\), and consider the \(k\)-lattice
		\[\Lambda_k:=\gen{A^{k-1},A^{k-2}B,\dots,B^{k-1}}_\Z+\gen{\frac17 A^{k-4}(A^3-2B^3)}_\Z.\]
		Let \(g_k\in \GL_k(\Q)\) be such that  \(\paren{\frac{A^{k-1}}{0!},\frac{A^{k-2}B}{1!},\dots,\frac{B^{k-1}}{(k-1)!}}g_k\) is an ordered \(\Z\)-basis for \(\Lambda_k\). 
        
        Then,  \(\calO_k^+\cap\calO_k^-=\emptyset,\) where
        \[\calO_k^+=\bkt{\varrho_{k-1}(\SL_2(\R))g_k}y_0,\;\;\calO_k^-=J\calO_k^+,\;\;J=\diag(1,-1,1,\dots,(-1)^{k-1}).\]
	\end{lemma}

\begin{proof}
It is more convenient to consider the contrapositive situation. We have \(\calO^+_k\cap\calO^-_k\neq\emptyset\) if and only if there exists \(f \in \SL_2(\R)\) such that
\[\bkt{\varrho_{k-1}(f) J}\in \bkt g \PGL_k(\Z) \bkt g^{-1}.\]

Notice that \([g] \PGL_k(\Z) [g^{-1}]\) is the stabilizer under the representation \eqref{eq:DefRep} of the homothety class \([L]\) of the \(k\)-lattice
\[L=7\gen{x^{k-1},x^{k-2}y,\dots,y^{k-1}}_\Z+\gen{x^{k-4}(x^3-2y^3)}_\Z.\]

By \eqref{eq:DefOfIrredRep}, the representing matrix of \(\begin{psmallmatrix}
    1&0\\0&-1
\end{psmallmatrix}\in \GL_2(\R)\) is \(J=\diag(1,-1,1,\dots,(-1)^{k-1})\). In particular, \(\bkt {\varrho_{k-1}(\GL_2(\R))}=\bkt {\varrho_{k-1}(\SL_2(\R))}\sqcup \bkt{\varrho_{k-1}(\SL_2(\R))J}.\) 

Thus, \(\calO^+_k\cap\calO^-_k\neq\emptyset\) if and only if there exists a \(h\in \GL_2(\R)\) with \(\det h<0\) such that \(\bkt{\varrho_{k-1}(h)}\in \bkt g \PGL_k(\Z) \bkt g^{-1},\) if and only if \(h\cdot [L]=[L]\) where the action is \eqref{eq:DefOfIrredAction}.

Assume that \(\bkt{\varrho_{k-1}(h)}\in \bkt g \PGL_k(\Z) \bkt g^{-1}\), for \(h\in \GL_2(\R)\). Then, we claim that \(\bkt h\in \PGL_2(\Q)\). Indeed, since \(\bkt {\varrho_{k-1}}\) is an algebraic immersion of \(\PGL_2(\R)\) defined over \(\Q\) onto a closed image, and preimages of \(\Q\) points under \(\Q\) immersions are also \(\Q\) points \cite[\S 1.4, pp.~7--8]{Milne2017AlgebraicGroups}. To see that \(\bkt{\varrho_{k-1}}\) is an immersion, note that the differential is nonsingular, so that the kernel is finite. A finite normal subgroup of a simple group is central, and the center of \(\PGL_2(\R)\) is trivial.

We conclude that \(\calO^+_k\cap\calO^-_k\neq\emptyset\) if and only if there exists \(A\in \GL_2(\Q)\) with \(\det A<0\) such that \(A\cdot L=tL\) for some \(t>0\), where the action is \eqref{eq:DefOfIrredAction}. The following local to global argument shows that this cannot occur. 

Choose \(A\) to be a primitive integral representative of its projective class, and write \(A\cdot L=tL\) with \(t>0\). Clearly, \(t\in\Q\). We will show that \(\det(A)>0\).

Denote \(M:=\Sym^{k-1}(\Z^2)\). For a prime \(p\), let \(\Q_p\) be the field of \(p\)-adic rational numbers, and let \(\Z_p\) be the ring of \(p\)-adic integers. Consider \(M_p=M\otimes_{\Z}\Z_p=\gen{x^{k-1},x^{k-2}y,\dots,y^{k-1}}_{\Z_p}\) and \(L_p=L\otimes_{\Z}\Z_p\).

By Cartan decomposition for \(\GL_2(\Q_p)\), there exist \(U,V\in \GL_2(\Z_p)\) and \(r_1,r_2\in \Z\) such that
\[
A=UD_pV,\qquad D_p=\begin{psmallmatrix}
    p^{r_1}&0\\0&p^{r_2}
\end{psmallmatrix}.
\]
Note that \(M_p=U\cdot M_p=V\cdot M_p\). Since \(A\) is primitive, \(\min(r_1,r_2)=0\).

Let \(p\neq 7\) be a prime. Then \(L_p=M_p\), so \(D_p\cdot M_p=tM_p\).
Since \(D_p\cdot x^{k-1}=p^{r_1(k-1)}x^{k-1}\) and \(D_p\cdot x^{k-2}y=p^{r_1(k-2)+r_2}x^{k-2}y\), it follows that
\[
r_1(k-1)=r_1(k-2)+r_2=v_p(t).
\]
Thus \(r_1=r_2=0\), and \(A\in \GL_2(\Z_p)\) for all \(p\neq 7\). In particular, \(\det A=\pm 7^q\) for some \(q\in\Z\).

We shall now show that \(\det A\) is a square in \(\Q_7\). This forces \(q\) to be even, and since \(-1\) is not a square in \(\Q_7\), it follows that \(\det A>0\).

Notice that \(7M_7\subseteq L_7\subseteq M_7\). These inclusions and \(A\cdot L_7=tL_7\) give
\[
7tM_7\subseteq A\cdot M_7\subseteq 7^{-1}tM_7.
\]
The same inclusions hold with \(D_7\) in place of \(A\). Comparing the \(x^{k-1}\) and \(y^{k-1}\) coordinates shows that both \((k-1)r_1\) and \((k-1)r_2\) lie between \(v_7(t)-1\) and \(v_7(t)+1\). Hence \((k-1)|r_1-r_2|\le2\). Since \(k\geq 4\) and \(\min(r_1,r_2)=0\), this again gives \(r_1=r_2=0\).

Thus \(A\in\GL_2(\Z_p)\) at every prime, and consequently \(A\in\GL_2(\Z)\). Then, \(A\cdot M=M\), so that action has determinant of absolute value \(1\). This gives that \(t=1\), and \(A\cdot L =L.\)

Consequently, \(A\) acts on the quotient
\[L_7/7M_7=\bigl\{\ell(x^{k-4}(x^3-2y^3))+7M_7:0\leq\ell\leq6\bigr\}\cong \gen{x^{k-4}(x^3-2y^3)}_{\mathbb F_7},\]
where \(\mathbb F_7 := \Z/7\Z.\) More concretely, consider the homogeneous polynomial function \(f(x,y)=x^{k-4}(x^3-2y^3)\) modulo \(7\), and denote \(\red_7(A)=A\mod 7.\) Then \(f((x,y)\red_7(A))=c f(x,y)\) for \(c\in \mathbb F^\times_7\). Consider the subgroup \(S\subset \GL_2(\mathbb F_7)\) which preserves \(\gen{x^{k-4}(x^3-2y^3)}_{\mathbb F_7}\).
Let \(\bkt S\) denote the image of \(S\) in \(\PGL_2(\mathbb F_7).\)

By the unique factorization theorem in \(\mathbb F_7[x,y]\), we conclude that \(\gen{x^3-2y^3}_{\mathbb F_7}\) is \(\bkt S\)-invariant. By \cite[Theorem~3.3]{StichtenothTopuzoglu2012}, it follows that \(\#\bkt S\,\mid \deg(T^3-2)=3\). Now consider the homomorphism \(\bkt\det:\bkt S\to\mathbb{F}_7^\times/(\mathbb F_7^\times)^2,\;s\mapsto \det(s)(\mathbb F_7^\times)^2\). We observe that \(\left|\mathbb{F}_7^\times/(\mathbb F_7^\times)^2\right|=2\). Thus, the only possibility is that \(\ker\bkt\det = \bkt S\), that is, the determinant of every \(s\in S\) is a square in \(\mathbb F_7\). By Hensel's lemma, \(\det(A)\) is a square in \(\Q_7\).
\end{proof}
\begin{remark}\label{rem:ExamplesDisjointOrbOrthCase}
    For \(\Lie(H)=\lie{so}_{m,m+1}(\R)\), we outline how to construct
    examples with \(g\in\GL_k(\Q)\) such that
    \(H[g]y_0\cap H[Jg]y_0=\emptyset\), where
    \(J=\diag(-1,I_{2m-1},-1)\).

    Recall that \(\SO_{m,m+1}(\R)\) has two connected components and that
    \(J\) lies in the nonidentity component. The two orbits are disjoint
    if and only if every element of \(\SO_{m,m+1}(\Q)\) preserving the
    lattice \(L=g\Z^k\) belongs to \(\SO^\circ_{m,m+1}(\R)\).

    The identity component is the kernel of the real spinor norm, which
    for rational isometries is obtained from the rational spinor norm
    (see \cite[Chapter~10]{Cassels1978} for the definition). The local argument for the determinant in the proof of Lemma~\ref{lem:DisjointOLambdaOrbsExample} can be adapted to the spinor norm.

    The idea is to find a rational lattice \(L\) that is unimodular at
    each odd prime \(p\ne7\) and admits an orthogonal \(\Z_7\)-basis
    whose norms are squares of pairwise distinct valuations.
    Unimodularity gives even valuation of the spinor norm at each odd
    prime \(p\ne7\). At \(7\), the step to verify is that every integral
    isometry of the chosen local lattice has square spinor norm.

    Once these local choices are realized in a single rational lattice
    and the assertion at \(7\) is verified, the rational spinor norm of
    an element preserving \(L\) is represented by one of \(1,-1,2,-2\).
    Since \(-1\) and \(-2\) are nonsquares in \(\Q_7\), the two negative
    classes are excluded, giving the desired disjointness.
\end{remark}

				\section{Projective \texorpdfstring{$(n-1)$}{(n-1)}-lattices and Frobenius numbers}\label{sec:7}
				We denote by \(\calP(X_{n-1}(\R^n))\) the space of probability measures on \(X_{n-1}(\R^n)\).
				
				Recall the map \(\pi:X_{n-1}(\R^n)\to\Gr_{n-1}(\R^n)\) defined by \(\pi([\Lambda]):=\gen{\Lambda}_{\R}\), and recall the disintegration formula,
				\[
				\mu=\int_{\Gr_{n-1}(\R^n)}\mu_{\V}\,d\nu(\V),\quad \nu=\pi_*\mu,
				\]
				where \(\mu_{\V}\) is a probability measure supported on the fiber \(X(\V):=\pi^{-1}(\V)\), defined for \(\nu\)-almost every \(\V\in\Gr_{n-1}(\R^n)\). Following \cite{SargentShapira2019}, we say that \(\mu\in \calP(X_{n-1}(\R^n))\) is a natural lift if, for \(\nu\)-almost every \(\V\), the measure \(\mu_{\V}\) is the unique \(\PGL_{n-1}(\R)\)-invariant probability measure on \(X(\V)\cong\PGL_{n-1}(\R)/\PGL_{n-1}(\Z)\).
				
				For \(\ba\in \Zprim^n,\) define
                \begin{equation}\label{eq:DefOfOrthogonalLattice}
                    \Lambda_{\ba}:=\ba^\perp \cap\Z^n.
                \end{equation}

				\begin{proposition}\label{prop:PushforwardTheta}
					Let \(n\geq3\). For \(i\in\N\), let \(\Omega_i\subset\Zprim^n\cap\R^n_{>0}\) be finite and let \(\mu_i\in\calP(X_{n-1}(\R^n))\) be the uniform probability measure with support \(\set{[\Lambda_{\ba}]:\ba\in\Omega_i}\). Assume that \(\lim_{i\to\infty}\mu_i=\mu\) in the weak-* topology, where \(\mu\) is a natural lift.
					
					Let \(\calD\subset\Gr_{n-1}^+(\R^n):=\set{\bv^\perp:\bv\in\R_{>0}^n}\) be a Borel set whose boundary has measure zero with respect to \(\nu:=\pi_*\mu\), and let \(\alpha\in\R\). Then
					\[
					\lim_{i\to\infty}\frac{\left|\set{\ba\in\Omega_i:\frac{F(\ba)}{(a_1\cdots a_n)^{1/(n-1)}}\geq\alpha,\ \ba^\perp\in\calD}\right|}{|\Omega_i|}
					=\nu(\calD)\Psi_n(\alpha),
					\]
					where \(F(\ba)\) denotes the Frobenius number \eqref{eq:DefFrobNum} and \(\Psi_n(\alpha)\) is given by \eqref{eq:DefComulativeFrob}.\end{proposition}
                We remark that for \(n=2\), the Frobenius numbers admit the formula \(F(a_1,a_2)=a_1a_2-a_1-a_2\), and the distribution of such numbers can be derived directly.

				\subsection{The projection-renormalization map}
                To prove Proposition~\ref{prop:PushforwardTheta}, we need to relate the Frobenius numbers with the image of orthogonal lattices in \(X_{n-1}(\R^n)\) via the work of \cite{Kannan1992}.
				
				Let \(\Pi:\R^n\to\R^{n-1}\) be the projection onto the first \(n-1\) coordinates.
                For every \(\bv\in\R_{>0}^n\), the restriction of \(\Pi\) to \(\bv^\perp\) is an isomorphism.
                 Denote
                \[
				\Gr_{n-1}^+(\R^n):=\set{\bv^\perp:\bv\in\R_{>0}^n}.
				\]
				Define the continuous map \(d_{\bullet}:\Gr_{n-1}^+(\R^n)\to\PGL_{n-1}(\R)\) by
				\[
				d_{\bv^\perp}:=\R^\times\diag(v_1,\dots,v_{n-1}),\quad \bv=(v_1,\dots,v_n)\in\R_{>0}^n.
				\]
				This is well defined, since replacing \(\bv\) by a positive multiple changes the diagonal matrix by a scalar.
				
				Let \[X_{n-1}^+(\R^n):=\pi^{-1}(\Gr_{n-1}^+(\R^n))=\set{[\Lambda]\in X_{n-1}(\R^n):\gen{\Lambda}_{\R}\in\Gr_{n-1}^+(\R^n)}.\]
				We obtain a continuous map 
				\begin{equation}\label{eq:DefTheProjNormMap}
					\Theta:X_{n-1}^+(\R^n)\to\PGL_{n-1}(\R)/\PGL_{n-1}(\Z),\;\;\Theta([\Lambda]):=d_{\gen{\Lambda}_{\R}}[\Pi(\Lambda)].
				\end{equation}

                The following is the key correspondence.
				\begin{proposition}\label{prop:KeyCorrespondneceFrobTheta}
					Let \(\ba\in\Zprim^n\cap\R_{>0}^n\). Then
					\[
					(\hat\chi\circ\Theta)([\Lambda_{\ba}])
					=\frac{F(\ba)+a_1+\dots+a_n}{(a_1\cdots a_n)^{1/(n-1)}},
					\]
					where \(F(\ba)\) is the Frobenius number \eqref{eq:DefFrobNum} and where \(\hat\chi\) is the normalized covering radius \eqref{eq:DefNormalizedCoveringRad}.
				\end{proposition}
				\begin{remark*}
For \(a_1,\dots,a_n\geq 2\), the covering-radius identity above
appears in a similar form  in \cite[Theorem~2.6]{Li2015}, building on
\cite{Kannan1992,AlievGruber2007}. Here we isolate the underlying
construction as the continuous map \(\Theta\) on
\(X_{n-1}^+(\R^n)\).
\end{remark*}
				We need the following preliminaries for the proof of Proposition~\ref{prop:KeyCorrespondneceFrobTheta}.
                
                For \(\tilde\ba:=(a_1,\dots,a_{n-1})\in\R_{>0}^{n-1}\), let
				\[
				\calS_{\tilde\ba}:=\set{\boldsymbol{x}\in\R_{\geq0}^{n-1}:\sum_{i=1}^{n-1}a_ix_i\leq1},
				\]
				and for a full-rank lattice, define the corresponding covering radius \(\tilde\Lambda\subset\R^{n-1}\), define
				\[
				\chi_{\tilde\ba}(\tilde\Lambda):=\inf\set{\alpha>0:\tilde\Lambda+\alpha\calS_{\tilde\ba}=\R^{n-1}}.
				\]
				
				Let \(\ba=(a_1,\dots,a_n)\in\Z^n\) such that \(a_n\neq0\). Then, 
				\[
				\Pi(\Lambda_{\ba})
				=\set{(x_1,\dots,x_{n-1})\in\Z^{n-1}:a_1x_1+\dots+a_{n-1}x_{n-1}\equiv0\pmod{a_n}},
				\]
				where \(\Lambda_{\ba}:=\ba^\perp\cap\Z^n\). Assuming further that \(\ba\in\Zprim^n\cap\R_{>0}^n\), we obtain by \cite[Theorem~2.5]{Kannan1992} the key identity
				\begin{equation}\label{eq:FrobNumerFormulaCoveringRad}
					\chi_{\tilde\ba}(\Pi(\Lambda_{\ba}))
					=F(a_1,\dots,a_n)+a_1+\dots+a_n.
				\end{equation}
				
				\begin{lemma}\label{lem:CovNormalizedProjectedLattice}
					Let \(\ba\in\Zprim^n\cap\R_{>0}^n\). Then
					\[
					\cov(d_{\tilde\ba}\Pi(\Lambda_{\ba}))=a_1\cdots a_n,
					\]
					where \(d_{\tilde\ba}:=\diag(a_1,\dots,a_{n-1})\).
				\end{lemma}
				
				\begin{proof}
					We first show that \(\cov(\Pi(\Lambda_{\ba}))=a_n\). Since the index of a sublattice of \(\Z^{n-1}\) equals its covolume, it is enough to prove that \(|\Z^{n-1}/\Pi(\Lambda_{\ba})|=a_n\).
					
					Consider the map \(\red_{\ba}:\Z^{n-1}\to\Z/(a_n\Z)\) defined by
					\[
					\red_{\ba}(x_1,\dots,x_{n-1})
					:=a_1x_1+\dots+a_{n-1}x_{n-1}\pmod{a_n}.
					\]
					Then \(\ker(\red_{\ba})=\Pi(\Lambda_{\ba})\). Since \(\ba\) is primitive, \(\red_{\ba}\) is onto. Therefore the index of \(\Pi(\Lambda_{\ba})\) in \(\Z^{n-1}\) is \(a_n\), and hence \(\cov(\Pi(\Lambda_{\ba}))=a_n\).
					
					Since \(\det(d_{\tilde\ba})=a_1\cdots a_{n-1}\), it follows that \(\cov(d_{\tilde\ba}\Pi(\Lambda_{\ba}))=a_1\cdots a_n\).
				\end{proof}
				
				\begin{proof}[Proof of Proposition~\ref{prop:KeyCorrespondneceFrobTheta}]
					Let \(\ba\in\Zprim^n\cap\R_{>0}^n\). By the definition of \(\hat\chi\circ\Theta\) and Lemma~\ref{lem:CovNormalizedProjectedLattice},
					\[
					(\hat\chi\circ\Theta)([\Lambda_{\ba}])
					=\frac{\chi(d_{\tilde\ba}\Pi(\Lambda_{\ba}))}
					{(a_1\cdots a_n)^{1/(n-1)}}.
					\]
					Observe that \(d_{\tilde\ba}\calS_{\tilde\ba}=\calS\). Consequently, for every full-rank lattice \(\tilde\Lambda\subset\R^{n-1}\),
					\begin{align*}
						\chi(d_{\tilde\ba}\tilde\Lambda)
						&=\inf\set{\alpha>0:d_{\tilde\ba}\tilde\Lambda+\alpha\calS=\R^{n-1}}\\
						&=\inf\set{\alpha>0:d_{\tilde\ba}(\tilde\Lambda+\alpha\calS_{\tilde\ba})=\R^{n-1}}\\
						&=\chi_{\tilde\ba}(\tilde\Lambda).
					\end{align*}
					Thus, by \eqref{eq:FrobNumerFormulaCoveringRad},
					\(\chi(d_{\tilde\ba}\Pi(\Lambda_{\ba}))=F(\ba)+a_1+\dots+a_n\), which proves the proposition.
				\end{proof}
				
				\subsection{Pushforwards: proof of Proposition~\ref{prop:PushforwardTheta}}
				For the following, we denote by \(\Omega_i\subset\Zprim^n\cap\R^n_{>0},\;i\in\N\) a sequence of finite subsets and for each \(i\) we let \(\mu_i\in\calP(X_{n-1}(\R^n))\) be the uniform probability measure supported on \(\set{[\Lambda_{\ba}]:\ba\in\Omega_i}\).
				\begin{lemma}\label{lem:RemovingCorrectionTerm}
					Assume that \(\lim_{i\to\infty}\mu_i=\mu\) in the weak-* topology, where \(\mu\) is a natural lift, and let \(\calD\subset\Gr_{n-1}^+(\R^n)\) be a Borel set satisfying \(\nu(\partial\calD)=0\), where \(\nu=\pi_*\mu\). Then, for every \(\epsilon>0\),
					\[
					\lim_{i\to\infty}
					\frac{\left|\set{\ba\in\Omega_i:
							\frac{a_1+\dots+a_n}{(a_1\cdots a_n)^{1/(n-1)}}\geq\epsilon,
							\ba^\perp\in\calD}\right|}{|\Omega_i|}=0.
					\]
				\end{lemma}
				
				\begin{proof}

					Identify \(\PGL_{n-1}(\R)/\PGL_{n-1}(\Z)\) with
					\(\hat X_{n-1}:=\SL_{n-1}(\R)/\SL_{n-1}(\Z)\), the space of unimodular lattices in \(\R^{n-1}\), by mapping \([\Lambda]\mapsto\cov(\Lambda)^{-1/(n-1)}\Lambda\).  By Lemma~\ref{lem:CovNormalizedProjectedLattice},
					\[
					L_{\ba}:=\frac{1}{(a_1\cdots a_n)^{1/(n-1)}}d_{\tilde\ba}\Pi(\Lambda_{\ba})
					\]
					is unimodular and represents \(\Theta([\Lambda_{\ba}])\) under this identification.
                    
					Given a full-rank lattice \(L\subset\R^{n-1}\), let
					\(L^*:=\set{x\in\R^{n-1}:x\cdot\ell\in\Z,\ \forall\ell\in L}\).
					We have
					\[
					L_{\ba}^*
					=(a_1\cdots a_n)^{1/(n-1)}
					d_{\tilde\ba}^{-1}\Pi(\Lambda_{\ba})^*.
					\]
					Notice that \(\Z^{n-1}\subset\Pi(\Lambda_{\ba})^*\) and that
					\(\frac{\tilde\ba}{a_n}\in\Pi(\Lambda_{\ba})^*\). Then,
					\[
					\frac{(a_1\cdots a_n)^{1/(n-1)}}{a_j}\be_j\in L_{\ba}^*
					\quad(1\leq j\leq n-1),
					\quad
					\frac{(a_1\cdots a_n)^{1/(n-1)}}{a_n}
					(\be_1+\dots+\be_{n-1})\in L_{\ba}^*.
					\]
					
					For \(1\leq j\leq n-1\) and \(\tau>0\), let
					\(\calM_{j,\tau}:=\set{\Lambda\in\hat X_{n-1}:\Lambda^*\cap(0,\tau]\be_j\neq\emptyset}\), and let
					\[
					\calM_{n,\tau}:=
					\set{\Lambda\in\hat X_{n-1}:
						\Lambda^*\cap(0,\tau](\be_1+\dots+\be_{n-1})\neq\emptyset}.
					\]
					If \(\frac{a_1+\dots+a_n}{(a_1\cdots a_n)^{1/(n-1)}}\geq\epsilon\), then
					\(\frac{a_j}{(a_1\cdots a_n)^{1/(n-1)}}\geq\frac{\epsilon}{n}\) for some \(1\leq j\leq n\). Consequently,
					\[
					\set{\ba\in\Omega_i:
						\frac{a_1+\dots+a_n}{(a_1\cdots a_n)^{1/(n-1)}}\geq\epsilon,
						\ba^\perp\in\calD}
					\subseteq
					\bigcup_{j=1}^n
					\set{\ba\in\Omega_i:
						L_{\ba}\in\calM_{j,n/\epsilon},
						\ba^\perp\in\calD}.
					\]
					
					We claim that \(\mu_{X_{n-1}}(\calM_{j,\tau})=0\) for every \(1\leq j\leq n\) and \(\tau>0\). The map \(\Lambda\mapsto\Lambda^*\) is a measure-preserving homeomorphism of \(\hat X_{n-1}\). Indeed, if \(\Lambda=g\Z^{n-1}\), then \(\Lambda^*=g^*\Z^{n-1}\), where \(g^*:=(g^{-1})^\top\), and \(g\mapsto g^*\) is an automorphism of \(\SL_{n-1}(\R)\) preserving \(\SL_{n-1}(\Z)\) and Haar measure. Hence, by Siegel's summation formula \cite{Siegel1945}, for \(1\leq j\leq n-1\),
					\[
					\begin{aligned}
						\mu_{X_{n-1}}(\calM_{j,\tau})
						&=\mu_{X_{n-1}}\set{\Lambda\in\hat X_{n-1}:
							\Lambda\cap(0,\tau]\be_j\neq\emptyset}\\
						&\leq\int_{\hat X_{n-1}}
						\#\bigl(\Lambda\cap(0,\tau]\be_j\bigr)\,
						d\mu_{X_{n-1}}(\Lambda)\\
						&=\vol_{n-1}\bigl((0,\tau]\be_j\bigr)=0.
					\end{aligned}
					\]
					The same argument with \(\be_1+\dots+\be_{n-1}\) in place of \(\be_j\) shows that \(\mu_{X_{n-1}}(\calM_{n,\tau})=0\).
					
					The sets \(\calM_{j,\tau}\) are closed. Indeed, suppose that \(\set{\Lambda_r}_{r\in \N}\subset \calM_{j,\tau}\) is a sequence converging to \(\Lambda\), where for each \(r\in \N\), \(t_r\be_j\in\Lambda_r^*\) with \(0<t_r\leq\tau\). After passing to a subsequence, \(t_r\to t\in[0,\tau]\). Since \(\Lambda_r^*\to\Lambda^*\), the lengths of the shortest nonzero vectors of \(\Lambda_r^*\) are bounded away from zero. Hence \(t>0\), and \(t\be_j\in\Lambda^*\). The same argument applies to \(\calM_{n,\tau}\).
					
					Since \(\mu\) is a natural lift, \(\Theta_*\mu_{\V}=\mu_{X_{n-1}}\) for \(\nu\)-almost every \(\V\in\Gr_{n-1}^+(\R^n)\). For fixed \(j\) and \(\tau\), the closure in \(X_{n-1}(\R^n)\) of
					\[
					\set{[\Lambda]\in X_{n-1}^+(\R^n):
						\pi([\Lambda])\in\calD,\ \Theta([\Lambda])\in\calM_{j,\tau}}
					\]
					is contained in
					\[
					\set{[\Lambda]\in X_{n-1}^+(\R^n):
						\pi([\Lambda])\in\overline{\calD},
						\Theta([\Lambda])\in\calM_{j,\tau}}
					\cup\pi^{-1}(\partial\calD).
					\]
					The first set has \(\mu\)-measure zero by disintegration and using
					\(\Theta_*\mu_{\V}=\mu_{X_{n-1}}\), while the second has \(\mu\)-measure
					\(\nu(\partial\calD)=0\). The Portmanteau theorem therefore gives
					\[
					\limsup_{i\to\infty}
					\frac{\left|\set{\ba\in\Omega_i:
							L_{\ba}\in\calM_{j,\tau},
							\ba^\perp\in\calD}\right|}{|\Omega_i|}=0.
					\]
					Applying this to \(\tau=n/\epsilon\) and summing over \(1\leq j\leq n\) proves the lemma.
				\end{proof}
				
				\begin{proof}[Proof of Proposition~\ref{prop:PushforwardTheta}]
					Assume that \(\lim_{i\to\infty}\mu_i=\mu\) in the weak-* topology, where \(\mu\) is a natural lift. We first show that
					\[
					\lim_{i\to\infty}
					\mu_i\set{[\Lambda]\in X_{n-1}^+(\R^n):
						\pi([\Lambda])\in\calD,
						(\hat\chi\circ\Theta)([\Lambda])\geq\alpha}
					=\nu(\calD)\Psi_n(\alpha).
					\]
					
					By \cite[Sections~5--6]{Marklof2010}, \(\hat\chi\) is continuous. Hence \(\hat\chi\circ\Theta\) is continuous on \(X_{n-1}^+(\R^n)\). Since \(\mu\) is a natural lift, uniqueness of the invariant probability measure gives \(\Theta_*\mu_{\V}=\mu_{X_{n-1}}\) for \(\nu\)-almost every \(\V\in\Gr_{n-1}^+(\R^n)\).
					
					Let
					\[
					A_\alpha:=\set{[\Lambda]\in X_{n-1}^+(\R^n):
						\pi([\Lambda])\in\calD,
						(\hat\chi\circ\Theta)([\Lambda])\geq\alpha}.
					\]
					By disintegration,
					\[
					\begin{aligned}
						\mu(A_\alpha)
						&=\int_{\calD}
						\mu_{\V}\set{[\Lambda]\in X(\V):
							(\hat\chi\circ\Theta)([\Lambda])\geq\alpha}\,d\nu(\V)\\
						&=\nu(\calD)\Psi_n(\alpha).
					\end{aligned}
					\]
					
					By \cite[Section~7, Lemma~7]{Marklof2010},
					\(\mu_{X_{n-1}}\set{\hat\chi=\alpha}=0\) for \(\alpha\in\R\). Continuity of \(\pi\), \(\Theta\), and \(\hat\chi\) gives
					\[
					\partial A_\alpha
					\subseteq
					\pi^{-1}(\partial\calD)
					\cup
					\set{[\Lambda]\in X_{n-1}^+(\R^n):
						\pi([\Lambda])\in\overline{\calD},
						(\hat\chi\circ\Theta)([\Lambda])=\alpha}.
					\]
					Both sets on the right have \(\mu\)-measure zero. Hence \(\mu(\partial A_\alpha)=0\), and the Portmanteau theorem gives
					\[
					\lim_{i\to\infty}\mu_i(A_\alpha)
					=\mu(A_\alpha)
					=\nu(\calD)\Psi_n(\alpha).
					\]
					By Proposition~\ref{prop:KeyCorrespondneceFrobTheta}, this is equivalent to
					\[
					\lim_{i\to\infty}
					\frac{\left|\set{\ba\in\Omega_i:
							\frac{F(\ba)+a_1+\dots+a_n}
							{(a_1\cdots a_n)^{1/(n-1)}}\geq\alpha,
							\ba^\perp\in\calD}\right|}{|\Omega_i|}
					=\nu(\calD)\Psi_n(\alpha).
					\]
					The proposition follows from Lemma~\ref{lem:RemovingCorrectionTerm}.
				\end{proof}
				
				\subsection{Proofs of Theorem~\ref{thm:frob_on_quad_var} and \ref{thm:1.2}}\label{sec:proof-of-shap-bers-applic} 
				Let us first prove Theorem~\ref{thm:frob_on_quad_var} as a consequence of Theorem 4.8 of \cite{BersudskyShapira2023} and Proposition~\ref{prop:PushforwardTheta}. 
                
                For a compact subset \(\calD \) of \(\set{\bx\in \R^d_{>0} : Q(\bx)=1}\) and a sequence $\{N_i\}_{i=1}^\infty$ satisfying the assumptions of Theorem~\ref{thm:frob_on_quad_var}, let
				\[ \Omega_i := \sqrt{N_i}\calD \cap \Z^d_{\text{prim}}.\]
				Our objective is to show that the sequence of probability measures
				\[\mu_i := \frac{1}{|\Omega_i|} \sum_{\ba \in \Omega_i} \delta_{[\Lambda_{\ba}]}\]
				is eventually well defined and converges in the weak-* topology to a natural lift supported on $X^+_{d-1}(\R^d)$. Once this is established, we apply Proposition~\ref{prop:PushforwardTheta} with the set \(\calD\) appearing there taken to be the whole space \(\Gr^+_{d-1}(\R^d),\) which completes the proof of Theorem~\ref{thm:frob_on_quad_var}.
				
				We first check that the  assumptions of Theorem 4.8 of \cite{BersudskyShapira2023} are satisfied for the sequence $\{N_i\}_{i=1}^\infty.$ Both $Q_1$ and $Q_2$ are integral, non-degenerate quadratic forms in $d\geq4$ variables such that $Q(\be_d)>0$ and $\mathrm{Stab}_{\SO_Q(\R)}(\be_d)\cong \SO_{d-1}(\R)$, and therefore satisfy the standing assumption in \cite[p.581]{BersudskyShapira2023}. When $Q=Q_2$,  Section 4.5 of \cite{BersudskyShapira2023} contains the proof that every $N_i\in \N$ has the $(Q_2,5)$ co-isotropic property. If $Q=Q_1,$ then any natural number $N_i$ can be represented as $Q(\bv)=N_i$ with $\bv\in\Z_{\text{prim}}^d$ unless $d=4$ and $8\mid N_i$ (see, e.g. \cite{CooperHirschhorn}). A similar argument to the case $Q=Q_2$ shows that every such $N_i$ has the $(Q_1,5)$ co-isotropic property.
				
				We now adopt some notation from \cite{BersudskyShapira2023}. Here, the ambient space is
				\[\mathcal{X}(\R):=\set{(\Lambda,\bv): \Lambda \text{ is a } (d-1)\text{-lattice in } \R^d,\ \bv\in\R^d,\ \bv \perp \Lambda,\ \mathrm{covol}(\Lambda)=\|\bv\|}\]
				as in \cite[p.592]{BersudskyShapira2023}.
				Following the definitions in \cite[p.598]{BersudskyShapira2023}, denote
				\[\calH_N(\R):= \{\bv \in \R^d: Q(\bv)=N\},\quad \calH_{N,\text{prim}}(\Z) := \{\ba \in \Zprim^d: Q(\ba)=N\},\]
				and 
				\[\mathcal{X}_N(\R):= \{(\Lambda,\bv)\in\mathcal{X}(\R): \bv\in \calH_N(\R)\}, \quad \mathcal{X}_N(\Z):=\{(\Lambda_{\ba},\ba): \ba \in \calH_{N,\text{prim}}(\Z)\}. \]
			Here,  $\Lambda_{\ba}=\ba^\perp \cap \Z^d.$ Since $Q(\be_d)=1,$ the scaling map $\pi_{\mathcal{X}_N}:\mathcal{X}_N(\R) \to \mathcal{X}_1(\R)$ defined in \cite[p.599]{BersudskyShapira2023} is given by
				\[\pi_{\mathcal{X}_N}(\Lambda, \bv) = \left(N^{-\frac1{2(d-1)}}\Lambda, N^{-\frac12}\bv\right).\]
				
				We now introduce the main objects appearing in Theorem 4.8 of \cite{BersudskyShapira2023}. Following \cite[p.604]{BersudskyShapira2023}, define
				\[h_i:=|\calH_{N_i,\text{prim}}(\Z)/\SO_Q(\Z)|,\quad \nu_i:=\frac1{h_i}\sum_{\ba\in\calH_{N_i,\text{prim}}(\Z)}\delta_{\pi_{\mathcal{X}_{N_i}}(\Lambda_{\ba},\ba)},\]
				where $h_i$ is finite by \cite[p.579]{BersudskyShapira2023} and positive since $\calH_{N_i,\text{prim}}(\Z)\neq \emptyset.$ Denote the projection map $\pi_{vec}:\mathcal{X}_1(\R) \to \calH_1(\R)$ as $\pi_{vec}(\Lambda,\bv)=\bv$, and from Section 4.3.2 of \cite{BersudskyShapira2023}, define
				\[\mu_{\mathcal{X}}:=\int_{\calH_1(\R)} \mu_{(\pi_{vec})^{-1}(\bv)} \dd\mu_{\calH}(\bv),\]
				where the fiber measure $\mu_{(\pi_{vec})^{-1}(\bv)}$ is induced by the invariant measure on $X_{d-1}$, and $\mu_{\calH}=\mu_{\calH_{Q(\be_d)}(\R)}$ is the $\SO_Q(\R)$-invariant measure on $\calH_1(\R)$ such that a measurable fundamental domain in $\calH_1(\R)$ for $\SO_Q(\Z)$ has measure equal to one.
				
				By Theorem 4.8 of \cite{BersudskyShapira2023}, we have
				\[\nu_i(f) \to\mu_{\mathcal{X}}(f), \quad \forall f \in C_c(\mathcal{X}_1(\R)).\]
				Let $E_{\calD}:= (\pi_{vec})^{-1}(\calD).$ Since $\mu_{\calH}(\partial\calD)=0$,
				\[\nu_i|_{E_{\calD}}(f) \to\mu_{\mathcal{X}}|_{E_{\calD}}(f), \quad \forall f \in C_c(\mathcal{X}_1(\R)).\]
				Here,
				\[\nu_i|_{E_{\calD}} = \frac1{h_i} \sum_{\ba \in \Omega_i}\delta_{\pi_{\mathcal{X}_{N_i}}(\Lambda_{\ba},\ba)} \quad \text{ and } \quad \mu_{\mathcal{X}}|_{E_{\calD}}= \int_{\calD} \mu_{(\pi_{vec})^{-1}(\bv)} \dd\mu_{\calH}(\bv),\]
				since $\Omega_i=\mathcal{H}_{N_i,\text{prim}}(\Z)\cap \sqrt{N_i}\calD.$
				
				We now translate the above result into our setting. Define a proper map
				\[r:E_{\calD} \to X_{d-1}(\R^d), \quad r(\Lambda,\bv)=[\Lambda].\]
				Then we have
				\[r_* (\nu_i|_{E_{\calD}}) =\frac1{h_i}\sum_{\ba \in \Omega_i}\delta_{[\Lambda_{\ba}]} \quad  \xrightarrow{weak-*} \quad r_* (\mu_{\mathcal{X}}|_{E_{\calD}})= \int_{\calD}\mu_{\bv^\perp} \ \dd\mu_{\calH}(\bv),\]
				where $\mu_{\bv^\perp}$ is the $\PGL_{d-1}(\R)$-invariant probability measure on $X_{d-1}(\bv^\perp).$
				
				We claim that
				\begin{equation}\label{eq:ratio_omegai_hi}
					\lim_{i\to\infty}\frac{|\Omega_i|}{h_i} = \mu_{\calH}(\calD).
				\end{equation}
				Once this is established, $\mu_i=\frac{h_i}{|\Omega_i|}r_* (\nu_i|_{E_{\calD}})$ is well defined for big enough $i$ and  
				\[\mu_i \xrightarrow{weak-*} \frac{1}{\mu_{\calH}(\calD)}\int_{\calD} \mu_{\bv^\perp} \ \dd\mu_{\calH}(\bv).\]
				Here, the right side is a probability measure with $\PGL_{d-1}(\R)$-invariant measure on each fiber, therefore it is a natural lift, which completes the proof of Theorem \ref{thm:frob_on_quad_var}.
				
				It remains to prove \eqref{eq:ratio_omegai_hi}. Choose a measurable fundamental domain $\calF$ in $\calH_1(\R)$ for $\SO_Q(\Z)$ with $\mu_{\calH}(\partial\calF)=0$ and $|\calF\cap \bv\SO_Q(\Z)|=1$ for every $\bv \in \calH_1(\R).$ For $E_{\calF}=(\pi_{vec})^{-1}(\calF),$
				\[\nu_i|_{E_{\calF}}(f) \to\mu_{\mathcal{X}}|_{E_{\calF}}(f), \quad \forall f \in C_c(\mathcal{X}_1(\R)).\]
				By definition of $h_i, \calF$ and $\mu_{\calH},$ for all $i,$
				\[\nu_i(E_{\calF})=\frac{\left|\left\{\ba \in \calH_{N_i,\text{prim}(\Z)}:\frac{\ba}{\sqrt{N_i}}\in \calF\right\}\right|}{|\calH_{N_i,\text{prim}(\Z)}/\SO_Q(\Z)|}=1 \quad \text{ and } \quad \mu_{\mathcal X}(E_{\calF})=\mu_{\calH}(\calF)=1.\]
				Therefore, $\nu_i|_{E_{\calF}}$ converges weakly to $\mu_{\mathcal X}|_{E_{\calF}},$ and thus $\pi_{vec}$ pushforwards the convergence
				\[(\pi_{vec})_* (\nu_i|_{E_{\calF}})= \frac1{h_i}\sum_{\ba \in \calH_{N_i,\text{prim}}(\Z)\cap \sqrt{N_i}\calF} \delta_{\frac{\ba}{\sqrt{N_i}}} \quad \xrightarrow{weak} \quad (\pi_{vec})_* (\mu_{\mathcal X}|_{E_{\calF}}) = \mu_{\calH}|_{\calF}.\]
				For $f\in C_c(\calH_{1}(\R)),$ if we let
				\[\hat{f}(\bv) := \sum_{\mathbf{w} \in \bv \SO_Q(\Z)} f(\mathbf{w}),\]
				then $\hat{f}$ is bounded and continuous $\mu_{\calH}$-almost everywhere. Hence,
				\[\int_{\calF} \hat{f} \ \dd(\pi_{vec})_* (\nu_i|_{E_{\calF}})  = \frac{1}{h_i}\sum_{\ba \in \calH_{N_i,\text{prim}}(\Z)} f({\textstyle \frac{\ba}{\sqrt{N_i}}})\quad \to\quad \int_{\calF} \hat f \ \dd \mu_{\calH}=\int_{\calH}f \ \dd \mu_{\calH}.\]
				Since $\calD$ is compact and $\mu_{\calH}(\partial\calD)=0$, the Portmanteau theorem gives \eqref{eq:ratio_omegai_hi}. \qed

                \vspace{ .1 in}
                \noindent \textbf{Proof of Theorem~\ref{thm:1.2}.}
				Now let us prove Theorem~\ref{thm:1.2} as a consequence of Theorems~\ref{thm:PrincipalSL2Equid},~\ref{thm:mainthm_g2sppgl} and Proposition~\ref{prop:PushforwardTheta}. Define
                \[\Omega_T:= \{\bv(a,b,c,d): a,b,c,d\in\N,\ ad-bc=1, \ \|(a,b,c,d)\|_p\leq T
                \},\]
				where $\|(a,b,c,d)\|_p=(|a|^p + |b|^p + |c|^p + |d|^p)^{1/p},$ and let
                \[\mu_T:=\frac{1}{|\Omega_T|} \sum_{\bv \in \Omega_T} \delta_{[\Lambda_{\bv}]}.\]
                We will show that $\mu_T$ converges in the weak-* topology to a natural lift measure on $X_n^+(\R^{n+1}),$ and that $|\Omega_T| \sim c_pT^2$. Proposition~\ref{prop:PushforwardTheta} with $\calD=\Gr_{n}^+(\R^{n+1})$ then yields
                \[\lim_{T\to \infty} \frac{\#\left\{\bv(a,b,c,d) \in \Omega_T: \frac{F(\bv(a,b,c,d))}{R(a,b,c,d)^{1/n}}\geq \alpha\right\}}{c_pT^2} = \Psi_{n+1}(\alpha),\]
                which proves Theorem~\ref{thm:1.2}. 
                
                In order to establish the desired weak-* convergence, we use Theorems~\ref{thm:PrincipalSL2Equid},~\ref{thm:mainthm_g2sppgl} with $\Gamma=\SL(2,\Z)$ which gives that probability measures
                \[\frac{1}{|\Gamma_T|} \sum_{\gamma \in \Gamma_T} \delta_{\gamma \cdot [\Lambda_0]}\]
                converge weak-* to a natural lift measure on $X_n(\VV_n),$ provided that $F_{[\Lambda_0]}=\PGL_n(\R)^\circ$. Thus, we need to relate $\Gamma_T$-orbit measures on $X_n(\VV_n)$ to $\mu_T$ on $X_n(\R^{n+1}).$

                Consider the basis of  $\VV_n$ given by
                \[p_i(x,y)=x^{n+1-i}y^{i-1},\quad i=1,\dots,n+1,\]
                and denote by $\gen{\cdot,\cdot}$ the unique inner product on \(\VV_n\) such that $p_1,\dots,p_{n+1}$ are orthonormal. 
                
                For $g\in\GL(\VV_n),$ define $g^\top$ to be the adjoint transformation satisfying
                \begin{equation}\label{eq:TransposeDef}
                    \gen{gv_1, v_2} = \gen{v_1,g^\top v_2},\quad v_1,v_2 \in \VV_n.
                \end{equation}
               We denote
               \begin{equation*}
                   g^*:=(g^\top)^{-1},\;\;g\in\GL(\VV_n),
               \end{equation*}
                and we observe that
                \begin{equation}\label{eq:InvarianceOfInnerProduct}
                    \gen{v_1,v_2}=\gen{gv_1,g^*v_2},\;\forall g\in \GL(\VV_n).
                \end{equation}
                For \(v\in \VV_n\), we denote \(v^\perp:=\set{p\in \VV_n:\,\gen{p,v}=0}.\) Then, by \eqref{eq:InvarianceOfInnerProduct},
                \begin{equation}\label{eq:PreservationOfOrth}
                    g v^\perp = (g^*v)^\perp.
                \end{equation}
                
                We define the irreducible representation
                \begin{equation}\label{eq:DualActioOnVecs}
                    \varrho^*:\GL_2(\R) \to \GL(\VV_n),\;\;\varrho^*(h)=(h\cdot )^*,\;\;h\in\GL_2(\R),
                \end{equation}
                where \(h\cdot\) denotes the resulting map in \(\GL(\VV_n)\) obtained via \eqref{eq:DefOfIrredAction}. 
                
                Clearly, 
                 $\SL(2,\Z)\cdot \Sym^n(\Z^2)=\Sym^n(\Z^2)$, where  $\Sym^n(\Z^2):=\gen{p_1,\dots,p_{n+1}}_\Z.$ Therefore, by \eqref{eq:PreservationOfOrth}, we obtain that \(\gamma\cdot \left(\Sym^n(\Z^2) \cap {\bv}^\perp\right) = \Sym^n(\Z^2) \cap (\varrho^*(\gamma)\bv)^\perp\) for all \(\bv \in \VV_n\) and \(\gamma\in\SL_2(\Z)\). In particular,
                \begin{equation}\label{eq:DualAction}
                    \gamma \cdot\Lambda_{\bv}=\Lambda_{\varrho^*(\gamma)\bv},\;\bv\in\Sym^n(\Z^2),\;\gamma\in\SL_2(\Z).
                \end{equation}
            
                Let $$\bv_0:=p_1+p_{n+1}=x^n+y^n.$$
                In the following, for $h\in \SL_2(\R)$, we denote by \(h^\top\) the usual matrix transpose. Let \(\gamma=\begin{psmallmatrix}
                    a&b\\c&d
                \end{psmallmatrix}\in\SL_2(\Z)\). Then, the $i$-th coefficient of $\varrho^*\left((\gamma^\top)^{-1}\right)\bv_0$ is given by
                \begin{align*}
                    \gen{p_i,\varrho^*\left((\gamma^\top)^{-1}\right)\bv_0} &= \gen{\gamma^\top \cdot p_i,\bv_0} = \gen{(ax+by)^{n+1-i}(cx+dy)^{i-1},x^n+y^n} \\
                    &= a^{n+1-i}c^{i-1} + b^{n+1-i}d^{i-1}.
                \end{align*}
                By defining the coefficient map $\iota:\VV_n \to \R^{n+1}$ as $\iota(\sum a_ip_i)=(a_1,\dots,a_{n+1})$ and the induced map $\iota:X_n(\VV_n) \to X_n(\R^{n+1}),$ we conclude that
                \[\bv(a,b,c,d) = \iota\left( \varrho^*\left((\gamma^\top)^{-1}\right)\bv_0\right) \in\Z_{\text{prim}}^{n+1}, \quad [\Lambda_{\bv(a,b,c,d)}] = \iota \left((\gamma^\top)^{-1}\cdot [\Lambda_{\bv_0}]\right).\]
                \begin{remark}\label{rem:GeomSeqAndSpecialPoints}
                    Using the above calculation, we see that
                    \[\iota\paren{\varrho^* \begin{psmallmatrix}
                        d&-c\\-b&a
                    \end{psmallmatrix} x^n}=(a^n,a^{n-1}c,\dots,c^n).\]
                    As mentioned in the introduction, the Frobenius numbers of such sequences have elementary formulae \cite{OngPonomarenko2008Frobenius}. We observe that the orthogonal lattice of \(x^n\) is \(\gen{x^{n-1}y,x^{n-2}y^2,\dots,y^n}_\Z\), which is \(\SL_2(\Z)\)-special, see Definition~\ref{def:GammaSpecialPoints}. Notice that the corresponding norm-ball averages converge to the measure on the \(\SL_2(\Z)\)-equivariant section, see Theorem~\ref{thm:MainTheoremGammaSpecial}.
                \end{remark}

                We now show that for even $n\geq4,$ the hyperplane $\gen{\Lambda_{\bv_0}}_{\R}$ is not an $F$-subspace for $F=\PSp(n,\R)$, which implies that $F_{[\Lambda_{\bv_0}]}=\PGL(n,\R)^\circ.$ Observe that
                \[\Lambda_{\bv_0} = \gen{p_1-p_{n+1},p_2,p_3,\dots,p_n}_\Z.\] By Proposition~\ref{prop:ExistanceOfStandardBasis}~\ref{enu:BiLinearObstruction}, it suffices to show that there is no nondegenerate symplectic form $B$ on $\R^n$ satisfying
                \begin{equation}\label{eq:IsotropyInPRoofOfApplication}
                    B(P^{(i-1)}(1,s), P^{(j-1)}(1,s)) = 0, \quad \forall 1\leq i,j\leq \frac{n}2,\ s \in \R,
                \end{equation}
                where $P(1,s)=\left(p_1(1,s)-p_{n+1}(1,s), p_2(1,s), \dots,p_n(1,s)\right) = (1-s^n,s,s^2,\dots,s^{n-1}).$ Assume for a contradiction that such $B$ exists. Denote by $\be_1,\dots,\be_n$ the standard basis for $\R^n.$ Setting $s=0,$ we have $P^{(i-1)}(1,0)=(i-1)!\be_i$ for $i<n,$ and thus
                \[B(\be_1,\be_j)=0, \quad \text{ for } \textstyle 1\leq j \leq \frac{n}2.\]
                We also have $B(\be_1, \be_{\frac{n}2+1})=0,$ since 
                \begin{align*}\frac{\dd}{\dd s}
                B(P(1,s),P^{(\frac{n}2-1)}(1,s))&=B(P^{(1)}(1,s),P^{(\frac{n}2-1)}(1,s))+B(P(1,s),P^{(\frac{n}2)}(1,s))\\&=B(P(1,s),P^{(\frac{n}2)}(1,s))=0.\end{align*}
                 For \(s\neq 0,\) define
                \[Q(s):=s^n P(1,1/s) = (s^n-1, s^{n-1},s^{n-2},\dots,s),\]
                and extend \(Q(s)\) and its derivatives to $s=0.$ Then $Q(0)=-\be_1$ and $Q^{(i)}(0)=i!\be_{n+1-i}$ for $1\leq i <n$. Since $Q^{(i)}(s)$ is a linear combination of  $f_0(s)P(1,1/s), \dots, f_{i-1}(s)P^{(i)}(1,1/s),$ we conclude by \eqref{eq:IsotropyInPRoofOfApplication}
                \[B(Q(s),Q^{(i)}(s))=0,\;\forall 1\leq i<\frac n2,\;s\in \R.\]
                Consequently, \(B(\be_1, \be_j) = 0\) for \(\frac{n}2 +2\leq j\leq n.\) Therefore, $B(\be_1,\be_j)=0$ for every $1\leq j \leq n,$ which contradicts the nondegeneracy of $B.$
                
                When $n=2,$ the subspace $\gen{x^2-y^2,xy}_\R$ is not Maximally Divisible. Therefore, $F_{[\Lambda_{\bv_0}]}=\PGL^\circ_n(\R)$ for any even $n\geq2,$ and thus Theorems~\ref{thm:PrincipalSL2Equid},~\ref{thm:mainthm_g2sppgl} give that
                \begin{equation}\label{eq:weakconv_proof_sl2rfrob}
                    \frac{1}{\#\Gamma_T} \sum_{\gamma \in \Gamma_T} \delta_{\gamma \cdot [\Lambda_{\bv_0}]} \xrightarrow{weak-*} \int_0^{\pi} \kappa(\theta)_*\mu_{\V^+_n} w(\theta)\dd\theta,
                \end{equation}
                where $\mu_{\V_n^+}$ is $\PGL(n,\R)$-invariant probability measure on $X_n(\V_n^+)$ for $\V_n^+=\gen{p_1,\dots,p_n}_\R.$

                We now analyze $\Gamma_T.$ Since the norm $\|\cdot\|_p$ on $\SL(2,\R)$ is preserved by transpose and inverse,
                \[\frac{1}{\#\Gamma_T} \sum_{\gamma \in \Gamma_T} \delta_{\gamma \cdot [\Lambda_{\bv_0}]}=\frac{1}{\#\Gamma_T} \sum_{\gamma \in \Gamma_T} \delta_{(\gamma^{-1})^\top \cdot [\Lambda_{\bv_0}]}.\]
                Define $J=\begin{psmallmatrix}0&-1\\1&0\end{psmallmatrix}$ and let
                \[\Gamma_T^0:=\left\{\begin{psmallmatrix} a&b \\ c&d \end{psmallmatrix} \in \Gamma_T : abcd=0\right\}, \quad \Gamma_T^+:=\left\{\begin{psmallmatrix} a&b \\ c&d \end{psmallmatrix} \in \Gamma_T : a,b,c,d \in \N \right\}.\]
                Then $\Gamma_T$ can be decomposed as
                \begin{equation}\label{eq:decomp_Gamma}
                \Gamma_T= \Gamma_T^0 \sqcup \widetilde\Gamma_T^+ \sqcup J(\widetilde\Gamma_T^+)\quad \text{ where } \quad  \widetilde\Gamma_T^+ = \Gamma_T^+ \sqcup (-\Gamma_T^+) \sqcup (\Gamma_T^+J) \sqcup (-\Gamma_T^+J).
                \end{equation}
                Here, each right translates consisting $\widetilde\Gamma_T^+$ has same orbit because $\varrho(J)$ and $\varrho(-I)$ fix $\Lambda_{\bv_0}.$ Additionally, the coefficients of $\varrho^* \left((\gamma^{-1})^\top\right) \bv_0$ are all positive for $\gamma \in \widetilde\Gamma_T^+$ while their signs alternate for $\gamma \in J(\widetilde \Gamma_T^+).$ Therefore, for $\gamma \in \Gamma_T\smallsetminus \Gamma_T^0$, $\iota\left(\left(\varrho^* \left((\gamma^{-1})^\top\right) \bv_0 \right)^\perp\right) \in \Gr^+_n(\R^{n+1})$ if and only if $\gamma \in \widetilde\Gamma_T^+,$ and thus for $X=\iota^{-1}(X^+_n(\R^{n+1})),$
                \[\frac1{\#(\Gamma_T\smallsetminus\Gamma_T^0)} \left.\left(\sum_{\gamma \in \widetilde\Gamma_T^+} \delta_{(\gamma^{-1})^\top \cdot [\Lambda_{\bv_0}]} + \sum_{\gamma \in J(\widetilde\Gamma_T^+)} \delta_{(\gamma^{-1})^\top \cdot [\Lambda_{\bv_0}]}\right)\right|_X = \frac1{2\#\Gamma_T^+} \sum_{\gamma \in \Gamma_T^+} \delta_{(\gamma^{-1})^\top \cdot [\Lambda_{\bv_0}]}.\]
                Observe also that
                \[\kappa(\theta) \cdot \V_n^+ = (\varrho^*(\kappa(\theta))p_{n+1})^\perp = ((-\sin\theta)^nx^n + (-\sin\theta)^{n-1} (\cos\theta) x^{n-1}y + \cdots + (\cos\theta)^ny^n)^\perp\]
                so that $\iota(\kappa(\theta) \cdot \V_n^+) \in \Gr_n^+(\R^{n+1})$ if and only if $\frac\pi2 < \theta < \pi.$
                Since $\frac{\# \Gamma_T^0}{\# \Gamma_T}\to0$ as $T\to \infty$ and $\partial(\iota^{-1}(X_n^+(\R^{n+1})))$ has measure zero, restricting \eqref{eq:weakconv_proof_sl2rfrob} gives
                \[ \frac1{2\#\Gamma_T^+} \sum_{\gamma \in \Gamma_T^+} \delta_{(\gamma^{-1})^\top \cdot [\Lambda_{\bv_0}]} \xrightarrow{weak-*} \int_{\pi/2}^{\pi} \kappa(\theta)_* \mu_{\V_n^+} w(\theta)\dd\theta.\]
                
                Since $\mathrm{Stab}^{\varrho^*}_{\SL(2,\Z)}(\bv_0)=\{\pm I, \pm J\},$ the map $\Gamma_T^+ \to \Omega_T: \begin{psmallmatrix} a&b \\ c&d \end{psmallmatrix} \mapsto \bv(a,b,c,d)$ is bijective.  Indeed, if  \(h=\begin{psmallmatrix}a&b\\c&d\end{psmallmatrix}\in\mathrm{SL}_2(\Z)\) satisfies \(\varrho^*((h^{-1})^\top)\bv_0=\bv_0\), then comparing the first and last coefficients gives \(a^n+b^n=c^n+d^n=1\). Since \(n\) is even, \(h\) is a determinant-one signed permutation matrix. Hence,
                \[\mu_T=\frac{1}{|\Omega_T|} \sum_{\bv \in \Omega_T} \delta_{[\Lambda_{\bv}]} =  \iota_*\left(\frac1{\#\Gamma_T^+} \sum_{\gamma \in \Gamma_T^+} \delta_{(\gamma^{-1})^\top \cdot [\Lambda_{\bv_0}]}\right).\]
                Noting that $\int_{\pi/2}^{\pi} w(\theta)\dd\theta=\frac12,$ we conclude that $\mu_T$ converges to a natural lift.

                We now calculate $c_p.$ By  e.g. \cite{GorodnikNevo},
                \[\lim_{T \to \infty}\frac{\vol(G_T)}{\#\Gamma_T} = \vol(G/\Gamma).\]
                Corollary~\ref{cor:volume-asymp-SL2} gives the asymptotic
                \[\lim_{T\to \infty}\frac{\#\Gamma_T}{T^2} = \frac{1}{\vol(G/\Gamma)} \cdot 2 \int_0^\pi \int_{-\infty}^\infty \|\kappa(\theta)\delta(\infty)\nu(s) \|_p^{-2} \ \dd s \dd \theta.\]
                Since $\kappa(\theta)\delta(\infty)\nu(s) = \begin{psmallmatrix}\cos\theta & s\cos\theta \\ \sin\theta & s\sin\theta \end{psmallmatrix},$ if we denote 
                \[I_p:=\int_0^\infty (1+s^p)^{-2/p}\dd s= \frac1p \beta \left(\frac1p,\frac1p\right) =\frac{\Gamma(1/p)^2}{p\Gamma(2/p)},\]
                where $\beta$ and $\Gamma$ denote the classical Beta and Gamma functions, then
                \[\int_0^\pi \int_{-\infty}^\infty \|\kappa(\theta)\delta(\infty)\nu(s) \|_p^{-2} \ \dd s \dd \theta = \int_0^\pi (|\cos\theta|^p+|\sin\theta|^p)^{-2/p} \dd \theta \int_{-\infty}^\infty (1+|s|^p)^{-2/p}\dd s = 4I_p^2. \]
                We recall that, using the standard fundamental domain for \(\SL_2(\Z)\) in \(\SL_2(\R)\), 
                \[\vol(G/\Gamma) = \int_0^\pi \int_{-\frac12}^{\frac12}\int_{-\infty}^{-\frac12 \log(1-s^2)}e^t\ \dd t \dd s \dd \theta = \frac{\pi^2}3.\]
                Therefore, as $\#\Gamma_T = \#\Gamma_T^0 + 8\#\Gamma_T^+$ from \eqref{eq:decomp_Gamma} and $\lim_{T\to\infty}\frac{\#\Gamma_T^0}{T^2}=0,$ we have
                \[c_p = \frac18 \cdot \frac3{\pi^2} \cdot 8|I_p|^2 = \frac{3\Gamma(1/p)^4}{\pi^2p^2\Gamma(2/p)^2}.\]
                \qed

					\appendix
					
					\section{Proving Proposition~\ref{prop:AdmissibleLieAlgebra}}\label{sec:DynkingAppend}
                    Below, \(\upE_{k-1},\,\upH_{k-1},\,\upS_{k-1}\) are as in \eqref{eq:DefPrincSL2Trip}.
                    
					The following simple observation, which we leave the reader to verify, will be used repeatedly. 
					\begin{lemma}\label{lem:CentralizerEHIsTrivial}
						If \(X\in \lie{pgl}_k(\R)\) commutes with \(\upH_{k-1}\) then \(X\) is a diagonal matrix.  The only diagonal matrix \(X\) commuting with \(\upE_{k-1}\) is the zero matrix. As a consequence,  the centralizer of \(\upE_{k-1},\upH_{k-1}\) in \(\lie{pgl}_k(\R)\) is trivial.
					\end{lemma}
					\begin{lemma}\label{lem:fIsSimpleAndContainsPrincSL2}
						Let \(\f\subset \lie{pgl}_k(\R)\) be a reductive subalgebra that contains \(\upE_{k-1}\) and \(\upH_{k-1}\). Then
						\begin{enumerate}
							\item \(\f\) is a split simple Lie algebra, with a split Cartan subalgebra \(
							\mathfrak t_\f:=\f\cap \mathfrak t\) of \(\mathfrak f\). 
							\item\label{enu:UniqAndContainmentOFrhok-1SL2} The element \(\upS_{k-1}:=\sum_{i=1}^{k-1}i(k-i)\be_{i+1,i}\) is the unique \(S\in\f\) such that \((\upE_{k-1},\upH_{k-1},S)\) is an \(\lie{sl}_2\)-triple. In particular, \(\f\) contains \(\varrho_{k-1}\lie{sl}_2(\R):=\gen{\upE_{k-1},\upH_{k-1},\upS_{k-1}}_\R\), which acts irreducibly on \(\R^k\).
						\end{enumerate}
					\end{lemma}
					\begin{proof}
						For a real Lie algebra $\mathfrak l$, we denote its complexification by $\mathfrak l_\C:=\mathfrak l\otimes_{\R} \C$. 
						
						We first prove that \(\f\) is semisimple and split over \(\R\). 
						A reductive algebra is semi-simple if its center is trivial. Indeed, the center of \(\f_{\C}\)  is trivial  by Lemma~\ref{lem:CentralizerEHIsTrivial}. 
						
						To see that \(\f\) is split, notice by Lemma~\ref{lem:CentralizerEHIsTrivial} that \(\set{X\in\f_\C:\,[X,\upH_{k-1}]=0} = \mathfrak t_{\C}\cap \f_\C\) is a maximal abelian subalgebra in \(\f_\C\), thus diagonalizable over \(\R\). 
						In particular, 
						\(
						\mathfrak t_\f:=\f\cap \mathfrak t\) is a split
						Cartan subalgebra of \(\mathfrak f\). 
						
						We claim that there is a unique \(\upS_{k-1}\in \f\) with which \((\upE_{k-1},\upH_{k-1},\upS_{k-1})\) is a \(\lie{sl}_2\)-triple. Indeed, the element
						$$\upS_{k-1}=\sum_{i=1}^{k-1}i(k-i)\be_{i+1,i}$$
						satisfies $[\upH_{k-1},\upS_{k-1}]=-2\upS_{k-1}$ and $[\upE_{k-1},\upS_{k-1}]=\upH_{k-1}.$ 
						Let us prove the uniqueness of $\upS_{k-1}.$
						First, it follows from 
						$[\upH_{k-1},X]=-2X $ that
						$X=\sum_{i=1}^{k-1}x_i \be_{i+1,i}, $ for some $x_1,\dots,x_{k-1}\in \R.$
						Then
						$$ [\upE_{k-1},X] = \diag \bigl(x_1,x_2-x_1,\ldots,x_{k-1}-x_{k-2},-x_{k-1}\bigr)=\upH_{k-1}$$
						implies that
						$ x_i=i(k-i) $ and uniqueness follows. 
						It remains to prove that \(\upS_{k-1}\in \mathfrak f\).
						Since \(\f\) is semi-simple, we have \(\lie{pgl}_k(\R)=\f\oplus \lie m\), where \(\ad(\f)\lie m\subset \lie m\). Thus, we obtain an \(\ad\)-equivariant linear projection
						$$ \pi:\lie{pgl}_k(\R)\longrightarrow\f, \quad \pi|_{\f}=\Id. $$
						By equivariance, for $Y\in\mathfrak f,\ X\in\lie{pgl}_k(\R),$
						$ \pi([Y,X])=[Y,\pi(X)]. $
						Since \(\upE_{k-1},\upH_{k-1}\in\f\), it follows that 
						$$[\upH_{k-1},\pi(\upS_{k-1})]=-2\pi(\upS_{k-1}), \quad [\upE_{k-1},\pi(\upS_{k-1})]=\upH_{k-1}. $$
						By uniqueness,
						$\pi(\upS_{k-1})=\upS_{k-1}.$
						Therefore
						$\upS_{k-1}\in\f$.
						
						Finally, we show that \(\f\) is  simple. Assume that \(\f=\lie{a}\oplus\lie{b}\), where \(\lie{a}\) and \(\lie{b}\) are ideals such that \([\lie{a},\lie{b}]=0\). Denote \(\lie{s}:=\gen{\upE_{k-1},\upH_{k-1},S_{k-1}}_\R\) and denote by \(\lie{s}_{\lie{a}}\) and  \(\lie{s}_{\lie{b}}\) the corresponding images of \(\lie s\) in \(\lie{a}\) and \(\lie{b}\), respectively. We claim that \(\lie s\) is contained in one of the ideals \(\lie a\) or \(\lie b\); then, by Lemma~\ref{lem:CentralizerEHIsTrivial}, the complementary ideal will be trivial, which completes the proof.  
						
						Since \(\lie{s}\) is  simple, we may  assume without loss of generality that \(\lie{s}_{\lie{a}}\cong\lie{s}.\)  Let \(H_{\lie a},\, E_{\lie a}\in \lie a\) be the images of \(\upH_{k-1},\,\upE_{k-1}\).  Observe that   \[[\upE_{k-1},E_{\lie a}]=[E_{\lie a},E_{\lie a}]+[E_{\lie b},E_{\lie a}]=0.\]
						Now, the centralizer of \(\upE_{k-1}\) in \(\lie{pgl}_k(\R)\) is   \(\langle\upE_{k-1},\upE_{k-1}^2,\dots,\upE_{k-1}^{k-1}\rangle_\R\), since \(\be_k\in \R^k\) is a cyclic vector for \(\upE_{k-1}\). Thus, \(E_{\lie{a}}\in \langle\upE_{k-1},\upE_{k-1}^2,\dots,\upE_{k-1}^{k-1}\rangle_\R \), say \(E_{\lie a}=\sum_{l=1}^{k-1} a_l\upE_{k-1}^l\).
						Since \([\upH_{k-1},\upE_{k-1}]=2\upE_{k-1}\), we have \([H_{\lie a},E_{\lie{a}}]=2E_{\lie a} \). As a consequence,
						\[[\upH_{k-1},E_{\lie a}]=[H_{\lie a},E_{\lie a}]+[H_{\lie b},E_{\lie a}]=2E_{\lie a}+0=2E_{\lie a}.\]
						Also, notice that \([\upH_{k-1},\upE_{k-1}]=2\upE_{k-1}\) implies
						\([\upH_{k-1},\upE_{k-1}^l]=2l\upE_{k-1}^l\). Thus, we conclude  that
						\[2\paren{\sum_{l=1}^{k-1} a_l\upE_{k-1}^l}=[\upH_{k-1},E_{\lie a}]= \sum_{l=1}^{k-1} 2l a_l\upE_{k-1}^l.\]
						Since \(\upE_{k-1},\upE_{k-1}^2,\dots,\upE_{k-1}^{k-1}\) are linearly independent matrices, we conclude that \(E_{\lie{a}}=c\upE_{k-1}\) for some \(c\neq 0\). In view of the direct sum \(\f=\lie a\oplus\lie b\), we get that \(c=1\), so that \(\upE_{k-1}\in \lie a\). A similar argument shows that also \(\upS_{k-1}\in \lie a\). Since \([\upE_{k-1},\upS_{k-1}]=\upH_{k-1}\), the proof is complete. 
					\end{proof}
					
					\begin{proof}[Proof of Proposition~\ref{prop:AdmissibleLieAlgebra}]
						By Lemma~\ref{lem:fIsSimpleAndContainsPrincSL2}, it is left to prove that \(\f\) is\stdlist{}. 
						
						As $\f$ is simple, split, and contains   $\lie s:=\varrho_{k-1}\lie{sl}_2(\R)=\gen{\upE_{k-1},\upH_{k-1},S_{k-1}}_\R$ with both \(\lie s\) and \(\lie{pgl}_k(\R)\) \stdlist{}, suppose that \(\lie s\subsetneq \lie f\subsetneq \lie{pgl}_k(\R)\). 
						
						Since \(\lie s_{\C}\) acts irreducibly on \(\C^k\),  \(\f_{\C}\) acts irreducibly on \(\C^k\). By Dynkin's classification of irreducible simple subalgebras of the classical Lie algebras \cite[Chapter~II, \S1, Theorem~2.3 and Table~5, pp.~258 and 365--369]{Dynkin1957MaximalSubgroups},
						it follows that \(\f\) is isomorphic to one of the  Lie algebras
						\begin{equation}\label{eq:fProperAlgebra}
							\lie{sp}_{2m}(\R)\text{ if \(k=2m\)},\quad \lie{so}_{m,m+1}(\R)\text{ if \(k=2m+1\)}, \quad \g_{2(2)}\text{ if \(k=7\)}.
						\end{equation}
						
						Consider the right multiplication action of \(\lie{pgl}_k(\C)\) on the vector space of row vectors \(\C^k\).  For each of the algebras \eqref{eq:fProperAlgebra}, by
						\cite[Chapter~VIII, \S7, no.~5, Remark~2 following Proposition~11
						and Proposition~12(ii), pp.~135--136]{Bourbaki2005}, there exists \(B\in \GL_k(\R)\)  such that
						\begin{equation}\label{eq:admissible-bilinear-invariance}
							XB+BX^{\mathsf T}=0,\; X\in\f.
						\end{equation}
						Taking $X=\upH_{k-1}$ gives \( 2(k+1-i-j)B_{ij}=0,\) so $B$ is anti-diagonal. Denote $b_i=B_{i,k+1-i}$.
						The $(i,k-i)$ entry of
						$\upE_{k-1}B+B\upE_{k-1}^{\mathsf T}=0$ gives \(  b_{i+1}+b_i=0\) for \(1\leq i<k.\)
						Consequently, there exists \(c\in \R\smallsetminus \set 0\), such that
						\[
						B=cQ_k,\;\;
						(Q_k)_{ij}=(-1)^{i-1}\delta_{i+j,k+1}.
						\]
						Thus, we conclude that \(\f\) is\stdlist{} if \(k=2m\) where \(\f\cong \lie{sp}_{2m}(\R)\), or if \(k=2m+1\) and \(\f\cong\lie{so}_{m,m+1}(\R)\). 
						
						Suppose finally that $k=7$ and $\f\cong\g_{2(2)}$.
						By \cite[\S22.3, pp.~352--357]{FultonHarris},   there exists a nonzero alternating three-form \(\Omega\in \bigwedge^3(\C^7)^*\)  
						\[
						0=(D_X\Omega)(u,v,w)
						:=\Omega(uX,v,w)+\Omega(u,vX,w)+\Omega(u,v,wX),\;\forall X\in \f_{\C}.
						\]
						Let \(\dd x_i\in (\C^7)^*\) be the covector of \(\be_i\), namely \(\dd x_i(\be_j)=\delta_{i,j}\), and write \(\Omega\) in the basis \(\dd x_{ijl}:=\dd x_i\wedge \dd x_j\wedge\dd x_l\), for \(1\leq i<j<l\leq k\). 
						
						Note that  \(D_{\upH_6}dx_j=(8-2j)dx_j,\) so that \(D_{\upH_6}(\dd x_i \wedge \dd x_j \wedge \dd x_l) = (24 - 2(i+j+l))(\dd x_i \wedge \dd x_j \wedge \dd x_l).\) 
						Since \(D_{\upH_6}\Omega=0\), the basis elements \(\dd x_{ijl}\) appearing in the decomposition of \(\Omega\) are such that \(i+j+l=12\). Namely, 
						\[
						\Omega=a\,\dd x_{147}+b\,dx_{156}+c\,dx_{237}
						+d\,dx_{246}+e\,dx_{345},
						\quad a,b,c,d,e\in\mathbb C.
						\]
						
						Now, note that \( D_{\upE_6}dx_j=dx_{j-1}\) for \(1<j\leq k\) and \(D_{\upE_6}dx_1=0\),
						so that
						\[
						D_{\upE_6}\Omega
						=(a+c)dx_{137}+(a+b+d)dx_{146}
						+(c+d)dx_{236}+(d+e)dx_{245}.
						\]
						Since $D_{\upE_6}\Omega=0$, we obtain
						\(
						(a,b,c,d,e)=a(1,-2,-1,1,-1),\;a\neq 0.
						\)
						Hence, \( \Omega=a\Phi,\) with \(
						\Phi=dx_{147}-2dx_{156}-dx_{237}+dx_{246}-dx_{345}.\) So \(\g_{2(2)}\) is also\stdlist{}, which completes the proof.
					\end{proof}
					
					\bibliographystyle{alpha}
					\bibliography{references}
				\end{document}